\documentclass[letterpaper,12pt]{article}
\usepackage{amsmath,amssymb,amsthm,graphicx,xcolor,mathrsfs}
\usepackage{setspace}
\usepackage{bbm, bm}
\usepackage{url}
\usepackage[colorlinks,linkcolor=black, citecolor=blue,urlcolor=blue]{hyperref}
\usepackage[numbers]{natbib}
\usepackage{varwidth}
\usepackage[ruled, linesnumbered, lined, commentsnumbered]{algorithm2e}
\usepackage{multirow}
\usepackage{cancel}
\usepackage[normalem]{ulem}
\usepackage{subfig, float}
\usepackage{rotating}
\usepackage{array}
\allowdisplaybreaks
\usepackage{fancyvrb}
\usepackage{tikz}

\newcommand{\PreserveBackslash}[1]{\let\temp=\\#1\let\\=\temp}
\newcolumntype{C}[1]{>{\PreserveBackslash\centering}p{#1}}
\renewcommand{\d}{\mathrm{d}}

\usepackage[latin1]{inputenc}

\numberwithin{equation}{section}
\theoremstyle{plain}
\newtheorem{Theorem}{Theorem}
\numberwithin{Theorem}{section}

\newtheorem{proposition}{Proposition}
\numberwithin{proposition}{section}

\newtheorem{Lemma}{Lemma}
\numberwithin{Lemma}{section}
{
	\theoremstyle{definition}
	\newtheorem{Definition}{Definition}
	\newtheorem{example}{Example}
	\numberwithin{example}{section}
	
	\numberwithin{fact}{section}

}

\newcommand{\rev}[1]{#1}
\newcommand{\furtherrev}[1]{#1}

\newcommand{\A}{\mathcal{A}}
\newcommand{\Pro}{\mathbb{P}}
\newcommand{\1}{\mathbbm{1}}

\def\bmu{{\boldsymbol{\mu}}}

\def\bbeta{{\boldsymbol{\beta}}}
\def\bth{{\boldsymbol{\theta}}}

\def\Pr{\mathbb P}

\DeclareMathOperator*{\argmin}{arg\,min}

\def\supp{{\rm supp}}

\def\R{{\mathbb R}}

\def\Var{{\rm Var}}

\def\E{\mathbb{E}}

\def\L{{\mathcal L}}

\begin{document}
	
	\title{Score Attack: A Lower Bound Technique for Optimal Differentially Private Learning}
	\author{T. Tony Cai\footnote{Department of Statistics and Data Science, the Wharton School, University of Pennsylvania,  tcai@wharton.upenn.edu. The research of T. Tony Cai was supported in part by NSF Grant DMS-2015259 and NIH grant R01-GM129781.} ,  Yichen Wang\footnote{Independent researcher, wangyichen2012@gmail.com.}, and  Linjun Zhang\footnote{Rutgers University, linjun.zhang@rutgers.edu. The research of Linjun Zhang was supported in part by NSF Grant DMS-2015378. }}

	\date{\today}
	
	\maketitle
	
\begin{abstract}
	Achieving optimal statistical performance while ensuring the privacy of personal data is a challenging yet crucial objective in modern data analysis. However, characterizing the optimality, particularly the minimax lower bound, under privacy constraints is technically difficult.
	
	To address this issue, we propose a novel approach called the score attack, which provides a lower bound on the differential-privacy-constrained minimax risk of parameter estimation. The score attack method is based on the tracing attack concept in differential privacy and can be applied to any statistical model with a well-defined score statistic. It can optimally lower bound the minimax risk of estimating unknown model parameters, up to a logarithmic factor, while ensuring differential privacy for a range of statistical problems. We demonstrate the effectiveness and optimality of this general method in various examples, such as the generalized linear model in both classical and high-dimensional sparse settings, the Bradley-Terry-Luce model for pairwise comparisons, and nonparametric regression over the Sobolev class.
\end{abstract}
	
  
\section{Introduction}
  
With the vast amount of data being generated by individuals, businesses, and governments, statistical and machine learning algorithms are widely employed to facilitate informed decision-making in domains such as healthcare, finance, public policy, transportation, education, and academic research. The extensive use of algorithms underscores the importance of safeguarding data privacy. As a result, the differential privacy framework \cite{dwork2006our, dwork2006calibrating} for privacy-preserving data processing has garnered substantial attention. Notably, the US Census Bureau utilized differentially private methods for the first time in the 2020 US Census \cite{hawes2020implementing} to publish demographic data.

 In essence, a differentially private algorithm protects data privacy by ensuring that an observer of the algorithm's output cannot ascertain the presence or absence of any individual record in the input \furtherrev{data set}.
The design and analysis of differentially private algorithms is a rapidly evolving research field, with many differentially private solutions available in the literature for essential statistical and machine learning problems. These include mean estimation \cite{barber2014privacy, kamath2018privately, kamath2020private, cai2021cost}, top-$k$ selection \cite{bafna2017price, steinke2017tight}, linear regression \cite{wang2018revisiting, cai2021cost}, multiple testing \cite{dwork2018differentially}, causal inference \furtherrev{\cite{lee2019privacy}}, and deep learning \cite{abadi2016deep, phan2016differential}. Achieving optimal statistical performance while preserving privacy is a challenging yet crucial objective in modern data analysis.
 
While desirable for many reasons, differential privacy imposes a constraint on algorithms and may compromise their accuracy in statistical inference. In the decision-theoretical framework, the accuracy of parameter estimation is often measured by the minimax risk, which is defined as the best possible worst-case performance among \textit{all} procedures. When the class of procedures considered is limited to differentially private ones, we arrive at the \textit{privacy-constrained} minimax risk, which represents the optimal statistical performance among all differentially private methods in the worst-case scenario. 

The difference between the unconstrained minimax risk and the privacy-constrained minimax risk quantifies the cost of differential privacy, or the amount of accuracy that is inevitably lost due to differential privacy, regardless of how well the differentially private algorithm is designed. Characterizing the minimax risk under privacy constraints is technically difficult, and there have been active efforts to quantify the cost of differential privacy, in such problems as mean estimation \cite{barber2014privacy, kamath2018privately, kamath2020private, cai2021cost}, top-$k$ selection \cite{bafna2017price, steinke2017tight}, linear regression \cite{cai2021cost}, and so on.

A key step in establishing minimax theory, whether constrained or unconstrained, is the derivation of minimax lower bounds. In the classical unconstrained setting, several effective lower bound techniques have been developed in the literature, including Le Cam's two-point argument, Assouad's Lemma, and Fano's Lemma. (See \cite{le2012asymptotic, tsybakov2009introduction} for more detailed discussions on minimax lower bound arguments.) However, these methods are not directly applicable to the privacy-constrained setting, and new technical tools are needed. 

In this paper, we introduce a general technique named the ``score attack'' to establish lower bounds on the privacy-constrained minimax risk. The method is applicable to any statistical model with a well-defined score statistic, which is simply the gradient of the log-likelihood function with respect to the model parameters. After presenting the technique in general terms in Section \ref{sec: setup}, we use it to derive precise privacy-constrained minimax lower bounds across four statistical models: the low-dimensional generalized linear models (GLMs), the Bradley-Terry-Luce model for pairwise comparisons, the high-dimensional sparse GLMs, and \furtherrev{nonparametric} regression over the Sobolev class. {For the sparse GLMs, we develop a finite-difference version of the score attack. It follows paths of single support swaps and telescopes exact likelihood increments; accurate support decoding gives completeness, while differential privacy gives soundness. This makes the score attack itself capture the discrete support entropy.}

\subsection{Main Results and Our Contribution}

\textbf{The score attack technique.}
The score attack generalizes tracing-adversary arguments developed from fingerprinting codes \cite{bun2014fingerprinting,dwork2015robust} and later used for Gaussian, Beta--Binomial, sparse-mean, and regression problems \citep{steinke2017tight,kamath2018privately,cai2021cost}. Those attacks were largely tailored to individual models. Our framework instead begins with the likelihood score and reduces the privacy lower bound to a soundness--completeness comparison together with a suitable prior. Section \ref{sec: setup} makes this principle explicit, and the four ensuing applications demonstrate its effectiveness at proving concrete lower bounds.

 \textbf{Optimal differentially private rates of convergence.} In this paper, we establish minimax rates of convergence, up to stated logarithmic factors and within the regimes of the individual theorems, for parameter estimation in low-dimensional GLMs, the BTL model, high-dimensional sparse GLMs, and \furtherrev{nonparametric} regression over the Sobolev class. \rev{We design algorithms that ensure differential privacy by leveraging the Laplace and Gaussian mechanisms \cite{dwork2006calibrating} and differentially private optimization methods \cite{bassily2014private,bassily2019private,chaudhuri2011differentially,kifer2012private}.} 
The main results are summarized as follows.
\begin{itemize}
	\item Low-dimensional GLMs: Theorem \ref{thm: low-dim glm lb} presents a minimax lower bound for estimating the parameters, \rev{with the privacy term obtained directly from the GLM score attack,} and Theorem \ref{thm: non-sparse glm upper bound} shows that this lower bound is achieved, up to logarithmic factors, by a \rev{projected noisy gradient descent algorithm on the same bounded parameter space}.
		\item BTL model for pairwise comparisons: Similarly, Theorem \ref{thm: ranking lower bound} establishes a minimax lower bound for parameter estimation \rev{through an edge score attack on an identifiable product submodel,} and Theorem \ref{thm: ranking upper bound} shows that this lower bound can be attained by a \rev{pure-private projected-Laplace objective-perturbation algorithm under edge-level privacy}.
	\item \rev{High-dimensional sparse GLMs: With $H=s^*\log(ed/s^*)$, Theorem \ref{thm: high-dim glm lb} gives the lower bound of $H/n+H^2/\{n^2\varepsilon^2\log(en)\}$. Specifically, Proposition \ref{prop: stratified sparse attack} develops a support-swap score attack whose exact telescoping identity, completeness, and privacy soundness imply the desired lower bound. Theorem \ref{thm: glm upper bound} analyzes \furtherrev{private iterative hard thresholding using the NoisyHT primitive} and proves the matching polynomial rate, with the stated logarithmic gap, by a corrected contraction argument centered at the truth.}
	\item \furtherrev{Nonparametric} regression over the Sobolev class: unlike the previous problems, where the number of parameters is finite, this problem deals with estimating an entire function with a differential privacy guarantee. Here, we establish a matching lower bound in Theorem \ref{thm: nonparametric integrated risk lower bound} and an upper bound in Theorem \ref{thm: nonparametric upper bound} for the minimax mean integrated squared risk. \rev{The privacy lower bound is obtained by applying the score attack to finite-dimensional Fourier submodels and then optimizing their dimension. Their polynomial rates in $(n,\varepsilon)$ match, and the Gaussian-mechanism upper bound has the displayed logarithmic gap.}
\end{itemize}

\subsection{Related Work}

\textbf{Lower bound techniques for $(\varepsilon, \delta)$-differential privacy.} The most closely related body of work concerns fingerprinting lemmas and tracing attacks \citep{tardos2008optimal, bun2014fingerprinting, steinke2017tight, kamath2018privately}, which can be viewed as special cases of the score attack technique in Gaussian and Beta-Binomial models. More recently, \cite{kamath2022new} extended these tracing attack techniques to exponential family models. In a further refinement, \cite{narayanan2023better} improved the analysis of tracing attacks, yielding stronger lower bounds for problems such as covariance matrix estimation and heavy-tailed mean estimation.

Another related line of research \cite{barber2014privacy, karwa2017finite, acharya2018differentially, acharya2021differentially} derives lower bounds on the privacy-constrained minimax risk using differentially private analogs of classical techniques such as Le Cam's, Fano's, and Assouad's inequalities. While these analogs retain the general applicability of their classical counterparts and have produced tight lower bounds in discrete distribution estimation \cite{acharya2018differentially, acharya2021differentially}, their effectiveness in broader classes of statistical problems remains an open question.

\textbf{Differentially private algorithms for various estimation problems.} There is a substantial body of literature on differentially private generalized linear models (GLMs), with a particular focus on logistic regression \cite{chaudhuri2009privacy, chaudhuri2011differentially, zhang2020privately, song2022distributed, song2021evading, avella2021privacy, avella2021differentially}. Notably, \cite{zhang2020privately} approached sparse logistic regression under differential privacy from the perspective of graphical models. While our work is inspired by these prior studies, it differs in its primary focus on the accuracy of parameter estimation, rather than on bounding the excess risk of the learned model.

Within the broader literature on private ranking, several studies have examined differentially private rank aggregation \cite{shang2014application, hay2017differentially, song2022distributed, li2022differentially, xu2023ranking}. In this paper, we study differentially private parameter estimation in the BTL model under squared $\ell_2$ loss as one of four applications of the score-attack framework.

Regarding \furtherrev{nonparametric} function estimation under differential privacy, \cite{wasserman2010statistical} and \cite{lei2011differentially} analyzed the convergence rates of noisy histogram estimators, though without addressing optimality or lower bounds. In contrast, \cite{hall2013differential} proposed general mechanisms for releasing differentially private functional data, while \cite{barber2014privacy} developed a minimax optimal differentially private histogram estimator for Lipschitz functions.

\textbf{Statistical estimation under local differential privacy.} Local differential privacy \cite{kasiviswanathan2011can} is a related but distinct model in which each individual randomizes before data aggregation. Its minimax theory includes general frameworks \cite{duchi2013local,duchi2018minimax}, optimal mechanisms for linear functionals \cite{rohde2020geometrizing}, and results for density estimation, prediction, and covariance estimation \cite{butucea2020local,kroll2020adaptive,sart2023density,gyorfi2022rate,amorino2023minimax}.

\subsection{Organization of the Paper}

Section \ref{sec: setup} develops the score-attack framework. Sections \ref{sec: GLMs}--\ref{sec: nonparametric} treat low-dimensional GLMs, the BTL model, sparse GLMs, and \furtherrev{nonparametric} regression, respectively; Section \ref{sec: discussion} discusses extensions. All proofs are provided in the supplementary material \cite{supplement}.

\subsection{Notation}

 For real-valued sequences $\{a_n\}, \{b_n\}$, we write $a_n \lesssim b_n$ if $a_n \leq cb_n$ for some universal constant $c \in (0, \infty)$, and $a_n \gtrsim b_n$ if $a_n \geq c'b_n$ for some universal constant $c' \in (0, \infty)$.  We say $a_n \asymp b_n$ if $a_n \lesssim b_n$ and $a_n \gtrsim b_n$. \furtherrev{Symbols such as $c,C,c_0,c_1,\ldots$ denote generic positive constants only when they have not been explicitly introduced as fixed model parameters; generic constants may change from line to line. An explicitly defined symbol retains that definition within its stated scope. Each theorem specifies whether its implicit constants are universal or may depend on fixed model parameters.}

For a vector $\bm v \in \R^d$ and a subset $S \subseteq [d]$,  $\bm v_S$ denotes the \furtherrev{zero-padded (or masked) restriction} of vector $\bm v$ to the index set $S$: the $i$th coordinate of $\bm v_S$ is equal to the $i$th coordinate of $\bm v$ if $i \in S$, and zero otherwise. Define $\supp(\bm v) := \{j \in [d]: v_j \neq 0\}$. For example, for 4-dimensional vectors $\bm a = (1, 0, 0, 1)$ and $\bm b = (-1, 0, 1, 2)$, we have $\supp(\bm a) = \{1, 4\}$, $\supp(\bm b) = \{1, 3, 4\}$, and $\bm b_{\supp(\bm a)} = (-1, 0, 0, 2)$. 

We use $\|\bm v\|_p$ to denote the vector $\ell_p$ norm for $ 1\leq p \leq \infty$, with an additional convention that $\|\bm v\|_0$ denotes the number of non-zero coordinates of $\bm v$. For a square matrix $\bm A$, $\lambda_j(\bm A)$ refers to its $j$th smallest eigenvalue, and $\lambda_{\max}(\bm A), \lambda_{\min}(\bm A)$ refer to its largest and smallest eigenvalues respectively. For a function $f: \R \to \R$, $\|f\|_\infty$ denotes the essential supremum of $|f|$. For $t \in \R$ and $R > 0$, let $\Pi_R(t)$ denote the projection of $t$ onto the closed interval $[-R, R]$. 

The symbol $n$ denotes sample size in Sections \ref{sec: GLMs}, \ref{sec: sparse GLMs}, and \ref{sec: nonparametric}, and the number of items in Section \ref{sec: BTL}. \furtherrev{Unless stated otherwise, $d,p,s^*,\varepsilon,$ and $\delta$ may depend on $n$; named structural constants, including Sobolev smoothness and radius parameters, are fixed.}

 
\section{The Score Attack}\label{sec: setup}
 
We first define the privacy-constrained minimax risk and then present the score attack before turning to the four applications.

\subsection{Differential Privacy and the Minimax Risk}\label{sec: privacy constrained minimax risk}

Differential privacy limits how much the output can reveal about any one record: for data sets $\bm X$ and $\bm X'$ that differ in a single datum (``adjacent data sets''), a \rev{differentially private} algorithm requires the distributions of $M(\bm X)$ and $M(\bm X')$ to be close. \furtherrev{Throughout the paper, unless a section specifies another protected unit, adjacency means replacement adjacency: the two data sets have the same size, agree in all but one record, and the remaining record in one data set replaces that in the other. Section \ref{sec: BTL} states its specialized edge-level convention explicitly.}
  
\begin{Definition}[Differential Privacy \cite{dwork2006calibrating}]\label{def: differential privacy}
	A randomized algorithm $M: \mathcal X^n \to \mathcal R$ is $(\varepsilon, \delta)$-differentially private if for every pair of \furtherrev{replacement-adjacent} data sets $\bm X, \bm X' \in \mathcal X^n$ and every measurable $S \subseteq \mathcal R$, 
	\begin{align*}
		\Pro\left(M(\bm X) \in S\right) \leq e^\varepsilon \cdot \Pro\left(M(\bm X') \in S\right) + \delta,
	\end{align*}
	where the probability measure $\Pro$ is induced by the randomness of $M$ only.
\end{Definition}

If an algorithm is $(\varepsilon, \delta)$-differentially private for small values of $\varepsilon, \delta \geq 0$, the distributions of $M(\bm X)$ and $M(\bm X')$ are almost indistinguishable. 
The popularity of differential privacy in applications partially lies in the ease of constructing differentially private algorithms. For example, adding random noise often suffices to achieve differential privacy for many non-private algorithms.
\begin{example}[The Laplace and Gaussian Mechanisms \cite{dwork2006calibrating, dwork2014algorithmic}] \label{fc: laplace and gaussian mechanisms}
	Let $M: \mathcal X^n \to \R^d$ be an algorithm that is not necessarily differentially private.
	\begin{itemize}
			\item Suppose $\sup_{\bm X, \bm X' \text{adjacent}} \|M(\bm X) - M(\bm X')\|_1 \leq B < \infty$. For $\bm w \in \R^d$ with coordinates $w_1,\ldots,w_d \stackrel{\text{i.i.d.}}{\sim}$ Laplace$(B/\varepsilon)$, $M(\bm X) + \bm w$ is $(\varepsilon, 0)$-differentially private. Indeed, write $\bm Y=M(\bm X)+\bm w$ and $\bm Y'=M(\bm X')+\bm w$. If $f_{\bm w}$ denotes the density of $\bm w$, then
		\begin{align*}
				\frac{f_{\bm Y}(t)}{f_{\bm Y'}(t)} = \frac{f_{\bm w}(t - M(\bm X))}{f_{\bm w}(t - M(\bm X'))} \leq \exp\left(\frac{\varepsilon \|M(\bm X) - M(\bm X')\|_1}{B}\right) \leq e^\varepsilon.
		\end{align*}
			\item If instead $\sup_{\bm X, \bm X' \text{ adjacent}} \|M(\bm X) - M(\bm X')\|_2 \leq B < \infty$, then, for $0<\varepsilon<1$, $0<\delta<1$, and $\bm w \sim N_d(\bm 0, \sigma^2\bm I)$ with $\sigma^2 = 2B^2\log(2/\delta)/\varepsilon^2$, the Gaussian mechanism $M(\bm X) + \bm w$ is $(\varepsilon, \delta)$-differentially private; see, for example, \cite{dwork2014algorithmic}.
	\end{itemize}
		That is, if a non-private algorithm's output is not too sensitive to changing any single datum in the input data set, perturbing the output with Laplace or Gaussian noise produces a differentially private algorithm. 
\end{example}

Differential privacy is a desirable property, but it is also a constraint that may come at the expense of statistical accuracy. The privacy-constrained minimax risk provides a natural way to quantify this cost. Its definition consists of the following elements.
\begin{itemize}
	\item $\{f_\bth: \bth \in \Theta\}$ is a family of statistical models supported on $\mathcal X$, and \furtherrev{$\mathcal D$ is an ambient decision space containing $\Theta$}.
	\item $\bm X = \{\bm x_1, \bm x_2, \cdots, \bm x_n\}$ is an i.i.d. sample drawn from $f_{\bth^*}$ for some unknown $\bth^* \in \Theta$, and \furtherrev{$M: \mathcal X^n \to \mathcal D$} is an estimator of $\bth^*$.
	\item \furtherrev{$\ell: \mathcal D \times \Theta \to \R_{+}$ is a loss-compatible distance and $\varphi: \R_{+} \to \R_{+}$ is increasing.}
\end{itemize} 
Then, the (statistical) risk of $M$ is given by \furtherrev{$\E\varphi(\ell(M(\bm X), \bth^*))$}, where the expectation is taken over the data distribution $f_{\bth^*}$ and the randomness of estimator $M$. Because the risk \furtherrev{$\E\varphi(\ell(M(\bm X), \bth^*))$} depends on the unknown $\bth^*$ and can be minimized by choosing $M(\bm X) \equiv \bth^*$, a more sensible measure of performance is the maximum risk over the entire class of distributions $\{f_\bth: \bth \in \Theta\}$, \furtherrev{$\sup_{\bth \in \Theta} \E\varphi(\ell(M(\bm X), \bth))$}.
The minimax risk of estimating $\bth \in \Theta$ is then given by
\begin{align} \label{eq: unconstrained minimax risk}
	\inf_{\furtherrev{M:\mathcal X^n\to\mathcal D}} \sup_{\bth \in \Theta} \E\furtherrev{\varphi}(\ell(M(\bm X), \bth)),
\end{align}
where the outermost infimum is taken over the class of all \furtherrev{$\mathcal D$-valued} estimators of $\bth$. By definition, this quantity characterizes the best possible worst-case performance that an estimator can hope to achieve over the class of models $\{f_\bth: \bth \in \Theta\}$.

In this paper, we study a \textit{privacy-constrained} minimax risk: \furtherrev{let $\mathcal M_{\varepsilon,\delta}(\mathcal D)$ be the collection of all $(\varepsilon,\delta)$-differentially private algorithms from $\mathcal X^n$ to $\mathcal D$}. We consider
\begin{align}\label{eq: privacy-constrained minimax risk}
	\inf_{M \in \furtherrev{\mathcal M_{\varepsilon, \delta}(\mathcal D)}} \sup_{\bth \in \Theta} \E\furtherrev{\varphi}(\ell(M(\bm X), \bth)).
\end{align}
\furtherrev{When the ambient decision space is clear from the application, we abbreviate $\mathcal M_{\varepsilon,\delta}(\mathcal D)$ as $\mathcal M_{\varepsilon,\delta}$.} Because $\mathcal M_{\varepsilon,\delta}(\mathcal D)$ is a subset of all \furtherrev{$\mathcal D$-valued} estimators, the privacy-constrained risk \eqref{eq: privacy-constrained minimax risk} is at least the unconstrained risk \eqref{eq: unconstrained minimax risk}. Their difference is the ``cost of privacy.''

Both risks are typically characterized by matching upper and lower bounds. Analyzing a concrete algorithm gives an upper bound, whereas a lower bound must apply simultaneously to all estimators. The score attack supplies this second direction under a privacy constraint.

\subsection{The Score Attack}\label{sec: general lb}

The score attack is a type of tracing attack \cite{bun2014fingerprinting, dwork2015robust, dwork2017exposed}. A tracing attack is an algorithm which takes a single ``candidate'' datum as input and attempts to infer whether this candidate belongs to a given data set or not, by comparing the candidate with some summary statistics computed from the data set. Statisticians may envision a tracing attack as a hypothesis test which rejects the null hypothesis that the candidate is out of the data set when some test statistic takes a large value. This hypothesis testing formulation motivates some desirable properties for a tracing attack.
\begin{itemize}
	\item Soundness (type I error control): if the candidate does not belong to the data set, the tracing attack is likely to \rev{take} small values. 
	\item Completeness (type II error control): if the candidate does belong, the tracing attack is likely to take large values.
\end{itemize}
For example, \cite{dwork2015robust, kamath2018privately, cai2021cost} show that, if the random sample $\bm X$ and the candidate $\bm z$ are drawn from a Gaussian distribution with mean $\bmu$, tracing attacks of the form $\langle M(\bm X) - \bmu, \bm z - \bmu \rangle$ \rev{are} sound and complete provided that $M(\bm X)$ is an accurate estimator of $\bmu$. 

It is this accuracy requirement that connects tracing attacks with risk lower bounds for differentially private algorithms: if an estimator $M(\bm X)$ is differentially private, it cannot possibly be too close to the estimand, or the existence of tracing attacks leads to a contradiction with the guarantees of differential privacy. Designing sound and complete tracing attacks, therefore, is crucial to the sharpness of privacy-constrained minimax lower bounds. Besides the Gaussian mean tracing attack mentioned above, there are some successful tracing attacks proposed for specific problems, such as top-$k$ selection \cite{steinke2017tight} or linear regression \cite{cai2021cost}, but a general recipe for the design and analysis of tracing attacks has not been available. 

The score attack is a form of tracing attack applicable to general parametric families of distributions. Given a parametric family $\{f_\bth(\bm x): \bth \in \Theta\}$ with $\Theta \subseteq \R^d$, the score statistic, or simply the score, is $S_\bth(\bm x) := \nabla_\bth \log f_\bth(\bm x)$. Suppose the densities have parameter-independent support, are differentiable in $\bth$, and are locally dominated so that differentiation can be passed under the integral. Then $\E_\bth S_\bth(\bm x)=\bm0$; if the second score moment is finite, $\mathcal I(\bth):=\E_\bth\{S_\bth(\bm x)S_\bth(\bm x)^\top\}=\Var_\bth\{S_\bth(\bm x)\}$ is the Fisher information matrix. These are the standard score regularity conditions; see, for example, \cite{lehmann2006theory}. \furtherrev{In this smooth parametric subsection, the decision space is $\mathcal D=\mathbb R^d$ and $M$ is $\mathbb R^d$-valued. If an application begins with a larger decision space, an explicitly stated post-processing map first produces the vector estimator to which the attack is applied.} The score attack is defined as
\begin{align}\label{eq: score attack}
	\A_\bth(\bm z, M(\bm X)) := \langle M(\bm X) - \bth, S_\bth(\bm z)   \rangle.
\end{align}
Large positive values of $\A_\bth(\bm z, M(\bm X))$ provide evidence that $\bm z$ belongs to $\bm X$. In particular, if $f_\bth(\bm x)$ is the density of $N(\bth, \bm I)$, the score attack coincides with the tracing attacks for Gaussian means studied in \citep{dwork2015robust, kamath2018privately, cai2021cost}.

{There is also an exact finite-difference version of the score identity. For two laws $P,Q$ with $Q\ll P$, define the likelihood-increment score $T_{P\to Q}(z)=dQ/dP(z)-1$. Then, for every $h$ that is integrable under both $P$ and $Q$,
\begin{align*}
\E_P\{h(Z)T_{P\to Q}(Z)\}=\E_Qh(Z)-\E_Ph(Z).
\end{align*}
Thus a chain of nearby parameter values produces a telescoping score attack without differentiating the parameter. Section \ref{sec: sparse GLMs} applies this identity along single support swaps. There, accurate decoding supplies completeness, differential privacy supplies soundness, and the disjoint decoding cells retain the full support entropy.}

The next theorem records the two properties needed later: soundness when the candidate is an independent replacement, and completeness when it is one of the observations.
\begin{Theorem}\label{thm: score attack general}
		Let $\bm X = \{\bm x_1, \bm x_2, \cdots, \bm x_n\}$ be an i.i.d. sample drawn from $f_\bth$. For each $i \in [n]$, let $\bm X'_{i}$ denote an adjacent data set of $\bm X$ obtained by replacing $\bm x_i$ with an independent copy $\bm x'_i \sim f_\bth$. Assume $\E_\bth S_\bth(\bm x)=\bm0$, $\mathcal I(\bth)=\E_\bth\{S_\bth(\bm x)S_\bth(\bm x)^\top\}$ is finite, and $\E_\bth\|M(\bm X)-\bth\|_2^2<\infty$.
	\begin{enumerate}
		\item Soundness: for each $i \in [n]$, 
		\begin{align}\label{eq: soundness general}
			\E \A_\bth(\bm x_i, M(\bm X'_i)) = 0; ~ \E |\A_\bth(\bm x_i, M(\bm X'_i))| \leq \sqrt{\E\|M(\bm X) - \bth\|_2^2}\sqrt{\lambda_{\max}(\mathcal I(\bth))}.
		\end{align}
		\item Completeness: write $f_\bth^{(n)}$ for the joint density of $\bm X$ with respect to a common dominating measure, and assume the family has common support. Suppose $\E_\bth|M_j(\bm X)|<\infty$ and, for every $j\in[d]$ and every $\bth$, $f_\bth^{(n)}(\bm x)$ is differentiable in $\theta_j$ and there are a neighborhood $U$ of $\bth$ and an integrable function $G_j$ such that
		\[
		\left|\E_M\{M_j(\bm x)\}\,\partial_{\theta_j}f_{\boldsymbol\vartheta}^{(n)}(\bm x)\right|\leq G_j(\bm x)
		\quad\text{for every }\boldsymbol\vartheta\in U.
		\]
		Here $\E_M$ denotes expectation over the internal randomness of $M$ for a fixed input. Then
		\begin{align}\label{eq: completeness general}
			\sum_{i \in [n]} \E \A_\bth (\bm x_i, M(\bm X)) = \sum_{j \in [d]} \frac{\partial}{\partial \theta_j} \E M(\bm X)_j.
		\end{align}
	\end{enumerate} 
\end{Theorem}
Theorem \ref{thm: score attack general} is proved in Section \ref{sec: proof of thm: score attack general} of the supplement. The special form of completeness for Gaussian and Beta--Binomial families is known as a fingerprinting lemma \citep{tardos2008optimal, bun2014fingerprinting, steinke2017tight, kamath2018privately}. To convert the theorem into a risk lower bound, we compare the attack evaluated on a sample point with the same attack after that point has been replaced independently. Differential privacy yields an inequality of the form
\[
\E \A_\bth (\bm x_i, M(\bm X))
\leq \E \A_\bth (\bm x_i, M(\bm X'_i))
+O(\varepsilon)\E |\A_\bth (\bm x_i, M(\bm X'_i))|,
\]
up to the explicit $\delta$ and tail terms stated in Proposition \ref{prop: score attack upper bound}. Combined with soundness \eqref{eq: soundness general}, this gives
\begin{align*}
	\sum_{i \in [n]} \E \A_\bth (\bm x_i, M(\bm X)) \leq \sqrt{\E\|M(\bm X) - \bth\|_2^2} \cdot n\sqrt{\lambda_{\max}(\mathcal I(\bth))} O(\varepsilon).
\end{align*}
Section \ref{sec: attack upper bound general} gives the precise privacy comparison.

A positive lower bound on $\sum_{i \in [n]} \E \A_\bth (\bm x_i, M(\bm X))$ can therefore be compared with the preceding upper bound to force a lower bound on the risk. Completeness supplies this second direction: averaging $\partial\E M_j/\partial\theta_j$ under a suitable prior and applying Stein's lemma relates the derivative to the estimation error. Section \ref{sec: attack lower bound general} formalizes this step \cite{stein1972bound, stein2004use}; related arguments for Gaussian mean estimation and top-$k$ selection appear in \citep{dwork2015robust, steinke2017between}. Together, the privacy comparison and the prior-averaged completeness bound form the lower-bound template used in the examples below.

\subsubsection{Score Attack and Differential Privacy}\label{sec: attack upper bound general}
When $\bm X'_i$ contains an independent replacement rather than the candidate $\bm x_i$, Theorem \ref{thm: score attack general} gives
\begin{align*}
	\E \A_\bth(\bm x_i, M(\bm X'_i)) = 0; ~ \E |\A_\bth(\bm x_i, M(\bm X'_i))| \leq \sqrt{\E\|M(\bm X) - \bth\|_2^2}\sqrt{\lambda_{\max}(\mathcal I(\bth))}.
\end{align*}
If $M$ is differentially private, the distribution of $M(\bm X'_i)$ is close to that of $M(\bm X)$; as a result, the inequalities above can be related to the case where the data set $\bm X$ does include the candidate $\bm x_i$.
\begin{proposition}\label{prop: score attack upper bound}
			Assume $\E_\bth S_\bth(\bm x)=\bm0$ and $\mathcal I(\bth)$ is finite. Suppose also that
			\begin{align*}
				\E_\bth\|M(\bm X)-\bth\|_2^2<\infty,
				\qquad
				\E_\bth|\A_\bth(\bm x_i,M(\bm X))|<\infty.
			\end{align*}
			If $M$ is an $(\varepsilon, \delta)$-differentially private algorithm with $0 < \varepsilon < 1$ and $0\leq\delta<1$, then for every $T > 0$, 
	\begin{align}\label{eq: score attack upper bound}
		\E \A_\bth(\bm x_i, M(\bm X))
		&\leq 2\varepsilon \sqrt{\E\|M(\bm X) - \bth\|_2^2}
		\sqrt{\lambda_{\max}(\mathcal I(\bth))} + 2\delta T \notag\\
			&\quad + \int_T^\infty \left\{\Pro\left(|\A_\bth(\bm x_i, M(\bm X))| > t \right)
			+\Pro\left(|\A_\bth(\bm x_i, M(\bm X_i'))| > t \right)\right\} \d t.
	\end{align}
\end{proposition}
Proposition \ref{prop: score attack upper bound} is proved in Section \ref{sec: proof of prop: score attack upper bound}. The quantity on the right side of \eqref{eq: score attack upper bound} is determined by the statistical model $f_\bth(\bm x)$ and the choice of $T$.

\subsubsection{Score Attack and Stein's Lemma}\label{sec: attack lower bound general}
Let $g(\bth)=\E_{\bm X|\bth}M(\bm X)$. Thus $g$ maps $\Theta$ to $\R^d$, and we are interested in bounding $\partial g_j(\bth)/\partial\theta_j$ from below. Stein's Lemma \cite{stein1972bound, stein2004use} is useful for this purpose.

\begin{Lemma}[Stein's Lemma]\label{lm: stein's lemma}
		Let $Z$ have a density $p$ supported on $[a,b]$, where $-\infty\leq a<b\leq\infty$. Suppose that $p$ is positive and continuously differentiable on $(a,b)$ and that $h:(a,b)\to\R$ is differentiable, $\E|h'(Z)|<\infty$, and \rev{$\E |h(Z)p'(Z)/p(Z)| < \infty$}. If the endpoint limits below exist, then
	\begin{align}\label{eq: stein's lemma}
			\E h'(Z) = \E\left[\frac{-h(Z)p'(Z)}{p(Z)}\right]
			+\lim_{t\uparrow b}h(t)p(t)-\lim_{t\downarrow a}h(t)p(t),
	\end{align}
		where the limits are understood along the support, including at infinite endpoints. In particular, if $p(z) = (2\pi)^{-1/2}e^{-z^2/2}$ and the boundary products vanish, then $\E h'(Z) = \E Z h(Z)$.
\end{Lemma}
Stein's Lemma implies that, by imposing appropriate prior distributions on $\bth$, one can obtain a lower bound for $\frac{\partial}{\partial\theta_j} g_j(\bth)$ on average over the prior distribution of $\bth$, as follows.
\begin{proposition}\label{prop: score attack stein's lemma}
		Let $\bth$ be distributed according to a product density $\bm \pi=\prod_{j\in[d]}\pi_j$. \rev{For each $j\in[d]$, suppose that, for $\pi_{-j}$-almost every $\bth_{-j}$, the one-dimensional map $t\mapsto g_j(t,\bth_{-j})$ and the density $\pi_j$ satisfy the regularity conditions in Lemma \ref{lm: stein's lemma}. Suppose also that the resulting integrands are jointly integrable, so that Fubini's theorem applies, and that
		\[
		\lim_{t\uparrow b_j}g_j(t,\bth_{-j})\pi_j(t)
		=\lim_{t\downarrow a_j}g_j(t,\bth_{-j})\pi_j(t)=0,
		\]
		where $[a_j,b_j]$ is the support of $\pi_j$.} Then
	\begin{align}\label{eq: score attack stein's lemma}
		\E_{\bm \pi} \left(\sum_{j \in [d]} \frac{\partial}{\partial \theta_j} g_j(\bth)\right)
		&\geq \E_{\bm \pi} \left(\sum_{j \in [d]} \frac{-\theta_j \pi'_j(\theta_j)}{\pi_j(\theta_j)}\right) \notag\\
		&\quad - \sqrt{\E_{\bm \pi}\E_{\bm X|\bth}\|M(\bm X) - \bth\|_2^2}
		\sqrt{\E_{\bm \pi} \left[\sum_{j \in [d]}\left(\frac{\pi'_j(\theta_j)}{\pi_j(\theta_j)}\right)^2\right]}.
	\end{align}
\end{proposition}
Proposition \ref{prop: score attack stein's lemma} is proved in Section \ref{sec: proof of prop: score attack stein's lemma}. In addition to the standard regularity conditions of Stein's Lemma, it assumes directly that the products $g_j\pi_j$ vanish at the endpoints. This condition holds, for example, when $g$ is bounded and the marginal prior densities vanish at their endpoints. If these products do not vanish, the two endpoint limits in \eqref{eq: stein's lemma} must be retained. \rev{The sparse-GLM construction in Section \ref{sec: sparse glm lower bound} instead uses the finite-difference score identity above, so it neither differentiates through support selection nor requires a product prior across coordinates.} The endpoint terms are retained explicitly when needed.

The right side of \eqref{eq: score attack stein's lemma} is especially useful
for minimax analysis.  After restricting attention to estimators satisfying
$\sup_{\bth \in \Theta}\E_{\bm X|\bth}
\|M(\bm X)-\bth\|_2^2<\furtherrev{R_0}$, the remaining quantities depend only on the chosen
prior $\pi$.
\begin{example}\label{ex: Gaussian stein's lemma}
		Let $\Theta=\R^d$ and let $\bm \pi$ be the density of $N(\bm 0, \bm I)$. For every estimator $M$ satisfying $\sup_{\bth \in \R^d}\E_{\bm X|\bth} \|M(\bm X) - \bth\|_2^2 < \furtherrev{R_0}$, \eqref{eq: score attack stein's lemma} reduces to
	\begin{align*}
			\E_{\bm \pi} \left(\sum_{j \in [d]} \frac{\partial}{\partial \theta_j} g_j(\bth)\right) \geq \sum_{j \in [d]} \E_{\pi_j} \theta_j^2 - \sqrt{\furtherrev{R_0}} \sqrt{\sum_{j \in [d]} \E_{\pi_j} \theta_j^2} = d - \sqrt{\furtherrev{R_0} d} \geq \frac d2,
		\end{align*}
		provided $d\geq\furtherrev{4R_0}$.
\end{example}

In view of the completeness property \eqref{eq: completeness general}, Proposition \ref{prop: score attack stein's lemma} suggests an \textit{average} lower bound for $\sum_{i \in [n]} \E \A_\bth (\bm x_i, M(\bm X))$ over a prior distribution $\bm \pi(\bth)$, with the form of the bound determined by the choice of $\bm \pi$. This parallels the familiar fact that the Bayes risk lower bounds the minimax risk, which is the final step in our minimax lower-bound arguments.

\subsubsection{From Score Attack to Lower Bounds}\label{sec: attack to lower bounds}
Theorem \ref{thm: score attack general} combined with Propositions \ref{prop: score attack upper bound} and \ref{prop: score attack stein's lemma} reveals the connection between the score attack and privacy-constrained minimax lower bounds.

Let $\bm \pi=\prod_{j=1}^d\pi_j$ be a product prior supported on the parameter space $\Theta$, and \rev{consider an estimator satisfying} $\E_{\bm X|\bth} \|M(\bm X) - \bth\|_2^2 < \furtherrev{R_0}$ for every $\bth \in \Theta$. \rev{In each application below, estimators outside this risk regime are handled separately; when projection is used, its validity and any constant-factor loss are stated explicitly.} The completeness part of Theorem \ref{thm: score attack general} and Proposition \ref{prop: score attack stein's lemma} imply that
\begin{align*}
	\sum_{i \in [n]} \E_{\bm \pi}\E_{\bm X|\bth} \A_\bth (\bm x_i, M(\bm X)) \geq \E_{\bm \pi} \left(\sum_{j \in [d]} \frac{-\theta_j \pi'_j(\theta_j)}{\pi_j(\theta_j)}\right) - \sqrt{\furtherrev{R_0}}\sqrt{  \E_{\bm \pi} \left[\sum_{j \in [d]}\left(\frac{\pi'_j(\theta_j)}{\pi_j(\theta_j)}\right)^2\right]}.
\end{align*}
Since Proposition \ref{prop: score attack upper bound} holds for every $\bth$, it follows that
	\begin{align*}
		&\sum_{i \in [n]} \E_{\bm \pi}\E_{\bm X|\bth} \A_\bth (\bm x_i, M(\bm X)) \\
			&\leq 2n\varepsilon \sqrt{\E_{\bm \pi} \E_{\bm X|\bth}\|M(\bm X) - \bth\|_2^2}
			\sqrt{\sup_{\bth\in\Theta}\lambda_{\max}(\mathcal I(\bth))} + 2n\delta T \\
			&\quad + \sum_{i \in [n]}\E_{\bm\pi}\int_T^\infty
			\left\{\Pro_{\bm X|\bth}\left(|\A_\bth(\bm x_i, M(\bm X))| > t \right)
			+\Pro_{\bm X|\bth}\left(|\A_\bth(\bm x_i, M(\bm X_i'))| > t \right)\right\}\d t.
\end{align*}
Comparing the two inequalities yields a lower bound on the Bayes risk
\begin{align*}
\inf_{M \in \mathcal M_{\varepsilon, \delta}} \E_{\bm \pi}
\E_{\bm X|\bth}\|M(\bm X)-\bth\|_2^2.
\end{align*}
Because the worst-case risk dominates the Bayes risk under every prior
supported on $\Theta$, the same argument lower bounds the privacy-constrained
minimax risk in \eqref{eq: privacy-constrained minimax risk}.

\subsection{The Utility of Score Attack}
The analysis in Section \ref{sec: general lb} reduces the privacy-constrained minimax lower bound \eqref{eq: privacy-constrained minimax risk} to analyzing the expectation of the score attack,
\begin{align*}
	\sum_{i \in [n]} \E_{\bm X|\bth} \A_\bth (\bm x_i, M(\bm X)).
\end{align*}
Specifically, the analysis of score attack consists of upper bounding the expectation via differential privacy, and lower bounding the expectation ``on average'' by choosing a prior over the parameter space $\Theta$.

\rev{This comparison is the privacy-specific engine used throughout the paper. In smooth parameter directions, the usual likelihood score and Stein's identity produce the obstruction to accurate private estimation. For a discrete parameter direction, the same architecture can instead be built from finite likelihood increments: a telescoping change-of-measure identity supplies the attack value, while decoding and privacy supply completeness and soundness.}

The remainder of the paper applies this template to the following problems.
\begin{itemize}
	\item Parameter estimation in classical models: the generalized linear model (Section \ref{sec: GLMs}), and the Bradley-Terry-Luce model (Section \ref{sec: BTL}).
	\item High-dimensional sparse parameter estimation (Section \ref{sec: sparse GLMs}).
	\item \furtherrev{Nonparametric} function estimation (Section \ref{sec: nonparametric}).
\end{itemize}
\rev{In each example, the score attack provides the model-specific privacy lower-bound calculation, and its sharpness is demonstrated by a concrete differentially private algorithm. The low-dimensional GLM, BTL, and \furtherrev{nonparametric} results use smooth scores. In the sparse example, Proposition \ref{prop: stratified sparse attack} uses a finite-difference score along support swaps and directly produces the support-entropy privacy term in Theorem \ref{thm: high-dim glm lb}.}


\section{The Generalized Linear Model}\label{sec: GLMs}
Generalized linear models (GLMs) are widely used in modern data-driven scientific research, with applications spanning genetics, metabolomics, finance, and econometrics. They also play a central role in many observational studies, where privacy concerns are often paramount.

As the first application of the score attack technique, we examine the privacy-constrained minimax risk for estimating parameters $\bbeta \in \R^d$ in the generalized linear model with scale parameter $\sigma$:
\begin{align}
		f_{\bbeta} (y\mid\bm x) = h(y, \sigma)\exp\left(\frac{y\bm x^\top \bbeta - \psi(\bm x^\top \bbeta)}{c(\sigma)}\right),\qquad \bm x \sim f_{\bm x}. \label{eq: glm definition}
\end{align}
We use an i.i.d. sample $\bm Z = \{\bm z_i\}_{i \in [n]} = \{(y_i, \bm x_i)\}_{i \in [n]}$ drawn from \eqref{eq: glm definition}. The functional form of the model, including the partition function $\psi$ and the normalizing factor $h$, is fixed and known; the sole parameter of interest is $\bbeta$.

{Throughout this section, we work on the fixed compact parameter space
\[
\Theta_0=\{\bbeta\in\mathbb R^d:\|\bbeta\|_2\leq b_0\},
\]
where $b_0>0$ is a model constant.  This makes the minimax comparison in the lower and upper bounds literal and permits a data-independent initialization for the upper-bound algorithm.}
{The scale $c(\sigma)$ is also treated as a fixed model constant, bounded above and away from zero; constants below may depend on these fixed bounds.}
\furtherrev{The ambient decision space is $\mathcal D=\mathbb R^d$, so $\mathcal M_{\varepsilon,\delta}$ in this section abbreviates $\mathcal M_{\varepsilon,\delta}(\mathbb R^d)$. The projection reduction used in the lower-bound proof is stated explicitly in Section \ref{sec: proof of thm: low-dim glm lb}.}

We say that the random design satisfies
\begin{itemize}
	\item [(SG)] Sub-Gaussian upper-tail condition:
	$\lambda_{\max}\{\E(\bm x\bm x^\top)\}\leq C$ and
	$\|\bm u^\top\bm x\|_{\psi_2}\leq K_{\bm x}\|\bm u\|_2$
	for every $\bm u\in\mathbb R^d$, for fixed constants
	$C,K_{\bm x}<\infty$.
\end{itemize}
Here $\|\cdot\|_{\psi_2}$ is the sub-Gaussian Orlicz norm. Condition (SG)
places no almost-sure bound on $\|\bm x\|_2$.

In Section \ref{sec: non-sparse glm lower bound}, we establish a minimax risk lower bound by the score attack. Section \ref{sec: non-sparse glm upper bound} gives a noisy gradient descent estimator that attains the same polynomial rate up to logarithmic factors.

\subsection{The Privacy-Constrained Minimax Lower Bound} \label{sec: non-sparse glm lower bound}
For the generalized linear model \eqref{eq: glm definition} and a candidate datum $(\tilde y, \tilde{\bm x})$, the score attack \eqref{eq: score attack} takes the form
\begin{align}\label{eq: low-dim glm attack}
	\A_{\bbeta} ((\tilde y, \tilde{\bm x}), M(\bm y, \bm X)) = \frac{1}{c(\sigma)} \big\langle M(\bm y, \bm X) - \bbeta, [\tilde y - \psi'(\tilde{\bm x}^\top\bbeta)] \tilde {\bm x} \big\rangle.
\end{align}
As outlined in Section \ref{sec: general lb}, we establish a privacy-constrained minimax lower bound for estimating $\bbeta$ by analyzing the sum of expectations $\sum_{i \in [n]} \E \A_\bbeta ((y_i, \bm x_i), M(\bm y, \bm X))$. When the reference to data $(\bm y, \bm X)$ and estimator $M$ is clear, we abbreviate $\A_\bbeta ((y_i, \bm x_i), M(\bm y, \bm X))$ as $A_i$.

We begin by upper-bounding $\sum_{i \in [n]} \E A_i$, which amounts to specializing the soundness part of Theorem \ref{thm: score attack general} and Proposition \ref{prop: score attack upper bound} to the GLM score attack \eqref{eq: low-dim glm attack}. 
\begin{proposition}\label{prop: low-dim glm attack soundness}
			Consider i.i.d. observations $(y_1, \bm x_1), \cdots, (y_n, \bm x_n)$ drawn from \eqref{eq: glm definition}. Suppose (SG) holds and $\|\psi^{''}\|_\infty < c_2 < \infty$. If the estimator $M$ takes values in $\Theta_0$ and is $(\varepsilon, \delta)$-differentially private with $0 < \varepsilon < 1$ and $0\leq\delta<1/2$, then
		\begin{align}\label{eq: low-dim glm attack soundness}
			\sum_{i \in [n]} \E_{\bm y, \bm X|\bbeta} A_i
			&\leq 2n\varepsilon\sqrt{\E\|M(\bm y, \bm X) - \bbeta\|^2_2}
			\sqrt{Cc_2/c(\sigma)} \notag\\
			&\quad+C_3n\delta\sqrt{c_2/c(\sigma)}
			\left\{\sqrt{d\ell_\delta}+\ell_\delta\right\},
		\end{align}
		where $\ell_\delta=\log(e/\delta)$ when $\delta>0$, the last term is interpreted as zero when $\delta=0$, and $C_3$ depends only on the fixed constants $b_0$ and $K_{\bm x}$.
\end{proposition}
Proposition \ref{prop: low-dim glm attack soundness} follows by computing the GLM Fisher information and choosing the truncation level $T$ in \eqref{eq: score attack upper bound}; Section \ref{sec: proof of prop: low-dim glm attack soundness} gives the calculation.  The complementary lower bound on the attack is as follows.

\begin{proposition}\label{prop: low-dim glm attack completeness}
		\rev{Let $V_1,\ldots,V_d$ be i.i.d. with density $p(v)=15(1-v^2)^2/16$ on $[-1,1]$, let $a=b_0/(2\sqrt d)$, and set $\beta_j=aV_j$. Thus the product prior $\bm\pi$ is supported on $\Theta_0$. There is a constant $c_0>0$, depending only on $b_0$, such that every estimator $M$ taking values in $\Theta_0$ and satisfying $\E_{\bm\pi}\E_{\bm y,\bm X\mid\bbeta}\|M(\bm y,\bm X)-\bbeta\|_2^2\leq c_0b_0^2$ obeys}
	\begin{align}\label{eq: low-dim glm attack completeness}
		\sum_{i \in [n]} \E_{\bm \pi}\E_{\bm y, \bm X|\bbeta} A_i \gtrsim d, 
	\end{align}
\rev{where $\bm \pi$ denotes the scaled, centered Beta product prior just defined.}
\end{proposition}
The proof of Proposition \ref{prop: low-dim glm attack completeness} applies Proposition \ref{prop: score attack stein's lemma} to the displayed product prior and is given in Section \ref{sec: proof of prop: low-dim glm attack completeness}.  Comparing the soundness and completeness bounds yields the following minimax result.

\begin{Theorem}\label{thm: low-dim glm lb}
		Consider i.i.d. observations $(y_1, \bm x_1), \cdots, (y_n, \bm x_n)$ drawn from \eqref{eq: glm definition}. Suppose (SG) holds and $\|\psi^{''}\|_\infty < c_2 < \infty$. If $d \lesssim n\varepsilon$, $0 < \varepsilon < 1$ and $\delta \lesssim n^{-(1+\gamma)}$ for some $\gamma > 0$, then
	\begin{align}\label{eq: low-dim glm lb}
		\rev{\inf_{M \in \mathcal M_{\varepsilon, \delta}}\sup_{\bbeta \in \Theta_0}} \E\|M(\bm y, \bm X)- \bbeta\|_2^2 \gtrsim c(\sigma) \left(\frac{d}{n} + \frac{d^2}{n^2\varepsilon^2}\right).
	\end{align}
	\furtherrev{All implicit constants may depend only on $\gamma$ and the fixed model constants $b_0,C,K_{\bm x},c_2$, and the fixed bounds on $c(\sigma)$, but not on $n,d,\varepsilon$, or $\delta$.}
\end{Theorem}
The first term in \eqref{eq: low-dim glm lb} is the non-private minimax risk lower bound, and the second term is the ``cost of differential privacy.'' \rev{The latter is obtained entirely from the soundness--completeness comparison for the GLM score attack; the separate Assouad argument in the supplement supplies only the classical $d/n$ term.} We show in the next section that the lower bound is attainable, up to a logarithmic term, by a noisy gradient descent algorithm.

The condition $d \lesssim n\varepsilon$ in Theorem \ref{thm: low-dim glm lb} restricts the result to a low-dimensional regime. Section \ref{sec: sparse GLMs} treats a separate high-dimensional sparse regime with $d\geq n$.

\subsection{Optimality of the Private GLM Lower Bound}\label{sec: non-sparse glm upper bound}
We consider minimizing the negative GLM log-likelihood
\begin{align*}
	\L_n(\bbeta; \bm Z) = \frac{1}{n}\sum_{i=1}^n \left(\psi(\bm x_i^\top \bbeta) - y_i\bm x_i^\top \bbeta\right)
\end{align*}
by a noisy gradient descent algorithm, first proposed by \cite{bassily2014private} in its generic form for arbitrary convex functions. The following algorithm specializes the generic algorithm to GLMs.

For $L>0$, define the Euclidean clipping map
\begin{align*}
	\mathsf C_L^{(2)}(\bm x)=\bm x\min\{1,L/\|\bm x\|_2\}.
\end{align*}
We also use a smooth saturation of the link score.  Fix an even
$C^1$ function $\chi:\mathbb R\to[0,1]$ such that
$\chi(t)=1$ for $|t|\leq1$, $\chi(t)=0$ for $|t|\geq2$, and
$\|\chi'\|_\infty<\infty$.  For $L_\psi\geq1$, define
\begin{align}
	\psi_{L_\psi}'(t)
	&=\psi'(0)+\int_0^t\psi''(u)\chi(u/L_\psi)\,\d u,
	&
	\psi_{L_\psi}(t)
	&=\psi(0)+\int_0^t\psi_{L_\psi}'(u)\,\d u.       \label{eq: saturated link score}
\end{align}
Under (G) below, $\psi_{L_\psi}$ is convex,
$\psi_{L_\psi}'(t)=\psi'(t)$ for $|t|\leq L_\psi$, and
\begin{align}
	\|\psi_{L_\psi}'\|_\infty
	\leq S_\psi:=|\psi'(0)|+2c_2L_\psi.              \label{eq: saturated score bound}
\end{align}
Both clipping levels below are deterministic. Thus the privacy guarantee
holds for every data set, not merely on a high-probability event under
the sampling model.

\begin{algorithm}[H]\label{algo: private glm}
	\SetAlgoLined
	\SetKwInOut{Input}{Input}
	\SetKwInOut{Output}{Output}
	\Input{Data set $\bm Z$, \rev{parameter set $\Theta_0$}, step size $\eta^0$, privacy parameters $\varepsilon, \delta$, noise scale $B$, number of iterations $T$, response-truncation level $R$, design-clipping radius $L_{\bm x}$, link-score saturation level $L_\psi$, initial value $\bbeta^0 \in \rev{\Theta_0}$.}
	Set $\overline{\bm x}_i=\mathsf C_{L_{\bm x}}^{(2)}(\bm x_i)$ for $i\in[n]$\;
	\For{$t$ in $0$ \KwTo $T-1$}{
		Generate $\bm w_t \in \R^d$ with $w_{t1}, w_{t2}, \cdots, w_{td} \stackrel{\text{i.i.d.}}{\sim} N\left(0, (\eta^0)^2 2B^2\frac{\log(2T/\delta)}{n^2(\varepsilon/T)^2}\right)$\;
		Compute $\bbeta^{t + 1} = {\Pi_{\Theta_0}\!\left[\bbeta^t - (\eta^0/n)\sum_{i=1}^n  (\psi_{L_\psi}'(\overline{\bm x}_i^\top \bbeta^t)-\Pi_{R}(y_i))\overline{\bm x}_i + \bm w_t\right]}$\;
	}
	\Output{$\bbeta^T$.}
	\caption{Differentially Private Generalized Linear Regression}
\end{algorithm}
The privacy and accuracy analysis uses the following assumptions.
\begin{itemize}
	\item [(D)] Well-conditioned sub-Gaussian design: (SG) holds, $\E\bm x=\bm0$, and the covariance matrix $\Sigma_{\bm x}=\E\bm x\bm x^\top$ satisfies $\Sigma_{\bm x}\succeq C^{-1}I_d$.
	\item [(G)] Curvature and smoothness: $\psi\in C^3(\mathbb R)$, $0<\psi''(t)\leq c_2$ for every $t\in\mathbb R$, and $\|\psi'''\|_\infty\leq c_3$, for fixed constants $c_2,c_3<\infty$.
\end{itemize}
\rev{Sub-Gaussian random design, well-conditioned covariance, and smooth positive conditional variance are standard finite-sample GLM conditions; see, for example, \cite{bach2010self,mei2018landscape,ostrovskii2021finite}.  We impose no almost-sure bound on the design, no dimension-free bound on $|\bm x^\top\bbeta|$, and no global bound on $\psi'$.  Deterministic design clipping and link-score saturation supply worst-case privacy sensitivity, while (D) supplies the high-probability design, predictor, and empirical-curvature events used for accuracy.}

\begin{proposition}\label{prop: low-dim glm empirical curvature}
	Under (D) and (G), there are constants $0<\alpha\leq\gamma<\infty$, $C_0<\infty$, and $L_{\psi,0}<\infty$, depending only on the fixed constants in these assumptions and on $b_0$, such that, for every deterministic $L_\psi\geq L_{\psi,0}$, whenever $n\geq C_0d\log(en)$,
	\begin{align}\label{eq: low-dim glm empirical curvature}
	\Pro\!\left(\alpha I_d\preceq \frac1n\sum_{i=1}^n\psi_{L_\psi}''(\bm x_i^\top\bbeta)\bm x_i\bm x_i^\top\preceq\gamma I_d
	\ \text{for every }\bbeta\in\Theta_0\right)\geq1-n^{-2}.
	\end{align}
\end{proposition}
The proof, given in Section \ref{sec: proof of prop: low-dim glm empirical curvature}, first obtains uniform population curvature from the unsaturated central region of the predictor and then proves uniform empirical-Hessian concentration with a covering argument.  The constants are uniform over $L_\psi\geq L_{\psi,0}$. This is the familiar population-curvature plus concentration route used in finite-sample GLM analysis \cite{mei2018landscape,ostrovskii2021finite}.

Basic composition gives overall $(\varepsilon,\delta)$-differential privacy once each of the $T$ updates is $(\varepsilon/T,\delta/T)$-differentially private.  Design clipping, response truncation, and \eqref{eq: saturated score bound} provide the required worst-case sensitivity bound.

\begin{proposition}\label{prop: non-sparse glm privacy}
	Under (G), for every deterministic $L_{\bm x},L_\psi,R<\infty$, choosing $B = 4(R + S_\psi)L_{\bm x}$ guarantees that Algorithm \ref{algo: private glm} is $(\varepsilon, \delta)$-differentially private.
\end{proposition}
Proposition \ref{prop: non-sparse glm privacy} is proved in Section \ref{sec: proof of prop: non-sparse glm privacy}. Although the privacy guarantee holds for any number of iterations $T$, choosing $T$ properly is crucial for accuracy because a larger value requires more privacy noise.

{Existing results on noisy gradient descent typically require $O(n)$ \cite{bassily2019private} or $O(n^2)$ \cite{bassily2014private} iterations for minimizing generic convex functions. On the event \eqref{eq: low-dim glm empirical curvature}, projected gradient descent is a contraction on the whole parameter set $\Theta_0$, so $O(\log n)$ iterations suffice.  Our proof works directly around the truth $\bbeta^*$ and retains both the empirical score and every privacy-noise vector; it does not require an empirical minimizer or a data-dependent warm start.}

\begin{Theorem}\label{thm: non-sparse glm upper bound}
		\rev{Let $\{(y_i, \bm x_i)\}_{i \in [n]}$ be an i.i.d. sample from the GLM \eqref{eq: glm definition}, let $\bbeta^*\in\Theta_0$, and suppose (D) and (G) hold. If $n\geq C_0d\log(en)$, initialize Algorithm \ref{algo: private glm} at $\bbeta^0=\bm0$. For $0<\varepsilon<1$ and $0<\delta\leq n^{-1}$, choose $\eta^0=2/(\alpha+\gamma)$, $L_{\bm x}=C_{\bm x}K_{\bm x}\sqrt{d+\log n}$, $L_\psi=L_{\psi,0}+C_\psi K_{\bm x}b_0\sqrt{\log n}$ for sufficiently large absolute constants $C_{\bm x},C_\psi$, $S_\psi=|\psi'(0)|+2c_2L_\psi$, $R=|\psi'(0)|+c_2L_\psi+\sqrt{6c_2c(\sigma)\log n}$, $B=4(R+S_\psi)L_{\bm x}$, and $T=\lceil\log(nb_0\vee n)/(-\log\rho)\rceil$, where $\rho=\max\{1/2,(\gamma-\alpha)/(\gamma+\alpha)\}$. Then Algorithm \ref{algo: private glm} is $(\varepsilon,\delta)$-differentially private and}
	\begin{align}
		\rev{\sup_{\bbeta^*\in\Theta_0}\E_{\bbeta^*}\|\bbeta^T - \bbeta^*\|^2_2 \lesssim c(\sigma)\left\{{\frac{d}{n}} + \frac{d(d+\log n) \log(1/\delta)\log^3 n}{n^2\varepsilon^2}\right\}} \label{eq: non-sparse glm upper bound}
	\end{align}
	\rev{where constants may depend on the fixed model constants and $b_0$, but not on $n,d,\varepsilon,$ or $\delta$.}
\end{Theorem}
Theorem \ref{thm: non-sparse glm upper bound} is proved in Section \ref{sec: proof of thm: non-sparse glm upper bound}. \rev{Since (D) implies (SG), both theorems apply under the upper bound's conditions. On the common parameter space $\Theta_0$, their bounds have matching polynomial dependence on $(n,d,\varepsilon)$ up to logarithmic factors: $d+\log n\asymp d$ when $d\gtrsim\log n$, and otherwise the additional ratio is at most $O(\log n)$.}

Another important implication of Theorems \ref{thm: low-dim glm lb} and \ref{thm: non-sparse glm upper bound} is the impact of differential privacy on the rate of convergence of estimating GLM parameters. As the first $O(d/n)$ term is the statistical rate of convergence in non-private estimation, the cost of differential privacy is negligible whenever $\varepsilon \gtrsim \rev{\sqrt{(d+\log n)\log(1/\delta)\log^3 n/n}}$. When $\varepsilon$ is less than this order, the rate of convergence is slower than its non-private \rev{counterpart}. \rev{Theorem \ref{thm: low-dim glm lb} is stated under $d\lesssim n\varepsilon$; claims for smaller $\varepsilon$ require a separate saturated lower-bound argument and are not asserted here.}

{The deterministic covariate clipping removes any almost-sure bounded-design assumption.  The smooth saturation \eqref{eq: saturated link score} similarly removes the global boundedness assumption on $\psi'$.  It bounds worst-case sensitivity while agreeing with the original score at the true predictors with high probability.  In contrast to hard clipping the gradient, the saturated score is the derivative of a convex potential and retains uniform empirical curvature, which is what allows the same contraction proof and rate.}


\section{The Bradley-Terry-Luce Model}\label{sec: BTL}

Pairwise comparisons arise in recommendation systems \cite{balakrishnan2012two}, sports tournaments \cite{masarotto2012ranking}, and education \cite{heldsinger2010using}. The Bradley--Terry--Luce (BTL) model assigns each item a latent strength and makes the probability of one item beating another depend on their strength difference. We study private estimation of these latent parameters under squared $\ell_2$ loss; a ranking can then be obtained by sorting the estimates.

Suppose there are $n$ items indexed by $[n] = \{1, 2,  \cdots, n\}$. We observe comparisons between pairs of items as follows.
\begin{itemize}
		\item Each pair of items indexed by $1 \leq i < j \leq n$ is compared independently with probability $0 < p < 1$. The $n$ items form a ``comparison graph'' in which an edge $(i,j)$ is present if and only if items $i$ and $j$ are compared. Let $\mathcal G$ denote its edge set.
	\item Each item $i$ is associated with a latent parameter $\theta_i \in [-1, 1]$. Given $\mathcal G$, the outcome of a comparison between items $i$ and $j$ is encoded by a Bernoulli random variable $Y_{ij}$ which takes the value $1$ if $i$ wins. The distribution of $Y_{ij}$ is independent of any other pair and determined by the latent parameters:
	\begin{align*}
		\Pro(Y_{ij} = 1) = \frac{e^{\theta_i}}{e^{\theta_i} + e^{\theta_j}}.
	\end{align*}
\end{itemize}
{
The goal is to estimate the latent parameter $\bth=(\theta_1,\ldots,\theta_n)$ from the observed outcomes $\{Y_{ij}\}_{(i,j)\in\mathcal G}$. The comparison graph is treated as public. We use \emph{edge-level privacy}: two data sets are adjacent when they have the same graph and differ in the outcome of one observed comparison. Thus the protected unit is one comparison record, rather than the collection of all comparisons incident to one item.

Because the BTL likelihood is invariant under adding a common constant, we work throughout with the identifiable parameter space
\[
\Theta=\{\bth\in\R^n:\|\bth\|_\infty\leq1,\ \bm1^\top\bth=0\}.
\]
\furtherrev{The ambient decision space is $\mathcal D=\mathbb R^n$. We write $\mathcal M^{\rm edge}_{\varepsilon,\delta}$ for the edge-level private class $\mathcal M^{\rm edge}_{\varepsilon,\delta}(\mathbb R^n)$.}
For an edge-level $(\varepsilon,\delta)$-differentially private estimator $M$, write
\[
R_\bth(M)=\E_{\mathcal G}\E_{\bm Y\mid\bth,\mathcal G}\|M(\bm Y,\mathcal G)-\bth\|_2^2.
\]
The quantity of interest is $\inf_{M\in\mathcal M^{\rm edge}_{\varepsilon,\delta}}\sup_{\bth\in\Theta}R_\bth(M)$.

\subsection{The Privacy-constrained Minimax Lower Bound}\label{sec: ranking lower bound}

We construct a score attack on an identifiable product submodel. Let $d_0=\lfloor n/2\rfloor$, let $\bm v_k=\bm e_{2k-1}-\bm e_{2k}$, and let $\bm V=(\bm v_1,\ldots,\bm v_{d_0})$. For $\bm\eta\in[-1,1]^{d_0}$, set $\bth(\bm\eta)=\bm V\bm\eta\in\Theta$ and, for any estimator $M$ of $\bth$, define the induced estimator $\widehat{\bm\eta}(M)=\bm V^\top M/2$. For $1\leq i<j\leq n$, write $\bm x_{ij}=\bm e_i-\bm e_j$ and
\[
q_{ij}(\bm\eta)=\{1+\exp(-\bm x_{ij}^\top\bm V\bm\eta)\}^{-1}.
\]
For lower-bound purposes there is no loss in replacing an arbitrary estimator by its Euclidean projection onto the closed convex set $\Theta$: projection preserves privacy and cannot increase its distance to any $\bth\in\Theta$. Throughout Propositions \ref{prop: ranking attack soundness} and \ref{prop: ranking attack completeness}, $M$ is therefore understood to take values in $\Theta$.
The edge score attack is
\begin{align*}
A_{ij}=\1\{(i,j)\in\mathcal G\}\left\langle\widehat{\bm\eta}(M(\bm Y,\mathcal G))-\bm\eta,
\{Y_{ij}-q_{ij}(\bm\eta)\}\bm V^\top\bm x_{ij}\right\rangle.
\end{align*}

The soundness bound uses edge-level differential privacy together with the size and largest Laplacian eigenvalue of the random comparison graph.
\begin{proposition}\label{prop: ranking attack soundness}
Suppose $p\geq C_0\log n/n$ for a sufficiently large constant $C_0$, $0<\varepsilon<1$, and $M$ is a $\Theta$-valued, edge-level $(\varepsilon,\delta)$-differentially private estimator. Then, for every $\bth(\bm\eta)$ in the paired submodel,
\begin{align}\label{eq: ranking attack soundness}
\E_{\mathcal G,\bm Y\mid\bm\eta}\sum_{1\leq i<j\leq n}A_{ij}
\leq C\varepsilon n^{3/2}p\sqrt{R_{\bth(\bm\eta)}(M)}+Cn^2p\delta+1.
\end{align}
\end{proposition}

For completeness, assign the independent density $\pi(t)=\1(|t|<1)(15/16)(1-t^2)^2$ to the coordinates of $\bm\eta$.
\begin{proposition}\label{prop: ranking attack completeness}
Suppose $M$ takes values in $\Theta$ and $\sup_{\bth\in\Theta}R_\bth(M)\leq c_0n$ for a sufficiently small absolute constant $c_0$. Under the preceding product prior on $\bm\eta$, there is an absolute constant $c>0$ such that
\begin{align}\label{eq: ranking attack completeness}
\E_{\bm\eta}\E_{\mathcal G,\bm Y\mid\bm\eta}\sum_{1\leq i<j\leq n}A_{ij}\geq cn.
\end{align}
\end{proposition}

Combining these propositions gives a lower bound on the same centered parameter space and under the same edge-level adjacency relation used by the upper bound.
\begin{Theorem}\label{thm: ranking lower bound}
Suppose $p\geq C_0\log n/n$, $0<\varepsilon<1$, and $0\leq\delta\leq c_0/(np)$ for sufficiently large $C_0$ and sufficiently small $c_0$. Then
\begin{align}\label{eq: ranking lower bound}
\inf_{M\in\mathcal M^{\rm edge}_{\varepsilon,\delta}}\sup_{\bth\in\Theta}R_\bth(M)
\gtrsim \frac1p+\min\left\{n,\frac{1}{np^2\varepsilon^2}\right\}.
\end{align}
\furtherrev{All implicit constants in this theorem are universal.}
\end{Theorem}
The first term is the non-private minimax risk. The second is the edge-privacy cost and is saturated by the diameter of $\Theta$ in the very-high-privacy regime. \rev{In the unsaturated regime, this privacy term follows directly by comparing soundness and completeness of the edge score attack above; thus the lower bound is a concrete instance of the general mechanism in Section \ref{sec: general lb}.}

\subsection{Optimality of the Private BTL Minimax Lower Bound}\label{sec: ranking upper bound}

We next give a matching pure-private objective-perturbation estimator. Write
$\phi(t)=\log(1+e^t)$ and define its globally strongly convex extension
\[
\bar\phi(t)=
\begin{cases}
\phi(-2)+\phi'(-2)(t+2)+\kappa_0(t+2)^2/2,&t<-2,\\
\phi(t),&|t|\leq2,\\
\phi(2)+\phi'(2)(t-2)+\kappa_0(t-2)^2/2,&t>2,
\end{cases}
\qquad
\kappa_0=\frac{e^{-2}}{(1+e^{-2})^2}.
\]
Thus $\bar\phi=\phi$ on every comparison margin generated by $\Theta$, while
$\bar\phi''(t)\geq\kappa_0$ for all $t\in\R$. Let
\[
\bar{\mathcal L}(\bth;y)=\sum_{(i,j)\in\mathcal G}
\left[-y_{ij}\bm x_{ij}^\top\bth+\bar\phi(\bm x_{ij}^\top\bth)\right]
\]
and let $\bm L_{\mathcal G}=\sum_{(i,j)\in\mathcal G}\bm x_{ij}\bm x_{ij}^\top$ be the graph Laplacian. Put $\bm P=\bm I-\bm1\bm1^\top/n$. If $\lambda_2(\bm L_{\mathcal G})<e^{-1}np$, set $\hat\bth=\bm0$. Otherwise draw independent random variables $W_1,\ldots,W_n$ with density
$w\mapsto(\varepsilon/4)\exp(-\varepsilon|w|/2)$, put $\bm B=\bm P\bm W$, and define
\begin{align}\label{eq: ranking MLE definition}
\widetilde\bth
&=\argmin_{\bth\in\bm1^\perp}
\{\bar{\mathcal L}(\bth;y)+\langle\bm B,\bth\rangle\},
\qquad
\hat\bth=\Pi_\Theta(\widetilde\bth).
\end{align}
Here $\Pi_\Theta$ denotes Euclidean projection. The projected Laplace perturbation
$\bm B$ lives in the centered subspace and satisfies
$\E\|\bm B\|_2^2=8(n-1)/\varepsilon^2$.

\begin{proposition}\label{prop: ranking MLE privacy}
For every public comparison graph $\mathcal G$ and every $\varepsilon>0$, the estimator in \eqref{eq: ranking MLE definition} is edge-level $(\varepsilon,0)$-differentially private.
\end{proposition}

\begin{proposition}\label{prop: ranking MLE accuracy}
If $p\geq C_1\log n/n$ for a sufficiently large absolute constant $C_1$ and $0<\varepsilon<1$, then
\begin{align*}
\sup_{\bth\in\Theta}R_\bth(\hat\bth)
\lesssim \frac1p+\min\left\{n,\frac{1}{np^2\varepsilon^2}\right\}.
\end{align*}
\end{proposition}

\begin{Theorem}\label{thm: ranking upper bound}
Under the conditions of Proposition \ref{prop: ranking MLE accuracy}, the estimator \eqref{eq: ranking MLE definition} is edge-level $(\varepsilon,0)$-differentially private and satisfies
\begin{align}\label{eq: ranking upper bound}
\sup_{\bth\in\Theta}R_\bth(\hat\bth)
\lesssim \frac1p+\min\left\{n,\frac{1}{np^2\varepsilon^2}\right\}.
\end{align}
\furtherrev{The implicit constant is universal.}
\end{Theorem}

Theorems \ref{thm: ranking lower bound} and \ref{thm: ranking upper bound} match up to absolute constants, and the upper bound in fact holds under pure differential privacy. In the unsaturated regime, edge privacy is negligible relative to the statistical error when $\varepsilon\gtrsim1/\sqrt{np}$.
}


\begingroup

\section{High-dimensional Sparse GLMs}\label{sec: sparse GLMs}

We now consider the generalized linear model
\begin{align}
	f_{\bbeta}(y\mid\bm x)
	=h(y,\sigma)\exp\left\{\frac{y\bm x^\top\bbeta-
	\psi(\bm x^\top\bbeta)}{c(\sigma)}\right\},
	\qquad \bm x\sim f_{\bm x},                              \label{eq: sparse glm definition}
\end{align}
when $d\geq n$ and the unknown coefficient vector is $s^*$-sparse.  For a sparsity level $s$ and radius $B>0$, write
\begin{align*}
	\Theta_s(B)=\{\bbeta\in\mathbb R^d:\|\bbeta\|_0\leq s,
	\|\bbeta\|_2\leq B\},
	\qquad H_s=s\log(ed/s).
\end{align*}
As in Section \ref{sec: GLMs}, $c(\sigma)$ is treated as a fixed model constant bounded above and away from zero; implicit constants may depend on these fixed bounds.
\furtherrev{Here the ambient decision space is $\mathcal D=\mathbb R^d$, and $\mathcal M_{\varepsilon,\delta}$ abbreviates $\mathcal M_{\varepsilon,\delta}(\mathbb R^d)$. In particular, the decoding cells below partition every possible estimator output, not merely the sparse parameter space.}
Because $\log\binom ds\asymp H_s$ for $s\leq d/2$, $H_s$ measures the combinatorial complexity of choosing the unknown support.  An ordinary likelihood score differentiates along continuous changes of $\bbeta$; on a fixed support, it sees an $s$-dimensional submodel but does not by itself retain the choice among $\binom ds$ possible supports.  To capture this discrete entropy, we connect supports by single-coordinate swaps and replace the derivative score by the exact likelihood increment $dQ/dP-1$.  Telescoping these increments produces a support-swap score attack with the same soundness--completeness structure as in Section \ref{sec: general lb}.

\subsection{A Support-Swap Score Lower Bound}\label{sec: sparse glm lower bound}

Let $P_{\bbeta}$ denote the law of one observation from \eqref{eq: sparse glm definition}.  Let $\mathcal V=\{S_1,\ldots,S_N\}$ be a constant-weight packing: every $S_v$ is an $s$-element subset of $[d]$, and distinct supports are separated in symmetric difference.  For an amplitude $a>0$, set
\begin{align*}
	\bbeta^v=a\bm 1_{S_v},\qquad P_v=P_{\bbeta^v},
	\qquad v\in[N].
\end{align*}
Here $\bm1_S$ is the coordinate indicator of $S$; in particular, $\|\bbeta^v\|_2=a\sqrt{s}$, so the packing lies in $\Theta_s(B)$ whenever $a\sqrt{s}\leq B$.

The Johnson graph $J(d,s)$ has all $s$-subsets of $[d]$ as its vertices, with two supports joined by an edge when one can be obtained from the other by deleting one coordinate and inserting one coordinate.  Thus every edge preserves sparsity, and the graph distance between $S$ and $S'$ is $|S\setminus S'|=|S'\setminus S|$.  For each ordered pair $(v,k)$, choose a shortest path
\begin{align*}
	S_{v,k,0}=S_k,S_{v,k,1},\ldots,S_{v,k,L_{vk}}=S_v
\end{align*}
	in $J(d,s)$.  Write
	$\bbeta_{v,k,\ell}=a\bm 1_{S_{v,k,\ell}}$ and
	$P_{v,k,\ell}=P_{\bbeta_{v,k,\ell}}$.  Because the parameter in
	\eqref{eq: sparse glm definition} enters through a positive exponential
	factor, all $P_{\bbeta}$ have the same support relative to the underlying
	base measure.  The finite-difference likelihood score is therefore well
	defined and is given by
\begin{align}
	T_{v,k,\ell}(z)
	&=\frac{dP_{v,k,\ell}}{dP_{v,k,\ell-1}}(z)-1\notag\\
	&=\exp\left[\frac{y\bm x^\top
	(\bbeta_{v,k,\ell}-\bbeta_{v,k,\ell-1})
	-\psi(\bm x^\top\bbeta_{v,k,\ell})
	+\psi(\bm x^\top\bbeta_{v,k,\ell-1})}{c(\sigma)}\right]-1. \label{eq: support swap likelihood score}
\end{align}
This is an exact likelihood-ratio increment from one support on the path to the next.

To turn estimation into decoding, assign every $\bm u\in\mathbb R^d$ to its nearest packing point by
\begin{align*}
\widehat v(\bm u)=\min\left\{j\in[N]:
\|\bm u-\bbeta^j\|_2=\min_{k\in[N]}\|\bm u-\bbeta^k\|_2\right\},
\qquad A_v=\{\bm u:\widehat v(\bm u)=v\}.
\end{align*}
The minimum index resolves distance ties deterministically, so $A_1,\ldots,A_N$ form a measurable partition of the estimator's output space.  For
\begin{align*}
	R_{v,k,\ell,i}
	=P_{v,k,\ell}^{\otimes i}\otimes
	P_{v,k,\ell-1}^{\otimes(n-i)},
\end{align*}
the first $i$ observations have law $P_{v,k,\ell}$ and the remaining observations have law $P_{v,k,\ell-1}$.  Thus $R_{v,k,\ell,i}$ interpolates between the two product experiments one observation at a time.  Define
\begin{align}
	\A_{v,k,\ell,i}(z_i,M(\bm Z))
	=\1\{M(\bm Z)\in A_v\}T_{v,k,\ell}(z_i).               \label{eq: support swap score attack}
\end{align}
The change-of-measure identity
\begin{align*}
	\E_{P}\!\left[h(Z)\left\{\frac{dQ}{dP}(Z)-1\right\}\right]
	=\E_Qh(Z)-\E_Ph(Z)
\end{align*}
is the exact finite-difference analogue of the likelihood-score identity.  Summing \eqref{eq: support swap score attack} over $i$ changes all $n$ observations from $P_{v,k,\ell-1}$ to $P_{v,k,\ell}$; summing next over $\ell$ moves from $P_k^{\otimes n}$ to $P_v^{\otimes n}$.  This double telescoping is the discrete counterpart of the derivative calculation in Section \ref{sec: general lb}.

\begin{proposition}[Support-swap score attack]\label{prop: stratified sparse attack}
	Suppose $N\geq4$ and $\|\bbeta^v-\bbeta^k\|_2\geq2r$ for $v\neq k$.  Let $Q_v$ be the law of $M(\bm Z)$ under $P_v^{\otimes n}$ and put
	\begin{align*}
		R=\max_{v\in[N]}\E_{P_v^{\otimes n}}
		\|M(\bm Z)-\bbeta^v\|_2^2.
	\end{align*}
	Then the averaged attack
	\begin{align*}
		\overline{\A}
		=\frac1{N^2}\sum_{v,k=1}^N
		\sum_{\ell=1}^{L_{vk}}\sum_{i=1}^n
		\E_{R_{v,k,\ell,i-1}}\A_{v,k,\ell,i}
	\end{align*}
	satisfies the exact identity and the completeness bound
	\begin{align}
		\overline{\A}
		&=\frac1{N^2}\sum_{v,k=1}^N\{Q_v(A_v)-Q_k(A_v)\},       \label{eq: support swap telescope}\\
		\overline{\A}
		&\geq1-\frac{R}{r^2}-\frac1N.                          \label{eq: support swap completeness}
	\end{align}
	Suppose in addition that every pair $P_v^{\otimes n},P_k^{\otimes n}$ has a coupling with expected Hamming distance at most $D$, and that $M$ is $(\varepsilon,\delta)$-differentially private.  With $m=\lfloor20D\rfloor$ and
	$\delta_m=\delta\sum_{j=0}^{m-1}e^{j\varepsilon}$, where the sum is zero when $m=0$, the same attack satisfies the soundness bound
	\begin{align}
		\overline{\A}
		\leq\frac{e^{m\varepsilon}-1}{N}+\frac1{20}+\delta_m.  \label{eq: support swap soundness}
	\end{align}
	Consequently, there are universal constants $c,C>0$ such that
	\begin{align*}
		0<\varepsilon\leq1,\qquad
		\varepsilon D\leq c\log N,\qquad
		\delta\leq c\varepsilon e^{-C\varepsilon D}
	\end{align*}
	imply $R\geq cr^2$.
\end{proposition}

The two bounds have a direct statistical interpretation.  If $R\ll r^2$, nearest-neighbor decoding succeeds under each $P_v^{\otimes n}$, and \eqref{eq: support swap completeness} makes the attack close to one.  For soundness, couple data sets from $P_v^{\otimes n}$ and $P_k^{\otimes n}$ and apply group privacy on the event that they differ in at most $m$ records.  Since the cells $A_v$ are disjoint, averaging the resulting decoder probabilities contributes the factor $1/N$ in \eqref{eq: support swap soundness}.  Thus accuracy makes the attack complete, while privacy makes it sound; balancing the two requires $\varepsilon D\lesssim\log N$.  Section \ref{sec: proof of prop: stratified sparse attack} gives the full proof.

\begin{Theorem}\label{thm: high-dim glm lb}
	Fix $B>0$, let $2\leq s^*\leq d/2$, and write
	\begin{align*}
		H=H_{s^*},\qquad \ell_n=\log(en).
	\end{align*}
	Suppose
	$\lambda_{\max}\{\E(\bm x\bm x^\top)\}\leq C$, and
	$\|\psi''\|_\infty\leq c_2$.  There are constants $c_0,c_1,c_2'>0$, depending only on these bounds, such that, if
	\begin{align*}
		0<\varepsilon\leq1,\qquad
		0\leq\delta\leq c_0\varepsilon
		\exp\{-c_1H/\sqrt{\ell_n}\},
	\end{align*}
	and
	\begin{align*}
		c(\sigma)\left\{\frac{H}{n}
		+\frac{H^2}{\ell_n n^2\varepsilon^2}\right\}
		\leq c_2'B^2,
	\end{align*}
	then
	\begin{align}
		\inf_{M\in\mathcal M_{\varepsilon,\delta}}
		\sup_{\bbeta\in\Theta_{s^*}(B)}
		\E_{\bbeta}\|M(\bm Z)-\bbeta\|_2^2
		\gtrsim c(\sigma)\left\{
		\frac{H}{n}
		+\frac{H^2}{\ell_n n^2\varepsilon^2}
		\right\}.                                             \label{eq: high-dim glm lb}
	\end{align}
	\furtherrev{The implicit constant may depend only on the displayed covariance and curvature bounds and on the fixed bounds for $c(\sigma)$, but not on $n,d,s^*,B,\varepsilon$, or $\delta$.}
\end{Theorem}

The proof, in Section \ref{sec: proof of thm: high-dim glm lb}, uses the same constant-weight packing for two separate arguments.  Its cardinality satisfies $\log N\asymp H$, and its squared separation is of order $a^2s^*$.  Choosing $a^2\asymp c(\sigma)H/(ns^*)$ and applying Fano's inequality gives the non-private term $c(\sigma)H/n$.

For the privacy term, choose
\begin{align*}
	a^2\asymp\frac{c(\sigma)H^2}{\ell_n n^2\varepsilon^2s^*}
\end{align*}
and note that a score-path calculation gives
$d_{\rm TV}(P_{\bbeta},P_{\bbeta'})\lesssim
\|\bbeta-\bbeta'\|_2/\sqrt{c(\sigma)}$.  Coordinatewise maximal couplings therefore produce data sets whose expected Hamming distance satisfies
$D\lesssim na\sqrt{s^*/c(\sigma)}$, and hence
$\varepsilon D\lesssim H/\sqrt{\ell_n}\leq C\log N$.
The condition on $\delta$ controls the corresponding group-privacy remainder in \eqref{eq: support swap soundness}.  Proposition \ref{prop: stratified sparse attack} then gives risk of order $a^2s^*=c(\sigma)H^2/\{\ell_n n^2\varepsilon^2\}$.  The factor $\ell_n$ is the cost of reducing $\varepsilon D$ from order $H$ to order $H/\sqrt{\ell_n}$; this weakens the sufficient condition on $\delta$ from an exponent of order $H$ to the exponent displayed in the theorem.  Thus the support-swap score attack supplies the privacy term, while Fano's inequality supplies only the non-private term.

\subsection{\texorpdfstring{\furtherrev{Private}}{Private} Iterative Hard Thresholding}\label{sec: sparse glm upper bound}

Ordinary iterative hard thresholding takes a gradient step and retains only the largest coordinates.  Its private analogue selects the coordinates sequentially using noisy magnitudes and then perturbs the selected values.  We analyze the resulting recursion directly around the true parameter $\bbeta^*$ under restricted Hessian bounds.

For a vector-valued statistic $\bm v(\bm Z)$ with $\ell_\infty$ sensitivity at most $\lambda$, the mechanism is as follows.

\begin{algorithm}[H]
	\caption{Noisy Hard Thresholding (NoisyHT)}\label{algo: noisy hard thresholding}
	\SetAlgoLined
	\SetKwInOut{Input}{Input}
	\SetKwInOut{Output}{Output}
	\Input{$\bm v\in\mathbb R^d$, sparsity $s$, privacy parameters $(\varepsilon,\delta)$, sensitivity $\lambda$.}
	Set $S=\varnothing$ and $b=2\lambda\sqrt{3s\log(1/\delta)}/\varepsilon$\;
	\For{$r=1$ \KwTo $s$}{
		Draw $w_{r1},\ldots,w_{rd}\stackrel{\mathrm{i.i.d.}}{\sim}\mathrm{Laplace}(b)$\;
		Append an index in $\arg\max_{j\notin S}\{|v_j|+w_{rj}\}$ to $S$, using deterministic tie breaking\;
	}
	Draw $\widetilde w_1,\ldots,\widetilde w_d\stackrel{\mathrm{i.i.d.}}{\sim}\mathrm{Laplace}(b)$\;
	\Output{$\bm v_S+\widetilde{\bm w}_S$.}
\end{algorithm}

\begin{Lemma}[\cite{dwork2018differentially,cai2021cost}]\label{lm: noisy hard thresholding privacy}
	If $\|\bm v(\bm Z)-\bm v(\bm Z')\|_\infty\leq\lambda$ for every pair of adjacent data sets, Algorithm \ref{algo: noisy hard thresholding} is $(\varepsilon,\delta)$-differentially private.
\end{Lemma}

For $L>0$, define coordinatewise clipping by
\begin{align*}
	[\mathsf C_L^{(\infty)}(\bm x)]_j
	=\operatorname{sign}(x_j)\min\{|x_j|,L\},\qquad j\in[d].
\end{align*}
For a deterministic level $L_{\bm x}$, set
$\overline{\bm x}_i=\mathsf C_{L_{\bm x}}^{(\infty)}(\bm x_i)$ and let
\begin{align*}
	\L_{n,L_\psi}^R(\bbeta)&=\frac1n\sum_{i=1}^n
	\{\psi_{L_\psi}(\bm x_i^\top\bbeta)-\Pi_R(y_i)\bm x_i^\top\bbeta\},\\
	\overline\L_{n,L_\psi}^R(\bbeta)&=\frac1n\sum_{i=1}^n
	\{\psi_{L_\psi}(\overline{\bm x}_i^\top\bbeta)
	-\Pi_R(y_i)\overline{\bm x}_i^\top\bbeta\}.
\end{align*}
	The design clipping and the saturated link score
	\eqref{eq: saturated link score} give a deterministic
	$\ell_\infty$-sensitivity bound for privacy. Condition (D) from
	Section \ref{sec: non-sparse glm upper bound} implies that no sample
	coordinate is clipped with high probability and supplies the
	restricted-curvature, predictor, and score-concentration bounds for
	accuracy.
	We append a final projection onto the unit Euclidean ball so that the estimator's risk is uniformly bounded; this post-processing step alters neither privacy nor sparsity.

\begin{algorithm}[H]
	\caption{Differentially Private Sparse Generalized Linear Regression}\label{algo: private sparse glm}
	\SetAlgoLined
	\SetKwInOut{Input}{Input}
	\SetKwInOut{Output}{Output}
	\SetKwFunction{NoisyHT}{NoisyHT}
	\Input{$\overline\L_{n,L_\psi}^R$, sparsity $s$, step size $\eta^0$, privacy parameters $(\varepsilon,\delta)$, noise bound $B_0$, iterations $T$, initial value $\bbeta^0$.}
	\For{$t=0$ \KwTo $T-1$}{
		$\bbeta^{t+1}=\NoisyHT\!\left(\bbeta^t-\eta^0\nabla\overline\L_{n,L_\psi}^R(\bbeta^t),
		s,\varepsilon/T,\delta/T,\eta^0B_0/n\right)$\;
	}
		Set $\widehat\bbeta=\Pi_{\{\bm u:\|\bm u\|_2\leq1\}}(\bbeta^T)$\;
		\Output{$\widehat\bbeta$.}
\end{algorithm}

We use the following restricted-curvature condition.  For the algorithmic sparsity $s\geq s^*$, let $\mathcal E_{\rm sc}(\alpha,\gamma;L_\psi)$ be the event that, for every $I\subset[d]$ with $|I|\leq2s+s^*$ and every $\bbeta$ on a line segment joining $\bbeta^*$ to an $s$-sparse vector in $\{\bm u:\|\bm u-\bbeta^*\|_2\leq2\}$,
\begin{align}
	\alpha\bm I_I\preceq
	[\nabla^2\L_{n,L_\psi}^R(\bbeta)]_{I,I}
	\preceq\gamma\bm I_I.                                  \label{eq: sparse glm curvature event}
\end{align}
Event \eqref{eq: sparse glm curvature event} is the uniform restricted-Hessian form of restricted strong convexity and smoothness needed by hard thresholding.  The next proposition derives this event from (D) and (G), rather than assuming it separately.

\begin{proposition}\label{prop: sparse glm empirical curvature}
	Suppose $d\geq n\geq2$, $1\leq s^*\leq s\leq d$, $\|\bbeta^*\|_0\leq s^*$, and $\|\bbeta^*\|_2\leq1$.  Under (D) and (G), there are fixed constants $0<\alpha\leq\gamma<\infty$, $C<\infty$, and $L_{\psi,0}^{\rm sp}<\infty$, depending only on the constants in (D) and (G), such that, for every deterministic $L_\psi\geq L_{\psi,0}^{\rm sp}$,
	\begin{align}
		n\geq C(s+s^*)\log(ed)
		\quad\Longrightarrow\quad
		\Pr\{\mathcal E_{\rm sc}(\alpha,\gamma;L_\psi)\}\geq1-n^{-3}.
		\label{eq: sparse glm empirical curvature probability}
	\end{align}
\end{proposition}

The proposition is uniform over all admissible $\bbeta^*$.  Its proof first establishes population curvature uniformly over all $(s+s^*)$-sparse vectors in the ball of radius three.  Sparse covering and Bernstein arguments then transfer that bound to every submatrix of size at most $2s+s^*$.  Section \ref{sec: proof of prop: sparse glm empirical curvature} gives the details.  This is the restricted analogue of Proposition \ref{prop: low-dim glm empirical curvature}; see also the standard restricted-curvature framework in \cite{loh2015regularized}.

\begin{Theorem}\label{thm: glm upper bound}
		Let $d\geq n\geq2$, $1\leq s^*\leq d$, $\|\bbeta^*\|_0\leq s^*$, and $\|\bbeta^*\|_2\leq1$.  Suppose (D) and (G) hold, and let $0<\alpha\leq\gamma<\infty$ be the constants in Proposition \ref{prop: sparse glm empirical curvature}.  For a sufficiently large
		universal constant $C_0$, assume
	\begin{align*}
		\left\lceil C_0(\gamma/\alpha)^2s^*\right\rceil\leq d,
	\end{align*}
	and set
	\begin{align*}
		s=\left\lceil C_0(\gamma/\alpha)^2s^*\right\rceil,
		\qquad
		\eta^0=2/(\alpha+\gamma),\qquad
		\bbeta^0=\bm0,
	\end{align*}
	$L_{\bm x}=C_{\bm x}K_{\bm x}\sqrt{\log(ed)}$,
	$L_\psi=L_{\psi,0}^{\rm sp}+C_\psi K_{\bm x}\sqrt{\log n}$,
	$S_\psi=|\psi'(0)|+2c_2L_\psi$, and
	$R=|\psi'(0)|+c_2L_\psi+C_R\sqrt{c(\sigma)\log n}$
	for sufficiently large absolute constants $C_{\bm x},C_\psi,C_R$,
	$B_0=4(R+S_\psi)L_{\bm x}$, \furtherrev{$q_0=1-3\alpha/(4\gamma)$, and
	\begin{align*}
		T=\left\lceil\frac{\log n}{-\log q_0}\right\rceil.
	\end{align*}}
	If $0<\varepsilon\leq1$, $0<\delta\leq n^{-1}$, and
	\begin{align*}
		n\gtrsim \sqrt{c(\sigma)}\,
		\frac{s^*\log d\sqrt{\log(ed)}
		\sqrt{\log(1/\delta)}\log^{3/2}n}{\varepsilon}
	\end{align*}
			with a sufficiently large implicit constant, then Algorithm \ref{algo: private sparse glm} is $(\varepsilon,\delta)$-differentially private and, with probability at least $1-Cn^{-3}$,
		\begin{align}
			\|\widehat\bbeta-\bbeta^*\|_2^2
		\lesssim c(\sigma)\left\{
		\frac{s^*\log d}{n}
		+\frac{(s^*\log d)^2\log(ed)
		\log(1/\delta)\log^3n}{n^2\varepsilon^2}
			\right\}.                                             \label{eq: glm upper bound}
		\end{align}
		Moreover,
		\begin{align}
		\E_{\bbeta^*}\|\widehat\bbeta-\bbeta^*\|_2^2
		&\lesssim c(\sigma)\left\{
		\frac{s^*\log d}{n}
		+\frac{(s^*\log d)^2\log(ed)
		\log(1/\delta)\log^3n}{n^2\varepsilon^2}
			\right\}.                                                \label{eq: sparse glm expected upper bound}
			\end{align}
		\furtherrev{All implicit constants and the constant $C$ in the probability bound may depend only on the fixed constants in (D) and (G), the fixed bounds on $c(\sigma)$, and the curvature ratio $\gamma/\alpha$, but not on $n,d,s^*,\varepsilon$, or $\delta$.}
\end{Theorem}

On the curvature, response no-truncation, and covariate no-clipping events used in the proof, Section \ref{sec: proof of thm: glm upper bound} establishes the recursion
\begin{align*}
	\|\bbeta^{t+1}-\bbeta^*\|_2^2
	\leq q\|\bbeta^t-\bbeta^*\|_2^2
	+C s\|\nabla\L_n(\bbeta^*)\|_\infty^2+C\mathcal W_t^2,
	\qquad q<1,
\end{align*}
where $\mathcal W_t$ is the NoisyHT perturbation.  Score concentration gives the first term in \eqref{eq: glm upper bound}; a simultaneous Laplace tail bound gives the second.

To compare the two theorems on the same parameter space, set $B=1$ in Theorem \ref{thm: high-dim glm lb}.  Suppose $s^*\leq d^{1-\kappa}$ for a fixed $\kappa>0$.  Then $H_{s^*}\asymp s^*\log d$, with constants depending only on $\kappa$, and the expected upper and lower bounds have the same polynomial dependence on $s^*$, $d$, $n$, and $\varepsilon$.  Their privacy terms differ by at most the factor $\log(ed)\log(1/\delta)\log^4n$: the covariate clipping contributes $\log(ed)$, one $\log(en)$ comes from the lower-bound tradeoff described after Theorem \ref{thm: high-dim glm lb}, and $\log(1/\delta)\log^3n$ comes from the private hard-thresholding algorithm.

The sufficient condition on $\delta$ in Theorem \ref{thm: high-dim glm lb} covers polynomially small $\delta$ when $H_{s^*}\lesssim\log^{3/2}(en)$. For larger support entropy the theorem requires a smaller value of $\delta$; we make no polynomial-$\delta$ matching claim in that regime.

\endgroup


\section{\texorpdfstring{\furtherrev{Nonparametric}}{Nonparametric} Function Estimation}\label{sec: nonparametric}
Although the score statistic is inherently a parametric concept, this section demonstrates that the score attack method can nonetheless yield optimal minimax lower bounds in \furtherrev{nonparametric} settings.

Consider $n$ pairs of random variables $\{(Y_i, X_i)\}_{i \in [n]}$ drawn i.i.d. from the model
	\begin{align*}
			Y_i = f(X_i) + \xi_i,\qquad X_i \sim U[0, 1],
	\end{align*}  
		where the noise term $\xi_i$ is independent of $X_i$ and follows $N(0, \sigma^2)$. Throughout this section, the noise level is a fixed model constant satisfying $0<\sigma_-\leq\sigma\leq\sigma_+<\infty$; implicit constants may depend on these fixed bounds, as well as on $\alpha$ and $C$. We estimate the unknown mean function $f: [0, 1] \to \R$ under $(\varepsilon, \delta)$-differential privacy and measure performance by the mean integrated squared risk (MISE),
	\begin{align*}
		R(\hat f, f) = \E \left[ \int_0^1 (\hat f(x) - f(x))^2 \d x \right],
	\end{align*}
	where the expectation is taken over the joint distribution of $\{(Y_i, X_i)\}_{i \in [n]}$. For a specified function class $\mathcal F$, the corresponding worst-case MISE is
	\begin{align*}
		R(\hat f, \mathcal F) = \sup_{f \in \mathcal F} 	R(\hat f, f) = \sup_{f \in \mathcal F}\E \left[ \int_0^1 (\hat f(x) - f(x))^2 \d x \right].
	\end{align*}
	That is, $R(\hat f, \mathcal F)$ measures the worst-case performance of $\hat f$ over the function class $\mathcal F$.
	{In this example, we take $\mathcal F$ to be the inhomogeneous periodic Sobolev class $\tilde W(\alpha, C)$ over $[0, 1]$: \furtherrev{throughout Section \ref{sec: nonparametric}, $\alpha\geq1$ is a fixed integer and $C>0$ is a fixed radius},
	\begin{align*}
		\tilde W(\alpha, C) = \left\{f:[0, 1] \to \R\ \middle|\
		\begin{array}{l}
		\displaystyle \int_0^1 \bigl(f(x)^2+(f^{(\alpha)}(x))^2\bigr)\d x\leq C^2,\\[2pt]
		 f^{(j)}(0)=f^{(j)}(1)\text{ for }0\leq j\leq\alpha-1
		\end{array}\right\}.
	\end{align*}
	The $L^2$ term controls the constant Fourier coefficient and makes the parameter set bounded.}
	\furtherrev{The ambient decision space is $\mathcal D=L^2[0,1]$. Accordingly, $\mathcal M_{\varepsilon,\delta}$ denotes the $\mathcal D$-valued private estimator class under the convention of Section \ref{sec: privacy constrained minimax risk}.} The privacy-constrained minimax risk of estimating $f$ is therefore
	\begin{align*}
		 \inf_{\hat f \in \mathcal M_{\varepsilon, \delta}} \sup_{f \in \tilde W(\alpha, C)} \E \left[ \int_0^1 (\hat f(x) - f(x))^2 \d x \right].
	\end{align*}
	
	Section \ref{sec: nonparametric lower bound} derives a lower bound by the score attack method, and Section \ref{sec: nonparametric upper bound} constructs an estimator with a matching risk upper bound.

	\subsection{The \texorpdfstring{\furtherrev{Nonparametric}}{Nonparametric} Minimax Lower Bound}\label{sec: nonparametric lower bound}
	Lower bounding the \furtherrev{nonparametric} privacy-constrained minimax risk is made easier by a sequence of reductions to parametric lower bound problems. The first step is to consider the orthogonal series expansion of $f \in \tilde W(\alpha, C)$ with respect to the Fourier basis
	\begin{align*}
		\rev{\varphi_1(t) = 1;\quad \varphi_{2k}(t) = \sqrt{2}\cos(2\pi k t),\quad \varphi_{2k+1}(t) = \sqrt{2}\sin(2\pi k t),\quad k = 1, 2, 3, \ldots.}
	\end{align*}
	We have $f = \sum_{j=1}^\infty \theta_j \varphi_j(x)$, where the Fourier coefficients are given by
\[
		\theta_j = \int_0^1 f(x)\varphi_j(x) \d x,\qquad j = 1, 2, 3, \cdots.
\]
	The Fourier coefficients allow a convenient representation of the periodic Sobolev class $\tilde W(\alpha, C)$: a function $f$ belongs to $\tilde W(\alpha, C)$ if and only if its Fourier coefficients belong to the ``Sobolev ellipsoid''
	{\begin{align}\label{eq: Sobolev ellipsoid definition}
		\Theta(\alpha, C) = \left\{\theta \in \R^{\mathbb Z^{+}}: \sum_{j=1}^\infty \tau_j^2 \theta_j^2 \leq C^2\right\},
	\end{align}
	where $\tau_1=1$ and $\tau_{2k}=\tau_{2k+1}=\sqrt{1+(2\pi k)^{2\alpha}}$ for $k\geq1$. Parseval's identity gives the equivalence between this ellipsoid and the preceding inhomogeneous Sobolev definition.} We can therefore define $\tilde W(\alpha, C)$ equivalently as
	\begin{align*}
		\tilde W(\alpha, C) = \left\{f = \sum_{j=1}^\infty \theta_j \varphi_j: \theta \in \Theta(\alpha, C)\right\}.
	\end{align*}
	
	This alternative definition of $\tilde W(\alpha, C)$ motivates a reduction from the original lower bound problem over an infinite-dimensional space, $\tilde W(\alpha, C)$, to a finite-dimensional lower bound problem. Specifically, for $k \in \mathbb N$, consider the $k$-dimensional subspace
	\begin{align*}
		\tilde W_k(\alpha, C) = \left\{f = \sum_{j=1}^\infty \theta_j \varphi_j: \theta \in \Theta(\alpha, C), \theta_j = 0 \text{ for every } j > k\right\}.
	\end{align*}
	It follows that $\tilde W_k(\alpha, C) \subseteq \tilde W(\alpha, C)$ for every $k$; in other words, for every $k$ we have
	\begin{align}\label{eq: nonparametric integrated risk reduction to finite}
		\inf_{\hat f \in \mathcal M_{\varepsilon, \delta}} \sup_{f \in \tilde W(\alpha, C)} \E \left[ \int_0^1 (\hat f(x) - f(x))^2 \d x \right] \geq \inf_{\hat f \in \mathcal M_{\varepsilon, \delta}} \sup_{f \in \tilde W_k(\alpha, C)} \E \left[ \int_0^1 (\hat f(x) - f(x))^2 \d x \right].
	\end{align}
	The next step is to find a minimax lower bound over each $k$-dimensional subspace, and then optimize $k$ to solve the original problem.
	
	\subsubsection{Finite-dimensional Minimax Lower Bounds via Score Attack}
	Once we focus on the $k$-dimensional subspace, the problem can be further simplified. For an estimator $\hat f$ and some $f \in \tilde W_k(\alpha, C)$, let $\{\hat\theta_j\}_{j \in \mathbb N}$ and $\{\theta_j\}_{j \in \mathbb N}$ be their respective Fourier coefficients. By the orthonormality of the Fourier basis, we have
	\begin{align}\label{eq: nonparametric integrated risk reduction to coefficient}
		\E \left[ \int_0^1 (\hat f(x) - f(x))^2 \d x \right] \geq \E \sum_{j=1}^k (\hat \theta_j - \theta_j)^2,
	\end{align}
	reducing the original problem into lower bounding the minimax mean squared risk of estimating a finite-dimensional parameter. Let $\Theta_k(\alpha, C)$ denote a finite-dimensional restriction of the Sobolev ellipsoid,
	\begin{align*}
		\Theta_k(\alpha, C) = \left\{\theta \in \R^k: \rev{\sum_{j=1}^k \tau_j^2 \theta_j^2 \leq C^2}\right\}, 
	\end{align*}
		The first $k$ Fourier coefficients of any private function estimator form a private vector estimator by post-processing; projecting this vector onto $\Theta_k(\alpha,C)$ cannot increase its distance to the true coefficient vector. We may therefore let $M(\bm X,\bm Y)$ be a $\Theta_k(\alpha,C)$-valued private estimator of $\bth=(\theta_1,\ldots,\theta_k)$. For $i\in[n]$, consider the score attack
	\begin{align*}
		\mathcal A(M(\bm X, \bm Y), (X_i, Y_i)) = \left\langle M(\bm X, \bm Y) - \bth,  \sigma^{-2}\left(Y_i - \sum_{j=1}^k \theta_j\varphi_j(X_i)\right)\bm\varphi(X_i)\right\rangle, 
	\end{align*}
		where $\bm\varphi$ denotes the vector-valued function $\bm\varphi: \R \to \R^k$, $\bm \varphi(x) = \left(\varphi_1(x), \varphi_2(x), \cdots, \varphi_k(x)\right)$.
		
		When the reference to $M$ and $(\bm X, \bm Y)$ is clear, write $A_{i} := \mathcal A(M(\bm X, \bm Y), (X_i, Y_i))$. To lower bound $\sup_{\bth \in \Theta_k(\alpha, C)} \E\|M(\bm X, \bm Y)-\bth\|_2^2$, we analyze $\sum_{i \in [n]} \E A_{i}$.
	
	\begin{proposition}\label{prop: nonparametric attack soundness}
			If $M$ is a $\Theta_k(\alpha,C)$-valued, $(\varepsilon, \delta)$-differentially private algorithm with $0 < \varepsilon < 1$ and $0<\delta<1/2$, then for every $\bth \in \Theta_k(\alpha, C)$,
		\begin{align}\label{eq: prop: nonparametric attack soundness}
				\sum_{i \in [n]} \E_{\bm X, \bm Y|\bth} A_{i} \leq \sigma^{-1}\left(2n\varepsilon \sqrt{\E_{\bm X, \bm Y|\bth}\|M(\bm X, \bm Y) - \bth\|_2^2} + C_1Cn\sqrt{k\log(1/\delta)}\delta\right),
			\end{align}
			where $C_1$ is universal.
	\end{proposition}
	The proof of Proposition \ref{prop: nonparametric attack soundness} is deferred to Section \ref{sec: proof of prop: nonparametric attack soundness}.
	
	{After upper bounding $\sum_{i \in [n]} \E_{\bm X, \bm Y|\bth} A_{i}$ at every $\bth \in \Theta_k(\alpha, C)$, we show that this sum is bounded away from zero on average under a product prior. Let
	\begin{align*}
		B^2=\frac{C^2}{2\sum_{j=1}^k\tau_j^2}\asymp k^{-(2\alpha+1)}
	\end{align*}
	and give each coordinate the density
	\begin{align*}
		q_B(t)=\frac{15}{16B}\left(1-\frac{t^2}{B^2}\right)^2\1(|t|<B).
	\end{align*}
	Because $|\theta_j|<B$ under this prior, $\sum_{j=1}^k\tau_j^2\theta_j^2<C^2/2$, so the prior is supported within $\Theta_k(\alpha,C)$. Moreover, $q_B$ vanishes at both endpoints, as required by Proposition \ref{prop: score attack stein's lemma}.}
	
	\begin{proposition}\label{prop: nonparametric attack completeness}
			\rev{Let $B^2=C^2/(2\sum_{j=1}^k\tau_j^2)$. Suppose $M$ takes values in $\Theta_k(\alpha,C)$ and satisfies}
				\begin{align*}
			\rev{\sup_{\bth \in \Theta_k(\alpha, C)} \E\|M(\bm X, \bm Y) - \bth\|_2^2 \leq \frac{kB^2}{40}.}
		\end{align*}
		\rev{If the coordinates of $\bth$ are independent with common density $q_B$,}
		then there is some constant $c > 0$ such that
		\begin{align}\label{eq: nonparametric attack completeness}
			\sum_{i \in [n]} \E_\bth \E_{\bm X, \bm Y|\bth} A_{i} > ck.
		\end{align}
	\end{proposition}
	The proposition is proved in Section \ref{sec: proof of prop: nonparametric attack completeness}. Together, Propositions \ref{prop: nonparametric attack soundness} and \ref{prop: nonparametric attack completeness} convert the two bounds on the expected attack into a finite-dimensional minimax lower bound.
	\begin{proposition}\label{prop: nonparametric coefficient lower bound}
			If $0 < \varepsilon < 1$ and $0 < \delta < cn^{-2}$ for a sufficiently small constant $c > 0$, then
		\begin{align}\label{eq: nonparametric coefficient lower bound}
			\inf_{M \in \furtherrev{\mathcal M_{\varepsilon, \delta}(\mathbb R^k)}} \sup_{\bth \in \Theta_k(\alpha, C)} \E\|M(\bm X, \bm Y)-\bth\|_2^2 \gtrsim \min\left(k^{-2\alpha}, \frac{k^2}{n^2\varepsilon^2}\right).
		\end{align}
	\end{proposition}
	The finite-dimensional lower bound is proved in Section \ref{sec: proof of prop: nonparametric coefficient lower bound}. Optimizing this bound over the truncation dimension $k$ recovers the \furtherrev{nonparametric} rate.
	
	\subsubsection{Optimizing the Finite-dimensional Lower Bounds}
	{By the reductions \eqref{eq: nonparametric integrated risk reduction to finite} and \eqref{eq: nonparametric integrated risk reduction to coefficient}, optimizing \eqref{eq: nonparametric coefficient lower bound} at $k\asymp(n\varepsilon)^{1/(\alpha+1)}$ gives the privacy term $(n\varepsilon)^{-2\alpha/(\alpha+1)}$. Independently, the usual non-private minimax lower bound over the same Sobolev ball is of order $n^{-2\alpha/(2\alpha+1)}$; it continues to hold after restricting the class of estimators to private estimators. Combining the two bounds and using $\max\{a,b\}\geq(a+b)/2$ gives the following result.}
	
	\begin{Theorem}\label{thm: nonparametric integrated risk lower bound}
			If $0 < \varepsilon < 1$, $0 < \delta < cn^{-2}$ for a sufficiently small constant $c > 0$, and $n\varepsilon \gtrsim 1$, then
		\begin{align}\label{eq: nonparametric integrated risk lower bound}
		\inf_{\hat f \in \mathcal M_{\varepsilon, \delta}} \sup_{f \in \tilde W(\alpha, C)} \E\left[\int_0^1 (\hat f(x)-f(x))^2 \d x\right] \gtrsim n^{-\frac{2\alpha}{2\alpha+1}} + (n\varepsilon)^{-\frac{2\alpha}{\alpha+1}}.
	\end{align}
	\furtherrev{The implicit constant may depend only on the fixed parameters $\alpha,C,\sigma_-$, and $\sigma_+$, but not on $n,\varepsilon$, or $\delta$.}
\end{Theorem}
		The first term can be recognized as the optimal MISE of function estimation in the periodic Sobolev class of order $\alpha$, and the second term is the cost of differential privacy. \rev{The second term is generated by the finite-dimensional score-attack bound and the optimization over $k$; only the first term is imported from the classical non-private theory. Thus the score attack remains the decisive privacy lower-bound step even after passing from the parametric Fourier coefficients back to the function class.} The next section shows the optimality of this \furtherrev{nonparametric} privacy-constrained lower bound by exhibiting an estimator with matching MISE up to the displayed logarithmic factors.
	
	\subsection{Optimality of the \texorpdfstring{\furtherrev{Nonparametric}}{Nonparametric} Lower Bound}\label{sec: nonparametric upper bound}
	Absent the differential privacy constraint, the $j$th Fourier coefficient of the mean function $f$ can be estimated by its empirical version, $\hat\theta_j = n^{-1}\sum_{i=1}^n Y_i\varphi_j(X_i)$, and the function $f$ is then estimated by $\hat f(x) = \sum_{j=1}^K \hat\theta_j \varphi_j(x)$ for some appropriately chosen $K$.
	
		We construct a private estimator of $f$ by privately estimating its Fourier coefficients; the resulting function estimator is private by post-processing. The sample mean $\hat\theta_j = n^{-1}\sum_{i=1}^n Y_i\varphi_j(X_i)$ lends itself naturally to an additive-noise mechanism, except that the Gaussian-distributed $Y_i$ are unbounded. We therefore truncate the responses before computing the empirical coefficients.
	
		Fix the number of released coefficients at $K$, and let $\bm\varphi$ denote the vector-valued function $\bm\varphi: \R \to \R^K$, $\bm \varphi(x) = \left(\varphi_1(x), \varphi_2(x), \cdots, \varphi_K(x)\right)$. The truncated empirical coefficient vector is
	\begin{align*}
		\frac{1}{n} \sum_{i=1}^n Y_i\1(|Y_i| \leq T) \cdot \bm \varphi(X_i).
	\end{align*}

	Over two adjacent data sets $D, D'$ with symmetric difference $\{(Y_i, X_i), (Y_i', X_i')\}$, their empirical coefficients differ by
	\begin{align*}
		\bm\Delta_{D, D'} = \frac{1}{n} \left(Y_i\1(|Y_i| \leq T) \cdot \bm \varphi(X_i) - Y'_i\1(|Y'_i| \leq T) \cdot \bm \varphi(X'_i)\right) \in \R^K.
	\end{align*}
	
		{Since $\sup_x\|\bm\varphi(x)\|_2^2\leq2K$, the replacement sensitivity of the truncated coefficient vector satisfies
	\[
	\|\bm\Delta_{D,D'}\|_2\leq\Delta_2:=\frac{2T\sqrt{2K}}n.
	\]
	For $0<\varepsilon<1$ and $0<\delta<1/2$, set
	\[
	\tau^2=\frac{2\Delta_2^2\log(1.25/\delta)}{\varepsilon^2}
	=\frac{16T^2K\log(1.25/\delta)}{n^2\varepsilon^2}.
	\]
	The Gaussian mechanism then shows that the following coefficient release is $(\varepsilon,\delta)$-differentially private:}
	\begin{align*}
		\tilde \bth_{K, T}  =  \frac{1}{n} \sum_{i=1}^n Y_i\1(|Y_i| \leq T) \cdot \bm \varphi(X_i) + \bm w, 
	\end{align*}
	\rev{where $\bm w\sim N(\bm0,\tau^2I_K)$. We use the associated function estimator}
	\begin{align}\label{eq: nonparametric estimator definition}
		\tilde f_{K, T} = \sum_{j=1}^K \left(\tilde \bth_{K, T}\right)_j\varphi_j(x).
	\end{align}
	
	\begin{Theorem} \label{thm: nonparametric upper bound}
			\rev{Let $r_{\alpha,C}$ be a uniform bound on $|f|$ over $\tilde W(\alpha,C)$. Suppose $n\geq2$, $0<\varepsilon<1$, $n\varepsilon\geq1$, and $0<\delta<1/2$. If $T=r_{\alpha,C}+2\sigma\sqrt{2\log n}$ and
		\[
		K=\max\left\{1,\left\lfloor c_1\min\left(n^{\frac{1}{2\alpha+1}},(n\varepsilon)^{\frac{1}{\alpha+1}}\right)\right\rfloor\right\}
		\]
		for some sufficiently small absolute constant $c_1\in(0,1]$, then $\tilde f_{K,T}$ is $(\varepsilon,\delta)$-differentially private and}
		\begin{align}\label{eq: nonparametric upper bound}
		\sup_{f \in \tilde W(\alpha, C)}\E \left[ \int_0^1 (\tilde f_{K, T}(x) - f(x))^2 \d x \right] \lesssim n^{-\frac{2\alpha}{2\alpha+1}} + \rev{(n\varepsilon)^{-\frac{2\alpha}{\alpha+1}} \log n\log(1/\delta)}.
	\end{align}
	\furtherrev{The implicit constant may depend only on the fixed parameters $\alpha,C,\sigma_-$, and $\sigma_+$, but not on $n,\varepsilon$, or $\delta$.}
\end{Theorem}
	Theorem \ref{thm: nonparametric upper bound} is proved in Section \ref{sec: proof of thm: nonparametric upper bound}. The risk upper bound \eqref{eq: nonparametric upper bound} matches the privacy-constrained minimax lower bound \eqref{eq: nonparametric integrated risk lower bound} \rev{in its polynomial dependence on $(n,\varepsilon)$, with a factor $\log n\log(1/\delta)$ in the upper bound. For polynomially small $\delta$, this is a logarithmic gap in $n$.} \rev{The lower bound is stated for $n\varepsilon\gtrsim1$; the regime $n\varepsilon=o(1)$ requires a separate saturated lower-bound argument.}
	
	The private \furtherrev{nonparametric} rate follows by reducing the function-estimation problem to a $k$-dimensional mean-estimation problem and choosing $k$ according to the smoothness of the mean function. \rev{The Gaussian release itself is not specific to the trigonometric basis, but a multivariate extension changes both the approximation dimension and the basis-norm sensitivity and therefore requires a separate rate calculation.}


\section{Discussion}\label{sec: discussion}

This paper introduced the score attack as a general method for deriving lower bounds on privacy-constrained minimax risk. Its soundness--completeness comparison is the common privacy-specific mechanism in all four developments: differential privacy upper bounds an average attack value, while a model-specific identity lower bounds the same value when estimation is accurate. The low-dimensional GLM, BTL, and \furtherrev{nonparametric} examples use the ordinary likelihood score and a Stein identity. For the sparse model, Proposition \ref{prop: stratified sparse attack} replaces an infinitesimal score by exact likelihood increments and telescopes them along paths of support swaps. Decoder completeness and privacy soundness then make the support entropy enter the lower bound directly. Together with the corresponding private estimators, these arguments determine the rates up to the logarithmic factors and parameter regimes stated in the individual theorems.

\textbf{Logarithmic gaps.} Some factors may result from truncation or composition and may be reducible by sharper algorithms. The parameter $\delta$ enters the soundness comparisons, and the support-swap attack avoids inflation by the packing cardinality. The lower bounds nevertheless do not recover the upper bounds' $\log(1/\delta)$ factors; a general characterization remains open \cite{steinke2017between,cai2021cost}. In the sparse problem, the additional $\log n$ in the lower rate weakens the sufficient condition on $\delta$ while preserving the full powers of $s^*$ and $\log(d/s^*)$.

\textbf{Pure differential privacy.} Every lower bound proved here for $(\varepsilon,\delta)$-differential privacy also applies to $(\varepsilon,0)$-differential privacy, since pure privacy further restricts the admissible estimators. Our upper-bound mechanisms generally use approximate privacy, except for the BTL estimator. Determining when pure and approximate privacy have different minimax rates remains an interesting problem.

\textbf{Beyond squared Euclidean loss.} The present framework is developed for squared $\ell_2$ loss. Other tasks may call for support recovery, ranking loss, confidence-set diameter, or testing error. Extending the soundness--completeness architecture to these losses would clarify whether a single score-based principle can play the role that Le Cam, Assouad, and Fano methods play across classical estimation and testing.

\textbf{Choice of priors for privacy-constrained estimation.} As in classical Bayes-minimax lower bounds, the prior must reflect the geometry that the desired rate is meant to expose. Smooth product priors are convenient for the coordinatewise Stein identities used in the low-dimensional GLM, BTL, and Fourier examples, and many nearby beta-type densities give the same rate. They do not, however, reveal the entropy of an unknown sparse support. The sparse analysis therefore uses a constant-weight packing, for which support swaps preserve sparsity and the finite likelihood increments telescope exactly. We do not claim that these choices are least favorable. Allowing a prior to depend on $(\varepsilon,\delta)$ may change the balance between completeness and the prior score, but it does not by itself remove the privacy-soundness remainder. Whether least-favorable priors exist for general privacy-constrained risks, and how to construct them, remain open questions.

\textbf{Practical membership inference attacks.} The score attack used in a lower-bound proof depends on the unknown true parameter and is therefore primarily a theoretical device. It can nevertheless suggest practical membership attacks if the parameter is replaced by an estimate computed from an independent public sample \cite{shokri2017membership}. For Gaussian mean estimation, this substitution recovers the tracing attack of \cite{homer2008resolving}. More generally, its effectiveness depends on whether soundness and completeness remain stable under the error of the public-sample surrogate.


\section{Omitted Proofs in Section \ref{sec: setup}} \label{sec: proofs}

This appendix proves Theorem \ref{thm: score attack general} and Propositions \ref{prop: score attack upper bound} and \ref{prop: score attack stein's lemma}.

\subsection{Proof of Theorem \ref{thm: score attack general}}\label{sec: proof of thm: score attack general}
\begin{proof}
	For soundness, we note that $\bm x_i$ and $M(\bm X'_i)$ are independent, and therefore
	\begin{align*}
			\E \mathcal A_\bth(\bm x_i, M(\bm X'_i)) = \E \langle M(\bm X'_i) - \bth, S_\bth(\bm x_i)\rangle = \langle \E M(\bm X'_i) - \bth,  \E S_\bth(\bm x_i)\rangle  = 0.
	\end{align*}
	The last equality holds because the score satisfies $\E S_\bth(\bm z) = \bm 0$ for any $\bm z \sim f_\bth$. For the first absolute moment, Jensen's inequality gives
	\begin{align*}
		&\E |\mathcal A_\bth(\bm x_i, M(\bm X'_i))| \leq \sqrt{\E \langle M(\bm X'_i) - \bth, S_\bth(\bm x_i)\rangle^2}  \\
		&\leq \sqrt{\E (M(\bm X'_i) - \bth)^\top (\Var S_\bth(\bm x_i))(M(\bm X'_i) - \bth)} \leq \sqrt{\E\|M(\bm X) - \bth\|^2_2} \sqrt{\lambda_{\max}(\mathcal I(\bth))}.
	\end{align*}
	
	For completeness, we first simplify
	\begin{align*}
		\sum_{i \in [n]} \E \A_\bth (\bm x_i, M(\bm X))
		= \E \Big\langle M(\bm X) - \bth, \sum_{i \in [n]} S_\bth(\bm x_i)\Big\rangle = \E \Big\langle M(\bm X), \sum_{i \in [n]} S_\bth(\bm x_i)\Big\rangle.
	\end{align*}
	By the definition of score and that $\bm x_1, \cdots, \bm x_n$ are i.i.d., $\sum_{i \in [n]} S_\bth(\bm x_i) = S_\bth(\bm x_1,\cdots, \bm x_n) = S_\bth(\bm X)$. It follows that
	\begin{align*}
		\E \Big\langle M(\bm X), \sum_{i \in [n]} S_\bth(\bm x_i)\Big\rangle = \E \Big\langle M(\bm X), S_\bth(\bm X)\Big\rangle = \sum_{j \in [d]} \E \left[M(\bm X)_j \frac{\partial}{\partial \theta_j} \log f_\bth^{(n)}(\bm X) \right].
	\end{align*}
	
		Let $\overline M_j(\bm x)=\E_M\{M_j(\bm x)\}$, where the expectation is over the internal randomness of $M$ for a fixed input. By the local domination assumption in the theorem, differentiation may be passed under the integral. Therefore, for every $j\in[d]$,
		\begin{align*}
			\E \left[M(\bm X)_j \frac{\partial}{\partial \theta_j} \log f_\bth^{(n)}(\bm X) \right]
			&=\int \overline M_j(\bm x)\,\partial_{\theta_j}f_\bth^{(n)}(\bm x)\,\d\bm x\\
			&=\frac{\partial}{\partial\theta_j}\int \overline M_j(\bm x)f_\bth^{(n)}(\bm x)\,\d\bm x
			= \frac{\partial}{\partial \theta_j}\E M(\bm X)_j.
		\end{align*}
\end{proof}

\subsection{Proof of Proposition \ref{prop: score attack upper bound}}\label{sec: proof of prop: score attack upper bound}
\begin{proof}
	Let $A_i := \A_\bth (\bm x_i, M(\bm X))$, $A'_i := \A_\bth (\bm x_i, M(\bm X'_i))$, and let $Z^+ = \max(Z, 0)$ and $Z^- = -\min(Z, 0)$ denote the positive and negative parts of a \rev{random variable} $Z$ respectively. We have
	\begin{align*}
		\E A_i = \E A_i^+ - \E A_i^- = \int_0^\infty \Pro(A_i^+ > t) \d t - \int_0^\infty \Pro(A_i^- > t) \d t.
	\end{align*}
	Condition on $\bm x_i$, its replacement $\bm x_i'$, and the remaining observations. Then $\bm X$ and $\bm X_i'$ are adjacent, and each event below is the same measurable function of the respective algorithm output. We may therefore apply differential privacy conditionally and then average over the data.
		For the positive part, if $0 < T < \infty$ and $0 < \varepsilon < 1$, differential privacy and $e^\varepsilon-1\leq2\varepsilon$ give
	\begin{align*}
		\int_0^\infty \Pro(A_i^+ > t) \d t &= \int_0^T \Pro(A_i^+ > t) \d t + \int_T^\infty \Pro(A_i^+ > t) \d t \\
		&\leq \int_0^T \left(e^\varepsilon\Pro({A'_i}^+ > t) + \delta\right)\d t + \int_T^\infty \Pro(A_i^+ > t) \d t \\
		&\leq \int_0^\infty \Pro({A'_i}^+ > t) \d t + 2\varepsilon\int_0^\infty \Pro({A'_i}^+ > t) \d t + \delta T + \int_T^\infty \Pro(|A_i| > t) \d t.
	\end{align*}
		For the negative part, reverse differential privacy gives
	\begin{align*}
		\int_0^\infty \Pro(A_i^- > t) \d t &= \int_0^T \Pro(A_i^- > t) \d t + \int_T^\infty \Pro(A_i^- > t) \d t \\
			& \geq \int_0^T \left(e^{-\varepsilon} \Pro({A'_i}^- > t) - \delta\right)\d t + \int_T^\infty \Pro(A_i^- > t) \d t \\
			&\geq \int_0^T \Pro({A'_i}^- > t) \d t - \varepsilon\int_0^T \Pro({A'_i}^- > t) \d t- \delta T \\
			&\geq \int_0^\infty \Pro({A'_i}^- > t) \d t - \varepsilon\int_0^\infty \Pro({A'_i}^- > t) \d t- \delta T
			-\int_T^\infty\Pro({A'_i}^->t)\,\d t.
	\end{align*}
	It then follows that
	\begin{align} \label{eq: score attack upper bound intermediate step}
		\E A_i &\leq \int_0^\infty \Pro({A'_i}^+ > t) \d t
		- \int_0^\infty \Pro({A'_i}^- > t) \d t \notag\\
			&\quad + 2\varepsilon \int_0^\infty \Pro(|A'_i| > t) \d t
			+ 2\delta T + \int_T^\infty \{\Pro(|A_i| > t)+\Pro(|A'_i|>t)\} \d t \notag \\
			&= \E A'_i + \rev{2\varepsilon\E|A'_i|} + 2\delta T
			+ \int_T^\infty \{\Pro(|A_i| > t)+\Pro(|A'_i|>t)\} \d t .
	\end{align}
	The proof is now complete by soundness \eqref{eq: soundness general}.
\end{proof}

\subsection{Proof of Proposition \ref{prop: score attack stein's lemma}}\label{sec: proof of prop: score attack stein's lemma}
\begin{proof}
	{
		Fix $j\in[d]$ and $\bth_{-j}$ outside a $\pi_{-j}$-null set. Apply Lemma \ref{lm: stein's lemma} to the one-dimensional function $h_j(t)=g_j(t,\bth_{-j})$. The assumed limits of $g_j(t,\bth_{-j})\pi_j(t)$ make both boundary terms vanish, and hence
	\begin{align*}
		\int \frac{\partial}{\partial t}g_j(t,\bth_{-j})\pi_j(t)\,\d t
		=\int -g_j(t,\bth_{-j})\frac{\pi'_j(t)}{\pi_j(t)}\pi_j(t)\,\d t.
	\end{align*}
	Integrating this identity with respect to $\pi_{-j}$ and using Fubini's theorem gives
	\begin{align*}
		\E_{\bm\pi}\frac{\partial}{\partial\theta_j}g_j(\bth)
		=\E_{\bm\pi}\left[-g_j(\bth)\frac{\pi'_j(\theta_j)}{\pi_j(\theta_j)}\right].
	\end{align*}
	Since $g(\bth)=\E_{\bm X|\bth}M(\bm X)$, Jensen's inequality yields
	$|g_j(\bth)-\theta_j|\leq\E_{\bm X|\bth}|M_j(\bm X)-\theta_j|$. Therefore,
	\begin{align*}
		\E_{\bm\pi}\frac{\partial}{\partial\theta_j}g_j(\bth)
		&\geq \E_{\bm\pi}\left[-\theta_j\frac{\pi'_j(\theta_j)}{\pi_j(\theta_j)}\right]\\
		&\quad-\E_{\bm\pi}\left[\E_{\bm X|\bth}|M_j(\bm X)-\theta_j|
		\left|\frac{\pi'_j(\theta_j)}{\pi_j(\theta_j)}\right|\right].
	\end{align*}
	Summing over $j$ and applying Cauchy--Schwarz proves
	\begin{align*}
		\E_{\bm\pi}\sum_{j\in[d]}\frac{\partial}{\partial\theta_j}g_j(\bth)
		&\geq \E_{\bm\pi}\sum_{j\in[d]}-\theta_j\frac{\pi'_j(\theta_j)}{\pi_j(\theta_j)}\\
		&\quad-\sqrt{\E_{\bm\pi}\E_{\bm X|\bth}\|M(\bm X)-\bth\|_2^2}
		\sqrt{\E_{\bm\pi}\sum_{j\in[d]}\left(\frac{\pi'_j(\theta_j)}{\pi_j(\theta_j)}\right)^2}.
	\end{align*}
	}
\end{proof}

	\bibliographystyle{plain}
	\bibliography{reference}	\label{lastpage}
	
	\newpage
	\appendix
	
\section{Omitted Proofs in Section \ref{sec: GLMs}}\label{sec: glm proofs}
\subsection{Proof of Proposition \ref{prop: low-dim glm attack soundness}}\label{sec: proof of prop: low-dim glm attack soundness}

\begin{proof}[Proof of Proposition \ref{prop: low-dim glm attack soundness}]
		In view of Theorem \ref{thm: score attack general}, we first calculate the score statistic of $f(y, \bm x)$ with respect to $\bbeta$ and the Fisher information matrix. The conditional support is parameter-independent, and (SG) together with the residual moment bound derived below gives finite score moments, so the standard exponential-family differentiation identities apply. We have
	\begin{align*}
		\frac{\partial}{\partial \bbeta} \log f(y, \bm x) &= \frac{\partial}{\partial \bbeta} \log \left(f_\bbeta(y|\bm x) f(\bm x)\right) = \frac{\partial}{\partial \bbeta} \log f_\bbeta(y|\bm x) \\
		&= \frac{\partial}{\partial \bbeta}\left(\frac{\bm x^\top \bbeta y - \psi(\bm x^\top \bbeta)}{c(\sigma)}\right) = \frac{[y - \psi'(\bm x^\top \bbeta)]\bm x}{c(\sigma)}.
	\end{align*}
	For the Fisher information, we have
	\begin{align*}
		\mathcal I(\bbeta) = -\E \left(\frac{\partial^2}{\partial \bbeta^2} \log f(y, \bm x) \right) = \E\left(\frac{\psi''(\bm x^\top \bbeta)}{c(\sigma)} \bm x \bm x^\top\right) \preceq \frac{c_2}{c(\sigma)} \E[\bm x \bm x^\top],
	\end{align*}
	where the last inequality holds by $\|\psi''\|_\infty \leq c_2$. We then have $\lambda_{\max}(\mathcal I(\bbeta)) \leq Cc_2 /c(\sigma)$ by $\lambda_{\max}(\E[\bm x \bm x^\top]) \leq C$. \rev{Write $A_i'$ for the corresponding independent-copy attack. The soundness part of Theorem \ref{thm: score attack general} gives $\E A_i'=0$ and, for every $i\in[n]$,}
	{
	\begin{align*}
		\E|A_i'|
		\leq \sqrt{\E\|M(\bm y,\bm X)-\bbeta\|_2^2}
		\sqrt{Cc_2/c(\sigma)}.
	\end{align*}}
	
		By Proposition \ref{prop: score attack upper bound}, we have
		\begin{align*}
			\E A_i &\leq 2\varepsilon \sqrt{\E\|M(\bm y, \bm X) - \bbeta\|^2_2} \sqrt{Cc_2/c(\sigma)} +  2\delta T \\
			&\quad+ \int_T^\infty \{\Pro(|A_i| > t)+\Pro(|A_i'|>t)\} \d t.
		\end{align*}
		We next control the two tail integrals without imposing an almost-sure
		bound on the design. Put
		\[
			Z_i=\frac{y_i-\psi'(\bm x_i^\top\bbeta)}{c(\sigma)},
			\qquad \nu^2=\frac{c_2}{c(\sigma)}.
		\]
		For the conditional law of the response, consider
		$f_{\theta}(y) = h(y, \sigma)\exp\left(\frac{y\theta - \psi(\theta)}{c(\sigma)}\right)$. Then
		\begin{align*}
			\E \exp\left(\frac{\lambda}{c(\sigma)} y\right)
		&= \int \exp\left(\frac{\lambda y}{c(\sigma)} \right)  h(y, \sigma)
		\exp\left(\frac{y\theta - \psi(\theta)}{c(\sigma)}\right) \d y\\
		&= \exp\left(\frac{\psi(\theta + \lambda) - \psi(\theta)}{c(\sigma)}\right).
	\end{align*}
		We may therefore compute the moment generating function of $Z_i$ conditional on $\bm x_i$:
	\begin{align}\label{eq: glm mgf upper bound}
		\log \E \exp\left(\lambda 
		\cdot \frac{y_i - \psi'(\bm x_i^\top\bbeta)}{c(\sigma)}\Big| \bm x_i\right) &= \frac{1}{c(\sigma)}\left(\psi(\bm x_i^\top\bbeta + \lambda) - \psi(\bm x_i^\top\bbeta) - \lambda\psi'(\bm x_i^\top\bbeta)\right) \notag \\
		&\leq \frac{1}{c(\sigma)} \cdot \frac{\lambda^2 \psi^{''}(\bm x_i^\top\bbeta + \tilde \lambda)}{2}
	\end{align}
			for some $\tilde \lambda$ between $0$ and $\lambda$. Because
			$\|\psi^{''}\|_\infty<c_2$,
			\[
				\E\{\exp(\lambda Z_i)\mid\bm x_i\}
				\leq\exp(\nu^2\lambda^2/2).
			\]
			Applying the Chernoff bound to both tails (for example,
			\cite{wainwright2019high}, equation (2.5)) gives, for every
			$u\geq1$,
			\begin{align}\label{eq: low-dim residual tail}
				\Pro\{|Z_i|>\sqrt{2}\nu\sqrt u\}\leq2e^{-u}.
			\end{align}

			Condition (SG) also gives a Euclidean-norm tail. To see this
			directly, let $\mathcal N$ be a $1/2$-net of the unit sphere
			with $|\mathcal N|\leq5^d$. Since
			$\|\bm x_i\|_2\leq2\max_{\bm v\in\mathcal N}
			|\bm v^\top\bm x_i|$, the sub-Gaussian tail bound and a union
			bound show that a constant $C_4$ satisfies
			\begin{align}\label{eq: low-dim design norm tail}
				\Pro\{\|\bm x_i\|_2>C_4K_{\bm x}(\sqrt d+\sqrt u)\}
				\leq2e^{-u},\qquad u\geq1.
			\end{align}
			Because both $M(\bm y,\bm X)$ and $\bbeta$ belong to
			$\Theta_0$,
			\[
				|A_i|\leq2b_0|Z_i|\,\|\bm x_i\|_2.
			\]
			The same pointwise inequality holds for $A_i'$ because the
			estimator evaluated on the adjacent data set also takes values
			in $\Theta_0$. Combining \eqref{eq: low-dim residual tail} and
			\eqref{eq: low-dim design norm tail} by a union bound gives a
			constant $a=C_5b_0K_{\bm x}\nu$ such that
			\begin{align}\label{eq: low-dim attack bernstein tail}
				\max\{\Pro(|A_i|>a(\sqrt{du}+u)),
				\Pro(|A_i'|>a(\sqrt{du}+u))\}
				\leq4e^{-u},\qquad u\geq1.
			\end{align}
			No independence between the estimator and $(y_i,\bm x_i)$ is
			used here; the deterministic range bound on the estimator
			makes the pointwise inequality valid.

			For $0<\delta<1/2$, let $L=\log(e/\delta)$ and set
			$T=a(\sqrt{dL}+L)$. The map
			$q(u)=a(\sqrt{du}+u)$ is increasing, so a change of variables
			in \eqref{eq: low-dim attack bernstein tail} yields
			\begin{align*}
				\int_T^\infty\Pro(|A_i|>t)\,\d t
				&\leq4\int_L^\infty e^{-u}q'(u)\,\d u\\
				&\leq C_6a\delta\left(\sqrt{d/L}+1\right).
			\end{align*}
			The same bound holds for $A_i'$. Since $L\geq1$, the two tail
			integrals and $2\delta T$ are together at most
			$C_7a\delta(\sqrt{dL}+L)$. Substitution into the privacy
			comparison and summation over $i$ give
			\begin{align}\label{eq: low-dim glm attack upper}
				\sum_{i \in [n]} \E A_i
				&\leq 2n\varepsilon
				\sqrt{\E\|M(\bm y, \bm X) - \bbeta\|^2_2}
				\sqrt{Cc_2/c(\sigma)}\notag\\
				&\quad+C_3n\delta\sqrt{c_2/c(\sigma)}
				\left\{\sqrt{d\ell_\delta}+\ell_\delta\right\},
			\end{align}
			where $\ell_\delta=\log(e/\delta)$ and the fixed factors
			$b_0,K_{\bm x}$ are absorbed into $C_3$. When $\delta=0$,
			let $T\to\infty$ in Proposition
			\ref{prop: score attack upper bound};
			\eqref{eq: low-dim attack bernstein tail} ensures
			integrability, and the last term in
			\eqref{eq: low-dim glm attack upper} is zero.
	
\end{proof}

\subsection{Proof of Proposition \ref{prop: low-dim glm attack completeness}}\label{sec: proof of prop: low-dim glm attack completeness}
\begin{proof}[Proof of Proposition \ref{prop: low-dim glm attack completeness}]
		{Put $g_j(\bbeta)=\E_{\bm y,\bm X\mid\bbeta}M_j(\bm y,\bm X)$. Because $M$ is bounded and the GLM has parameter-independent support, the local domination needed to differentiate $g_j$ follows from the exponential-family likelihood and the bounded score moments above. Completeness of the score attack and integration by parts under the product prior therefore give
	\begin{align*}
		\sum_{i=1}^n\E_{\bm\pi}\E_{\bm y,\bm X\mid\bbeta}A_i
	&=\sum_{j=1}^d\E_{\bm\pi}\partial_j g_j(\bbeta)\\
	&=d-\E_{\bm\pi}\sum_{j=1}^d\{g_j(\bbeta)-\beta_j\}\frac{p_a'(\beta_j)}{p_a(\beta_j)},
	\end{align*}
		where $p_a(t)=15\{1-(t/a)^2\}^2/(16a)$ on $[-a,a]$. Because $M$ takes values in the compact set $\Theta_0$, $g$ is bounded; hence $g_j(t,\bbeta_{-j})p_a(t)\to0$ at both endpoints. A direct calculation gives
	\[
	\int_{-a}^a\left\{\frac{p_a'(t)}{p_a(t)}\right\}^2p_a(t)\,\d t=\frac{10}{a^2}.
	\]
	Consequently, Cauchy--Schwarz and Jensen's inequality yield
	\begin{align*}
		\sum_{i=1}^n\E_{\bm\pi}\E_{\bm y,\bm X\mid\bbeta}A_i
	&\geq d-\left(\E_{\bm\pi}\|g(\bbeta)-\bbeta\|_2^2\right)^{1/2}\frac{\sqrt{10d}}a\\
	&\geq d-\frac{2\sqrt{10}\,d}{b_0}
		\left(\E_{\bm\pi}\E_{\bm y,\bm X\mid\bbeta}\|M-\bbeta\|_2^2\right)^{1/2}.
	\end{align*}
	The last term is at most $d/2$ when $c_0$ is chosen sufficiently small. This proves the claim.}
\end{proof}

\subsection{Proof of Theorem \ref{thm: low-dim glm lb}}\label{sec: proof of thm: low-dim glm lb}
\begin{proof}[Proof of Theorem \ref{thm: low-dim glm lb}]
	\rev{We prove both terms on $\Theta_0$.}
	
	{\furtherrev{Begin with an arbitrary $M\in\mathcal M_{\varepsilon,\delta}(\mathbb R^d)$.} Let $P_{\Theta_0}$ denote Euclidean projection onto $\Theta_0$. For every $\bbeta\in\Theta_0$, the projection theorem gives
	\[
	\|P_{\Theta_0}(M(\bm y,\bm X))-\bbeta\|_2\leq \|M(\bm y,\bm X)-\bbeta\|_2.
	\]
	Moreover, $P_{\Theta_0}\circ M$ remains $(\varepsilon,\delta)$-differentially private and has squared error at most $4b_0^2$. \furtherrev{Thus post-processing maps the full ambient estimator class to $\Theta_0$ without increasing risk, which justifies the range condition in Propositions \ref{prop: low-dim glm attack soundness} and \ref{prop: low-dim glm attack completeness}.}

	Write $\mathcal R(M)=\sup_{\bbeta\in\Theta_0}\E_{\bbeta}\|M-\bbeta\|_2^2$. If $\mathcal R(M)>c_0b_0^2$, the desired privacy term already follows, after reducing the implicit constant, from $d\lesssim n\varepsilon$. Otherwise Proposition \ref{prop: low-dim glm attack completeness} applies. Combining it with Proposition \ref{prop: low-dim glm attack soundness} and averaging over $\bm\pi$ gives}
	\begin{align*}
		d &\lesssim \sum_{i \in [n]}\E_{\bm \pi} \E_{\bm y, \bm X|\bbeta} A_i\\
		&\leq 2n\varepsilon\E_{\bm \pi}
		\sqrt{\E_{\bm y, \bm X|\bbeta} \|M(\bm y, \bm X) - \bbeta\|^2_2}
		\sqrt{Cc_2/c(\sigma)}\\
		&\quad + C_3n\delta\sqrt{c_2/c(\sigma)}
		\left\{\sqrt{d\ell_\delta}+\ell_\delta\right\}.
	\end{align*}
	The final term is zero when $\delta=0$. In that case the desired
	comparison follows immediately, so for the remainder of this argument
	assume $\delta>0$.
	It follows that
	\begin{align*}
		&2n\varepsilon\E_{\bm \pi}
		\sqrt{\E_{\bm y, \bm X|\bbeta} \|M(\bm y, \bm X) - \bbeta\|^2_2}
		\sqrt{Cc_2/c(\sigma)}\\
			&\qquad\gtrsim d - C_3n\delta\sqrt{c_2/c(\sigma)}
			\left\{\sqrt{d\ell_\delta}+\ell_\delta\right\}.
	\end{align*}
		For $d\geq1$,
		\[
			\frac{\sqrt{d\ell_\delta}+\ell_\delta}{d}
			\leq\sqrt{\ell_\delta}+\ell_\delta.
		\]
		On the range $\delta\lesssim n^{-(1+\gamma)}$, the functions
		$\delta\sqrt{\log(e/\delta)}$ and
		$\delta\log(e/\delta)$ are maximized, up to fixed constants, at the
		largest allowed value of $\delta$. Hence
		\[
			n\delta\{\sqrt{\ell_\delta}+\ell_\delta\}
			\leq C_\gamma c_\delta n^{-\gamma}
			\{\sqrt{\log(en)}+\log(en)\},
		\]
		where $c_\delta$ is the implicit constant in the condition on
		$\delta$. The right side is uniformly bounded in $n\geq2$ and can
		be made sufficiently small by choosing $c_\delta$ small enough.
		Thus the preceding display is bounded below by a fixed positive
		multiple of $d$. Jensen's inequality then gives
	\begin{align*}
		\E_{\bm \pi}\E_{\bm y, \bm X|\bbeta} \|M(\bm y, \bm X) - \bbeta\|^2_2 \gtrsim \frac{c(\sigma)d^2}{n^2\varepsilon^2}.
	\end{align*}
	Because the sup-risk is always greater than the Bayes risk, we have 
	\begin{align*}
		\rev{\sup_{\bbeta \in \Theta_0}}\E_{\bm y, \bm X|\bbeta} \|M(\bm y, \bm X) - \bbeta\|^2_2 \gtrsim \frac{c(\sigma)d^2}{n^2\varepsilon^2}.
	\end{align*}
	{For completeness, the non-private term also holds on $\Theta_0$. Indeed, take a hypercube $\bbeta_v=h v$, $v\in\{-1,1\}^d$, with $h=c\sqrt{c(\sigma)/n}$ and $c>0$ sufficiently small. The condition $d\lesssim n\varepsilon\leq n$ ensures that this cube is contained in $\Theta_0$. Taylor's theorem, the stated bound on $\psi''$, and the covariance upper bound give, for neighboring vertices,
	\[
	D_{\rm KL}(P_{\bbeta_v}^{\otimes n}\|P_{\bbeta_{v'}}^{\otimes n})
	\leq \frac{nc_2}{2c(\sigma)}\E\{\bm x^\top(\bbeta_v-\bbeta_{v'})\}^2\leq C_0
	\]
	with $C_0<1$ after decreasing $c$. Assouad's lemma therefore yields a lower bound of order $dh^2\asymp c(\sigma)d/n$. Combining the two terms completes the proof.}
\end{proof}

\subsection{Proof of Proposition \ref{prop: low-dim glm empirical curvature}}\label{sec: proof of prop: low-dim glm empirical curvature}

\begin{proof}[Proof of Proposition \ref{prop: low-dim glm empirical curvature}]
	Fix a deterministic $L_\psi\geq L_{\psi,0}$, where the lower
	threshold $L_{\psi,0}$ is chosen in Step 1 below.  For
	$\bbeta\in\Theta_0$, define the saturated-score population and empirical
	Hessians
	\[
	H_{L_\psi}(\bbeta)
	=\E\{\psi_{L_\psi}''(\bm x^\top\bbeta)\bm x\bm x^\top\},
	\qquad
	H_{n,L_\psi}(\bbeta)=\frac1n\sum_{i=1}^n
	\psi_{L_\psi}''(\bm x_i^\top\bbeta)\bm x_i\bm x_i^\top.
	\]
	We first show that $H_{L_\psi}(\bbeta)$ is uniformly well conditioned
	and then transfer this property to $H_{n,L_\psi}(\bbeta)$, with
	constants independent of $L_\psi$.

	\textit{Step 1: uniform population curvature.}
	The sub-Gaussian part of (D) implies that there are universal constants $a_1,a_2>0$ such that, for every unit vector $\bm u$,
	\[
	\E|\bm u^\top\bm x|^4\leq a_1K_{\bm x}^4,
	\qquad
	\Pro(|\bm x^\top\bbeta|>L)
	\leq2\exp\!\left\{-\frac{a_2L^2}{K_{\bm x}^2b_0^2}\right\}
	\quad(\bbeta\in\Theta_0).
	\]
	The second inequality is immediate for $\bbeta=\bm0$ and follows by scaling otherwise.  Hence Cauchy--Schwarz gives, uniformly over unit $\bm u$ and $\bbeta\in\Theta_0$,
	\begin{align*}
	\E\{(\bm u^\top\bm x)^2\1(|\bm x^\top\bbeta|>L)\}
	&\leq \{\E(\bm u^\top\bm x)^4\}^{1/2}
	\Pro(|\bm x^\top\bbeta|>L)^{1/2}\\
	&\leq \sqrt{2a_1}K_{\bm x}^2
	\exp\!\left\{-\frac{a_2L^2}{2K_{\bm x}^2b_0^2}\right\}.
	\end{align*}
	Choose a fixed $L_0<\infty$, depending only on $(C,K_{\bm x},b_0)$, so that the last display is at most $(2C)^{-1}$, and set
	$L_{\psi,0}=\max\{1,L_0\}$. By (G) and continuity,
	\[
	m_0:=\inf_{|t|\leq L_0}\psi''(t)>0.
	\]
	For $L_\psi\geq L_{\psi,0}$, equation
	\eqref{eq: saturated link score} gives
	$\psi_{L_\psi}''(t)=\psi''(t)$ whenever $|t|\leq L_0$ and
	$0\leq\psi_{L_\psi}''(t)\leq c_2$ everywhere.  Using
	$\E(\bm u^\top\bm x)^2\geq C^{-1}$, we obtain
	\begin{align*}
	\bm u^\top H_{L_\psi}(\bbeta)\bm u
	&\geq m_0\E\{(\bm u^\top\bm x)^2
	\1(|\bm x^\top\bbeta|\leq L_0)\}
	\geq \frac{m_0}{2C}=: \alpha_0.
	\end{align*}
	The upper bound follows directly from (D) and (G):
	\[
	\bm u^\top H_{L_\psi}(\bbeta)\bm u
	\leq c_2\E(\bm u^\top\bm x)^2
	\leq c_2C=:\gamma_0.
	\]
	Thus
	\begin{equation}\label{eq: low-dim glm population curvature}
	\alpha_0I_d\preceq H_{L_\psi}(\bbeta)\preceq\gamma_0I_d
	\qquad\text{for every }\bbeta\in\Theta_0
	\text{ and }L_\psi\geq L_{\psi,0}.
	\end{equation}
	The lower bound uses only the central region on which the saturated and
	original scores coincide; it does not impose a dimension-free
	samplewise bound on the linear predictors.

	\textit{Step 2: concentration at a fixed parameter.}
	For fixed $\bbeta\in\Theta_0$ and unit $\bm u$, put
	\[
	Z_i(\bm u,\bbeta)
	=\psi_{L_\psi}''(\bm x_i^\top\bbeta)(\bm u^\top\bm x_i)^2.
	\]
	Because $0\leq\psi_{L_\psi}''\leq c_2$ and
	$\bm u^\top\bm x_i$ has sub-Gaussian norm at most $K_{\bm x}$, the centered variable
	$Z_i(\bm u,\bbeta)-\E Z_i(\bm u,\bbeta)$ has sub-exponential norm at most $a_3c_2K_{\bm x}^2$, where $a_3$ is universal.  Bernstein's inequality therefore yields
	\begin{align}\label{eq: fixed glm hessian concentration}
	&\Pro\!\left(\left|\bm u^\top\{H_{n,L_\psi}(\bbeta)-H_{L_\psi}(\bbeta)\}\bm u\right|>q\right)\notag\\
	&\qquad\leq2\exp\!\left[-a_4n
	\min\left\{\frac{q^2}{c_2^2K_{\bm x}^4},
	\frac{q}{c_2K_{\bm x}^2}\right\}\right]
	\end{align}
	for a universal $a_4>0$.

	\textit{Step 3: uniformity over directions and parameters.}
	A standard $1/2$-net argument for the Euclidean sphere and the sub-Gaussian tail in (D) show that, for a sufficiently large fixed $A$,
	\begin{equation}\label{eq: glm design norm event}
	\mathcal E_R:=\left\{\max_{1\leq i\leq n}\|\bm x_i\|_2\leq
	R_n\right\},
	\qquad R_n=AK_{\bm x}\sqrt{d+\log n},
	\end{equation}
	has probability at least $1-n^{-3}$.  The same argument, followed by integration of the tail, gives
	\begin{equation}\label{eq: glm design third moment}
	\E\|\bm x\|_2^3\leq a_5K_{\bm x}^3d^{3/2}
	\end{equation}
	for a universal $a_5$.

	Let $\mathcal N_u$ be a $1/4$-net of the unit sphere with $|\mathcal N_u|\leq9^d$.  Set
	\[
	\widetilde c_3=c_3+c_2\|\chi'\|_\infty,
	\qquad
	L_n=\widetilde c_3\{R_n^3+a_5K_{\bm x}^3d^{3/2}\},
	\qquad
	r_n=\min\left\{b_0,\frac{\alpha_0}{8L_n}\right\},
	\]
	with $r_n=b_0$ when $L_n=0$, and let $\mathcal N_\beta$ be an $r_n$-net of $\Theta_0$ satisfying
	$|\mathcal N_\beta|\leq(1+2b_0/r_n)^d$.
	Apply \eqref{eq: fixed glm hessian concentration} with $q=\alpha_0/16$ and take a union bound over
	$\mathcal N_u\times\mathcal N_\beta$.  Since, for any symmetric matrix $A$,
	\[
	\|A\|_{\rm op}\leq2\max_{\bm u\in\mathcal N_u}|\bm u^\top A\bm u|,
	\]
	we obtain
	\begin{equation}\label{eq: glm hessian net event}
	\max_{\bbeta\in\mathcal N_\beta}
	\|H_{n,L_\psi}(\bbeta)-H_{L_\psi}(\bbeta)\|_{\rm op}
	\leq\frac{\alpha_0}{8}
	\end{equation}
	except on an event of probability at most
	\[
	2\exp\left\{d\log9+d\log(1+2b_0/r_n)-a_6n\right\},
	\]
	where $a_6>0$ depends only on the fixed constants in (D), (G), the
	fixed cutoff $\chi$, and $b_0$.

	It remains to pass from the parameter net to all of $\Theta_0$.
	By \eqref{eq: saturated link score}, for $L_\psi\geq1$,
	\[
		\left|\frac{\d}{\d t}\psi_{L_\psi}''(t)\right|
		=\left|\psi'''(t)\chi(t/L_\psi)
		+\frac{\psi''(t)}{L_\psi}\chi'(t/L_\psi)\right|
		\leq\widetilde c_3.
	\]
	The mean-value theorem therefore implies
	\begin{align*}
	\|H_{n,L_\psi}(\bbeta)-H_{n,L_\psi}(\widetilde\bbeta)\|_{\rm op}
	&\leq \widetilde c_3\|\bbeta-\widetilde\bbeta\|_2
	\frac1n\sum_{i=1}^n\|\bm x_i\|_2^3,\\
	\|H_{L_\psi}(\bbeta)-H_{L_\psi}(\widetilde\bbeta)\|_{\rm op}
	&\leq \widetilde c_3\|\bbeta-\widetilde\bbeta\|_2\E\|\bm x\|_2^3.
	\end{align*}
	On $\mathcal E_R$, choosing $\widetilde\bbeta\in\mathcal N_\beta$ with
	$\|\bbeta-\widetilde\bbeta\|_2\leq r_n$ and using \eqref{eq: glm design third moment} shows that the sum of these two differences is at most
	$L_nr_n\leq\alpha_0/8$.  Consequently, \eqref{eq: glm hessian net event} gives
	\begin{equation}\label{eq: uniform glm hessian concentration}
	\sup_{\bbeta\in\Theta_0}
	\|H_{n,L_\psi}(\bbeta)-H_{L_\psi}(\bbeta)\|_{\rm op}
	\leq\frac{\alpha_0}{4}.
	\end{equation}

	Finally, if $n\geq C_0d\log(en)$ with $C_0$ sufficiently large, then $d\leq n$ and
	$\log(1+2b_0/r_n)\lesssim\log(en)$.  The exceptional probability in \eqref{eq: glm hessian net event} is therefore at most $n^{-3}$.  Combining it with \eqref{eq: glm design norm event} shows that \eqref{eq: uniform glm hessian concentration} holds with probability at least $1-2n^{-3}\geq1-n^{-2}$.  Together with \eqref{eq: low-dim glm population curvature}, this proves \eqref{eq: low-dim glm empirical curvature}, for example with
	\[
	\alpha=\frac{\alpha_0}{2},
	\qquad
	\gamma=\gamma_0+\frac{\alpha_0}{2}.
	\]
\end{proof}

\subsection{Proof of Proposition \ref{prop: non-sparse glm privacy}}\label{sec: proof of prop: non-sparse glm privacy}

\begin{proof}[Proof of Proposition \ref{prop: non-sparse glm privacy}]
	Fix deterministic levels $L_{\bm x},L_\psi,$ and $R$. For one datum
	$z=(y,\bm x)$, write $\overline{\bm x}=\mathsf C_{L_{\bm x}}^{(2)}(\bm x)$ and define its clipped gradient contribution at $\bbeta$ by
	\begin{align*}
		\bm h_{\bbeta}(z)
		=\{\psi_{L_\psi}'(\overline{\bm x}^{\top}\bbeta)-\Pi_R(y)\}
		\overline{\bm x}.
	\end{align*}
	The clipping map is defined for every $\bm x\in\mathbb R^d$ and satisfies
	$\|\overline{\bm x}\|_2\leq L_{\bm x}$. Equation
	\eqref{eq: saturated score bound} and $|\Pi_R(y)|\leq R$ therefore
	give the deterministic bound
	\begin{align}
		\|\bm h_{\bbeta}(z)\|_2
		&\leq
		\{|\psi_{L_\psi}'(\overline{\bm x}^{\top}\bbeta)|+|\Pi_R(y)|\}
		\|\overline{\bm x}\|_2\notag\\
		&\leq(R+S_\psi)L_{\bm x}.                    \label{eq: clipped low dim contribution bound}
	\end{align}
	This inequality is uniform in $z$ and $\bbeta$; no sampling event is being used.

	Condition on all private releases before iteration $t$.  The value of
	$\bbeta^t$ is then fixed.  Let $q_t(\bm Z)$ denote the unprojected update before adding $\bm w_t$.  If replacement-adjacent data sets $\bm Z$ and $\bm Z'$ differ by replacing $z$ with $z'$, all common summands cancel, and \eqref{eq: clipped low dim contribution bound} yields
	\begin{align*}
	\|q_t(\bm Z)-q_t(\bm Z')\|_2
	&=\frac{\eta^0}{n}
	\|\bm h_{\bbeta^t}(z)-\bm h_{\bbeta^t}(z')\|_2\\
	&\leq\frac{\eta^0}{n}
	\{\|\bm h_{\bbeta^t}(z)\|_2+
	\|\bm h_{\bbeta^t}(z')\|_2\}\\
	&\leq\frac{2\eta^0(R+S_\psi)L_{\bm x}}{n}
	\leq\frac{\eta^0B}{n},
	\end{align*}
	where the last inequality uses $B=4(R+S_\psi)L_{\bm x}$.
	Algorithm \ref{algo: private glm} adds independent Gaussian coordinates with variance
	\begin{align*}
	(\eta^0)^2\frac{2B^2\log(2T/\delta)}
	{n^2(\varepsilon/T)^2}.
	\end{align*}
	The Gaussian mechanism, Example \ref{fc: laplace and gaussian mechanisms}, therefore makes this conditional release $(\varepsilon/T,\delta/T)$-differentially private.  Projection onto $\Theta_0$ is post-processing.  Since the calculation holds for every realization of the previous releases, adaptive basic composition over $T$ iterations proves $(\varepsilon,\delta)$-differential privacy.
\end{proof}

\subsection{Proof of Theorem \ref{thm: non-sparse glm upper bound}} \label{sec: proof of thm: non-sparse glm upper bound}

{Let
\[
\overline{\bm x}_i=\mathsf C_{L_{\bm x}}^{(2)}(\bm x_i),
\qquad
\overline\L_{n,L_\psi}^R(\bbeta)=\frac1n\sum_{i=1}^n
\{\psi_{L_\psi}(\overline{\bm x}_i^\top\bbeta)
-\Pi_R(y_i)\overline{\bm x}_i^\top\bbeta\}
\]
be the loss actually used by Algorithm \ref{algo: private glm}.  Also
write
\[
\L_{n,L_\psi}(\bbeta)=\frac1n\sum_{i=1}^n
\{\psi_{L_\psi}(\bm x_i^\top\bbeta)-y_i\bm x_i^\top\bbeta\}
\]
for the saturated-score loss with original covariates and responses,
and write \(\L_n\) for the original empirical negative log-likelihood.
Choose
\[
\eta^0=\frac{2}{\alpha+\gamma},\qquad
\rho=\max\left\{\frac12,\frac{\gamma-\alpha}{\gamma+\alpha}\right\}<1,\qquad
L_\psi=L_{\psi,0}+C_\psi K_{\bm x}b_0\sqrt{\log n},
\]
and, for sufficiently large absolute constants \(C_{\bm x},C_\psi\), set
\[
L_{\bm x}=C_{\bm x}K_{\bm x}\sqrt{d+\log n},
\qquad
S_\psi=|\psi'(0)|+2c_2L_\psi,
\]
and
\[
R=|\psi'(0)|+c_2L_\psi+\sqrt{6c_2c(\sigma)\log n},
\qquad B=4(R+S_\psi)L_{\bm x}.
\]
Finally, take
\[
T=\left\lceil\frac{\log(nb_0\vee n)}{-\log\rho}\right\rceil.
\]
Thus $T\asymp\log n$ and $\rho^{2T}b_0^2\leq n^{-2}$.

\begin{proof}[Proof of Theorem \ref{thm: non-sparse glm upper bound}]
	Proposition \ref{prop: non-sparse glm privacy}, applied with the
	deterministic value of \(L_{\bm x}\) above, proves the asserted privacy
	guarantee.  This part of the argument holds for every data set.  We now
	control the risk under the sampling model.

	Define the data events
	\begin{align*}
	\mathcal E_{\bm x}
	&=\{\|\bm x_i\|_2\leq L_{\bm x},\ 1\leq i\leq n\},\\
	\mathcal E_{\rm curv}
	&=\left\{\alpha I_d\preceq\nabla^2\L_{n,L_\psi}(\bbeta)
		\preceq\gamma I_d\ \text{for every }\bbeta\in\Theta_0\right\},\\
	\mathcal E_{\rm pred}
	&=\{|\bm x_i^\top\bbeta^*|\leq L_\psi,\ 1\leq i\leq n\},\\
	\mathcal E_{\rm tr}
	&=\{\Pi_R(y_i)=y_i,\ 1\leq i\leq n\}.
	\end{align*}
	We first verify the probability of \(\mathcal E_{\bm x}\) without
	appealing to an almost-sure design bound.  Let \(\mathcal N\) be a
	\(1/2\)-net of the Euclidean unit sphere with
	\(|\mathcal N|\leq5^d\).  The usual net inequality gives, for every
	\(\bm x\in\mathbb R^d\),
	\[
	\|\bm x\|_2\leq2\max_{\bm u\in\mathcal N}|\bm u^\top\bm x|.
	\]
	For \(\bm x=\bm0\) the claim is immediate.  Otherwise, if
	\(\bm v=\bm x/\|\bm x\|_2\) and
	\(\|\bm u-\bm v\|_2\leq1/2\), then
	\(\|\bm x\|_2=\bm v^\top\bm x
	\leq\bm u^\top\bm x+\|\bm u-\bm v\|_2\|\bm x\|_2\).
	Assumption (D) and a union bound over the net therefore imply
	\begin{align*}
	\Pr\{\|\bm x\|_2>L_{\bm x}\}
	&\leq\sum_{\bm u\in\mathcal N}
	\Pr\{|\bm u^\top\bm x|>L_{\bm x}/2\}\\
	&\leq2\cdot5^d
	\exp\left\{-\frac{cL_{\bm x}^2}{4K_{\bm x}^2}\right\}.
	\end{align*}
	Choosing \(C_{\bm x}\) sufficiently large makes the last display at
	most \(n^{-5}\), uniformly for \(d\geq1\) and \(n\geq2\).  A second
	union bound, now over the \(n\) observations, gives
	\begin{equation}
	\Pr(\mathcal E_{\bm x}^c)\leq n^{-4}.             \label{eq: low dim no covariate clipping}
	\end{equation}
	On \(\mathcal E_{\bm x}\), the definition of Euclidean clipping gives
	\(\overline{\bm x}_i=\bm x_i\) for every \(i\). Consequently,
	\begin{equation}
	\overline\L_{n,L_\psi}^R(\bbeta)=\L_{n,L_\psi}^R(\bbeta)
	\quad\text{for every }\bbeta\in\Theta_0,          \label{eq: low dim clipped loss coupling}
	\end{equation}
	where \(\L_{n,L_\psi}^R\) denotes the saturated-score,
	response-truncated loss with the original covariates.

	The Hessian of \(\L_{n,L_\psi}\) does not involve the responses.
	Proposition
	\ref{prop: low-dim glm empirical curvature} therefore gives
	\[
	\sup_{\bbeta^*\in\Theta_0}
	\Pro_{\bbeta^*}(\mathcal E_{\rm curv}^c)\leq n^{-2}.
	\]
	Assumption (D) also gives, uniformly over
	$\bbeta^*\in\Theta_0$,
	\begin{align*}
	\Pr(\mathcal E_{\rm pred}^c)
	&\leq2n\exp\left\{
	-\frac{cL_\psi^2}{K_{\bm x}^2b_0^2}\right\}
	\leq n^{-3},
	\end{align*}
	after choosing $C_\psi$ sufficiently large; the assertion is
	immediate when $\bbeta^*=\bm0$.

	Moreover, \eqref{eq: glm mgf upper bound} and (G) imply,
	conditionally on $\bm x_i$,
	\[
	\Pro_{\bbeta^*}\!\left(
	|y_i-\psi'(\bm x_i^\top\bbeta^*)|>u\mid\bm x_i\right)
	\leq2\exp\!\left\{-\frac{u^2}{2c_2c(\sigma)}\right\}.
	\]
	On $\mathcal E_{\rm pred}$, the mean-value theorem and (G) give
	\[
	|\psi'(\bm x_i^\top\bbeta^*)|
	\leq|\psi'(0)|+c_2L_\psi.
	\]
	Our choice of $R$ and a union bound over the residual tails therefore
	give
	\[
	\sup_{\bbeta^*\in\Theta_0}
	\Pr_{\bbeta^*}(\mathcal E_{\rm pred}^c\cup
	\mathcal E_{\rm tr}^c)\leq n^{-3}+2n^{-2}.
	\]
	On $\mathcal E_{\rm pred}$,
	$\psi_{L_\psi}'(\bm x_i^\top\bbeta^*)
	=\psi'(\bm x_i^\top\bbeta^*)$ for every $i$. Thus, on
	\[
	\mathcal E=\mathcal E_{\bm x}\cap\mathcal E_{\rm curv}
	\cap\mathcal E_{\rm pred}\cap\mathcal E_{\rm tr},
	\]
	the loss used by the algorithm is uniformly strongly convex and
	smooth on $\Theta_0$, and its gradient at $\bbeta^*$ equals the
	original empirical GLM score. The preceding probability bounds give
	\begin{equation}
	\sup_{\bbeta^*\in\Theta_0}\Pr_{\bbeta^*}(\mathcal E^c)
	\leq3n^{-2}+n^{-3}+n^{-4}.                           \label{eq: low dim good event probability}
	\end{equation}

	Fix a data set in $\mathcal E$ and put
	$e_t=\bbeta^t-\bbeta^*$ and $g_*=\nabla\L_n(\bbeta^*)$.
	Because $\Theta_0$ is convex and every iterate is projected onto it, the integrated Hessian
	\[
	H_t=\int_0^1\nabla^2\L_{n,L_\psi}
	\{\bbeta^*+u(\bbeta^t-\bbeta^*)\}\,\d u
	\]
	has all eigenvalues in $[\alpha,\gamma]$.  Hence
	$\|I_d-\eta^0H_t\|_{\rm op}\leq\rho$.  Non-expansiveness of Euclidean projection and $\bbeta^*\in\Theta_0$ now give the one-step inequality
	\begin{align}\label{eq: truth centered glm contraction}
	\|e_{t+1}\|_2
	&\leq\|(I_d-\eta^0H_t)e_t-\eta^0g_*+\bm w_t\|_2\\
	&\leq\rho\|e_t\|_2+\eta^0\|g_*\|_2+\|\bm w_t\|_2.\notag
	\end{align}
		Iterating \eqref{eq: truth centered glm contraction} first gives
		\[
		\|e_T\|_2\leq \rho^T\|e_0\|_2
		+\eta^0\|g_*\|_2\sum_{j=0}^{T-1}\rho^j
		+\sum_{k=0}^{T-1}\rho^{T-1-k}\|\bm w_k\|_2.
		\]
		Using $\|e_0\|_2\leq b_0$, $\sum_{j=0}^{T-1}\rho^j\leq(1-\rho)^{-1}$, and $(a+b+c)^2\leq3(a^2+b^2+c^2)$ yields
	\begin{align}\label{eq: glm expected contraction}
	\E_{\bm w}\|e_T\|_2^2
	&\leq 3\rho^{2T}b_0^2+
	\frac{3(\eta^0)^2}{(1-\rho)^2}\|g_*\|_2^2
	+\frac{3}{(1-\rho)^2}\E\|\bm w_0\|_2^2.
	\end{align}
	Indeed,
	\[
	\E\left(\sum_{k=0}^{T-1}\rho^{T-1-k}\|\bm w_k\|_2\right)^2
	\leq\left(\sum_{j\geq0}\rho^j\right)
	\left(\sum_{j\geq0}\rho^j\E\|\bm w_0\|_2^2\right)
	=\frac{\E\|\bm w_0\|_2^2}{(1-\rho)^2}.
	\]
	This calculation keeps the privacy noise from every iteration and does not invoke a noiseless descent inequality for the noisy iterate.

	The empirical score has mean zero and its cross terms vanish.  Since
	$\operatorname{var}(y\mid\bm x)=c(\sigma)\psi''(\bm x^\top\bbeta^*)$, assumptions (D) and (G) imply
	\begin{align}\label{eq: glm score second moment}
	\E_{\bbeta^*}\|g_*\|_2^2
	&=\frac1n\E_{\bbeta^*}\!\left[
	\{y-\psi'(\bm x^\top\bbeta^*)\}^2\|\bm x\|_2^2\right]\\
	&\leq \frac{c(\sigma)c_2}{n}\E\|\bm x\|_2^2
	\lesssim c(\sigma)\frac dn.\notag
	\end{align}
	For one privacy-noise vector, the variance specified in Algorithm \ref{algo: private glm} gives
	\begin{align}\label{eq: glm privacy noise second moment}
	\E\|\bm w_0\|_2^2
	&=2(\eta^0)^2B^2\frac{dT^2\log(2T/\delta)}
	{n^2\varepsilon^2}\\
	&\lesssim c(\sigma)
	\frac{d(d+\log n)\log(1/\delta)\log^3 n}
	{n^2\varepsilon^2}.\notag
	\end{align}
	Here the first equality follows by summing the common coordinate
	variance over \(d\) coordinates.  Also,
	\[
	B^2=16(R+S_\psi)^2L_{\bm x}^2
	\lesssim c(\sigma)(d+\log n)\log n,
	\]
	because $L_\psi\lesssim\sqrt{\log n}$ and all model constants,
	including $\psi'(0)$, are fixed. Fixed additive constants are absorbed
	using the stipulated fixed upper and lower bounds on \(c(\sigma)\). We
	also used
	$T\lesssim\log n$ and
	$\log(2T/\delta)\lesssim\log(1/\delta)$ for $\delta\leq n^{-1}$.

	Finally, both the algorithmic output and \(\bbeta^*\) belong to
	\(\Theta_0\), so projection gives the deterministic bound
	$\|\bbeta^T-\bbeta^*\|_2^2\leq4b_0^2$.  Integrating
	\eqref{eq: glm expected contraction}, applying
	\eqref{eq: glm score second moment}--\eqref{eq: glm privacy noise second moment}, and using \eqref{eq: low dim good event probability} gives
	\begin{align*}
	\sup_{\bbeta^*\in\Theta_0}\E_{\bbeta^*}
	\|\bbeta^T-\bbeta^*\|_2^2
	&\lesssim n^{-2}+c(\sigma)\frac dn\\
	&\quad+c(\sigma)
	\frac{d(d+\log n)\log(1/\delta)\log^3 n}
	{n^2\varepsilon^2}.
	\end{align*}
	Since \(d\geq1\) and \(n\geq2\), the term \(n^{-2}\) is absorbed
	by \(d/n\).  This proves \eqref{eq: non-sparse glm upper bound}.
\end{proof}}

	\section{Omitted Proofs in Section \ref{sec: BTL}}\label{sec: BTL proofs}

{
\subsection{Proof of Proposition \ref{prop: ranking attack soundness}}\label{sec: proof of prop: ranking attack soundness}
\begin{proof}[Proof of Proposition \ref{prop: ranking attack soundness}]
By the range condition in the proposition, $M$ takes values in $\Theta$. In particular, each induced coordinate $\widehat\eta_k(M)=(M_{2k-1}-M_{2k})/2$ belongs to $[-1,1]$.

Fix the public graph $\mathcal G$ and an observed edge $e=(i,j)$. Let $\bm Y^{(e)}$ be obtained by replacing $Y_{ij}$ with an independent copy $Y'_{ij}$ having the same conditional BTL distribution, and define
\[
A'_e=\left\langle\widehat{\bm\eta}(M(\bm Y^{(e)},\mathcal G))-\bm\eta,
\{Y_{ij}-q_{ij}(\bm\eta)\}\bm V^\top\bm x_{ij}\right\rangle.
\]
The output $M(\bm Y^{(e)},\mathcal G)$ is independent of $Y_{ij}$, so $\E(A'_e\mid\mathcal G)=0$. The vector $\bm V^\top\bm x_{ij}$ has either one nonzero coordinate of magnitude two or two nonzero coordinates of magnitude one. Hence $|A_e|\leq4$ and $|A'_e|\leq4$. Applying the privacy comparison in \eqref{eq: score attack upper bound intermediate step} with $T=4$ gives
\begin{align}\label{eq: ranking attack soundness privacy bound}
\E(A_e\mid\mathcal G)\leq2\varepsilon\E(|A'_e|\mid\mathcal G)+8\delta.
\end{align}

Let $m_{\mathcal G}=|\mathcal G|$ and
\[
R_{\mathcal G}(\bm\eta)=\E_{\bm Y\mid\bm\eta,\mathcal G}
\|M(\bm Y,\mathcal G)-\bm V\bm\eta\|_2^2.
\]
Since replacing one coordinate by an independent copy preserves the joint distribution of the data, the random vector
$\widehat{\bm\eta}(M(\bm Y^{(e)},\mathcal G))-\bm\eta$ has the same marginal distribution as
$\widehat{\bm\eta}(M(\bm Y,\mathcal G))-\bm\eta$. Conditional Cauchy--Schwarz therefore gives
\begin{align*}
\sum_{e\in\mathcal G}\E(|A'_e|\mid\mathcal G)
&\leq\sqrt{m_{\mathcal G}\sum_{e\in\mathcal G}
\E\left[\left\{(\widehat{\bm\eta}-\bm\eta)^\top\bm V^\top\bm x_e\right\}^2\middle|\mathcal G\right]}\\
&=\sqrt{m_{\mathcal G}\E\left[(\widehat{\bm\eta}-\bm\eta)^\top
\bm V^\top\bm L_{\mathcal G}\bm V(\widehat{\bm\eta}-\bm\eta)\middle|\mathcal G\right]}\\
&\leq\sqrt{2m_{\mathcal G}\lambda_{\max}(\bm L_{\mathcal G})
\E(\|\widehat{\bm\eta}-\bm\eta\|_2^2\mid\mathcal G)}\\
&\leq\sqrt{m_{\mathcal G}\lambda_{\max}(\bm L_{\mathcal G})R_{\mathcal G}(\bm\eta)}.
\end{align*}
The last two inequalities use $\|\bm V\|_{\rm op}^2=2$ and
\begin{align}\label{eq: ranking induced risk}
\|\widehat{\bm\eta}-\bm\eta\|_2^2
=\frac14\|\bm V^\top\{M-\bm V\bm\eta\}\|_2^2
\leq\frac12\|M-\bm V\bm\eta\|_2^2.
\end{align}

Let $\mathcal H$ be the graph event
\[
m_{\mathcal G}\leq n^2p,\qquad \lambda_{\max}(\bm L_{\mathcal G})\leq4np.
\]
Standard Chernoff bounds for the number of edges and the maximum degree, together with
$\lambda_{\max}(\bm L_{\mathcal G})\leq2\max_i\deg_{\mathcal G}(i)$, show that
$\Pro(\mathcal H^c)\leq n^{-3}$ when $p\geq C_0\log n/n$ and $C_0$ is sufficiently large. On $\mathcal H$, the preceding display is at most
$2n^{3/2}p\sqrt{R_{\mathcal G}(\bm\eta)}$. Averaging over $\mathcal G$ and applying Cauchy--Schwarz once more gives
\[
\E_{\mathcal G}\left[\1_{\mathcal H}\sum_{e\in\mathcal G}\E(|A'_e|\mid\mathcal G)\right]
\leq2n^{3/2}p\sqrt{R_{\bth(\bm\eta)}(M)}.
\]
On $\mathcal H^c$, the deterministic bound $|A'_e|\leq4$ gives a contribution at most $2n^2\Pro(\mathcal H^c)\leq2/n$. Finally,
$\E m_{\mathcal G}=\binom n2p$. Summing \eqref{eq: ranking attack soundness privacy bound} over the observed edges and averaging over the graph proves \eqref{eq: ranking attack soundness} after adjusting constants.
\end{proof}

\subsection{Proof of Proposition \ref{prop: ranking attack completeness}}\label{sec: proof of prop: ranking attack completeness}
\begin{proof}[Proof of Proposition \ref{prop: ranking attack completeness}]
For a fixed public graph, let $f_{\bm\eta}(\bm Y\mid\mathcal G)$ be the likelihood under $\bth=\bm V\bm\eta$. By construction,
\[
\nabla_{\bm\eta}\log f_{\bm\eta}(\bm Y\mid\mathcal G)
=\sum_{(i,j)\in\mathcal G}\{Y_{ij}-q_{ij}(\bm\eta)\}\bm V^\top\bm x_{ij}.
\]
Consequently,
\begin{align*}
\E_{\bm Y\mid\bm\eta,\mathcal G}\sum_{i<j}A_{ij}
&=\E\left\langle\widehat{\bm\eta}(M)-\bm\eta,
\nabla_{\bm\eta}\log f_{\bm\eta}(\bm Y\mid\mathcal G)\right\rangle\\
&=\sum_{k=1}^{d_0}\frac{\partial}{\partial\eta_k}
\E_{\bm Y\mid\bm\eta,\mathcal G}\widehat\eta_k(M).
\end{align*}
The finite Bernoulli likelihood is smooth in $\bm\eta$, and $M$ is bounded, so differentiation may be passed through the finite sum defining the conditional expectation. The last equality therefore follows from the score identity; the term involving the nonrandom vector $\bm\eta$ has expectation zero. Since the graph distribution does not depend on $\bm\eta$, the same identity holds after averaging over $\mathcal G$.

Let $g(\bm\eta)=\E_{\mathcal G,\bm Y\mid\bm\eta}\widehat{\bm\eta}(M)$ and let the coordinates of $\bm\eta$ be independent with density
$\pi(t)=\1(|t|<1)(15/16)(1-t^2)^2$. Because $M$ is $\Theta$-valued, $g$ is bounded; hence $g_k(t,\bm\eta_{-k})\pi(t)\to0$ at both endpoints. Proposition \ref{prop: score attack stein's lemma}, applied in dimension $d_0$, therefore yields
\begin{align*}
\E_{\bm\eta,\mathcal G,\bm Y}\sum_{i<j}A_{ij}
&\geq\sum_{k=1}^{d_0}\E\left[-\eta_k\frac{\pi'(\eta_k)}{\pi(\eta_k)}\right]
-\sqrt{\E\|\widehat{\bm\eta}(M)-\bm\eta\|_2^2}
\sqrt{\sum_{k=1}^{d_0}\E\left\{\frac{\pi'(\eta_k)}{\pi(\eta_k)}\right\}^2}\\
&=d_0-\sqrt{10d_0\,\E\|\widehat{\bm\eta}(M)-\bm\eta\|_2^2}.
\end{align*}
Here direct integration gives
$\E[-\eta\pi'(\eta)/\pi(\eta)]=1$ and
$\E\{\pi'(\eta)/\pi(\eta)\}^2=10$. By \eqref{eq: ranking induced risk}, the assumed risk bound implies
$\E\|\widehat{\bm\eta}(M)-\bm\eta\|_2^2\leq c_0n/2$. Since $d_0\asymp n$, choosing $c_0$ sufficiently small makes the last display at least $cn$.
\end{proof}

\subsection{Proof of Theorem \ref{thm: ranking lower bound}}\label{sec: proof of thm: ranking lower bound}
\begin{proof}[Proof of Theorem \ref{thm: ranking lower bound}]
The non-private term $1/p$ follows from the standard BTL minimax lower bound \cite{negahban2017rank, shah2015estimation} on the centered parameter space.

For the privacy term, let $M$ be any edge-level $(\varepsilon,\delta)$-differentially private estimator \furtherrev{with values in the ambient space $\mathbb R^n$}. Set
$\widetilde M=\Pi_\Theta\circ M$, then $\widetilde M$ remains edge-level $(\varepsilon,\delta)$-differentially private and
\[
R_\bth(\widetilde M)\leq R_\bth(M),\qquad \bth\in\Theta.
\]
It therefore suffices to establish the lower bound for the $\Theta$-valued estimator $\widetilde M$, which we relabel as $M$. \furtherrev{This post-processing reduction is precisely what permits the range conditions in the two score-attack propositions while the minimax infimum still ranges over all $\mathbb R^n$-valued private estimators.} If
$\sup_{\bth\in\Theta}R_\bth(M)\leq c_0n$, average Propositions \ref{prop: ranking attack soundness} and \ref{prop: ranking attack completeness} over the paired prior to obtain
\[
cn\leq C\varepsilon n^{3/2}p
\sqrt{\E_{\bm\eta}R_{\bth(\bm\eta)}(M)}+Cn^2p\delta+1.
\]
The assumption $\delta\leq c_1/(np)$, with $c_1$ sufficiently small, absorbs the second and third terms into the left side. Hence
\[
\sup_{\bth\in\Theta}R_\bth(M)
\geq\E_{\bm\eta}R_{\bth(\bm\eta)}(M)
\gtrsim\frac{1}{np^2\varepsilon^2}.
\]
If the assumed $c_0n$ risk bound fails, the risk is already of order $n$. These two cases give
\[
\sup_{\bth\in\Theta}R_\bth(M)
\gtrsim\min\left\{n,\frac{1}{np^2\varepsilon^2}\right\}.
\]
Taking the maximum with the non-private lower bound, and using that the maximum of two nonnegative quantities is at least half their sum, proves \eqref{eq: ranking lower bound}.
\end{proof}

\subsection{Proof of Proposition \ref{prop: ranking MLE privacy}}\label{sec: proof of prop: ranking MLE privacy}
\begin{proof}[Proof of Proposition \ref{prop: ranking MLE privacy}]
The graph is public and is identical in adjacent data sets, so the spectral branch in \eqref{eq: ranking MLE definition} is data independent. If
$\lambda_2(\bm L_{\mathcal G})<e^{-1}np$, the estimator outputs $\bm0$ and is trivially private. Suppose otherwise.

For every $\bth\in\bm1^\perp$,
\[
\nabla^2\bar{\mathcal L}(\bth;y)
=\sum_{(i,j)\in\mathcal G}
\bar\phi''(\bm x_{ij}^\top\bth)
\bm x_{ij}\bm x_{ij}^\top
\succeq\kappa_0\bm L_{\mathcal G}.
\]
Thus the extended loss is
$\mu_{\mathcal G}$-strongly convex and coercive on $\bm1^\perp$, where
$\mu_{\mathcal G}=\kappa_0\lambda_2(\bm L_{\mathcal G})$.
For every $\bm b\in\bm1^\perp$, the function
$\bth\mapsto\bar{\mathcal L}(\bth;y)+\langle\bm b,\bth\rangle$ is therefore coercive and strongly convex, so it has a unique minimizer satisfying
$-\nabla\bar{\mathcal L}(\bth;y)=\bm b$. This proves surjectivity of the map below. Strong monotonicity of $\nabla\bar{\mathcal L}$ (equivalently, strong anti-monotonicity of the map below) proves injectivity and gives a Lipschitz inverse. Consequently, for each data set $y$, the map
\[
\bm T_y:\bm1^\perp\to\bm1^\perp,\qquad
\bm T_y(\bth)=-\nabla\bar{\mathcal L}(\bth;y),
\]
is a homeomorphism, and $\widetilde\bth\in O$ if and only if
$\bm B\in\bm T_y(O)$ for every Borel set $O\subseteq\bm1^\perp$.

We first record the translation inequality obeyed by the projected Laplace noise.
For any measurable $S\subseteq\bm1^\perp$ and any $\bm s\in\bm1^\perp$, the product Laplace density of $\bm W$ gives
\begin{align}\label{eq: projected Laplace translation}
\Pro(\bm B\in S)
&\leq\exp\left(\frac{\varepsilon}{2}\|\bm s\|_1\right)
\Pro(\bm B\in S+\bm s).
\end{align}
Indeed, apply the ordinary product-Laplace density comparison to
$\{\bm w:\bm P\bm w\in S\}$ and its translate by $\bm s$; the identity
$\bm P\bm s=\bm s$ turns the latter set into
$\{\bm w:\bm P\bm w\in S+\bm s\}$.

If $y,y'$ differ only in the outcome on edge $e=(i,j)$, then
\[
\bm T_{y'}(\bth)=\bm T_y(\bth)+\bm s,
\qquad
\bm s=(y'_{ij}-y_{ij})\bm x_{ij},
\qquad
\|\bm s\|_1\leq2.
\]
Taking $S=\bm T_y(O)$ in \eqref{eq: projected Laplace translation} therefore yields
\[
\Pro_y(\widetilde\bth\in O)
\leq e^\varepsilon\Pro_{y'}(\widetilde\bth\in O).
\]
The same argument with $y$ and $y'$ interchanged proves pure
edge-level differential privacy. The final projection onto $\Theta$ preserves privacy by post-processing.
\end{proof}

\subsection{Proof of Proposition \ref{prop: ranking MLE accuracy}}\label{sec: proof of prop: ranking MLE accuracy}
\begin{proof}[Proof of Proposition \ref{prop: ranking MLE accuracy}]
Let
\[
\mathcal H_0=\left\{\lambda_2(\bm L_{\mathcal G})\geq e^{-1}np,
\ m_{\mathcal G}\leq n^2p\right\}.
\]
To verify the graph bound, note that
$\E\bm L_{\mathcal G}=np\bm P$. Matrix Bernstein's inequality
\cite{tropp2015introduction}, applied to
$\sum_{i<j}(\1\{(i,j)\in\mathcal G\}-p)\bm x_{ij}\bm x_{ij}^\top$,
gives
\[
\|\bm L_{\mathcal G}-np\bm P\|_{\rm op}
\lesssim\sqrt{np\log n}+\log n
\]
with probability at least $1-n^{-2}$. A Chernoff bound gives
$m_{\mathcal G}\leq n^2p$ with the same probability order. Hence, after increasing
$C_1$ if needed, $\Pro(\mathcal H_0^c)\leq n^{-2}$ whenever
$p\geq C_1\log n/n$.

On $\mathcal H_0$, put $\mu_{\mathcal G}=\kappa_0\lambda_2(\bm L_{\mathcal G})$.
Because every true comparison margin satisfies
$|\bm x_{ij}^\top\bth|\leq2$, the extended and ordinary logistic losses have the same gradient at the true $\bth$. In particular,
$\nabla\bar{\mathcal L}(\bth;Y)$ is the ordinary BTL score.
The first-order condition for \eqref{eq: ranking MLE definition} and strong convexity give
\[
\mu_{\mathcal G}\|\widetilde\bth-\bth\|_2^2
\leq
\left\langle-\bm B-\nabla\bar{\mathcal L}(\bth;Y),
\widetilde\bth-\bth\right\rangle,
\qquad
\|\widetilde\bth-\bth\|_2
\leq\frac{\|\bm B+\nabla\bar{\mathcal L}(\bth;Y)\|_2}{\mu_{\mathcal G}}.
\]
Conditional on $\mathcal G$, the score has mean zero and independent edge contributions. Since
$\Var(Y_{ij}\mid\mathcal G)\leq1/4$ and $\|\bm x_{ij}\|_2^2=2$,
\[
\E\left(\|\nabla\bar{\mathcal L}(\bth;Y)\|_2^2\mid\mathcal G\right)
\leq\frac{m_{\mathcal G}}2.
\]
The perturbation is independent and centered. Since the coordinates of $\bm W$
have variance $8/\varepsilon^2$ and $\bm P$ is an orthogonal projector,
\[
\E\|\bm B\|_2^2
=\E(\bm W^\top\bm P\bm W)
=\frac{8(n-1)}{\varepsilon^2}.
\]
It follows that, on $\mathcal H_0$,
\[
\E(\|\widetilde\bth-\bth\|_2^2\mid\mathcal G)
\leq
\frac{m_{\mathcal G}/2+8(n-1)/\varepsilon^2}
{\kappa_0^2\lambda_2^2(\bm L_{\mathcal G})}
\lesssim\frac1p+\frac{1}{np^2\varepsilon^2}.
\]

Projection onto the closed convex set $\Theta$ cannot increase distance to the true parameter, so the same conditional risk bound holds for $\hat\bth$. On $\mathcal H_0^c$, both $\hat\bth$ and $\bth$ belong to $\Theta$, so their squared distance is at most $4n$; this event contributes only $O(1/n)$ to the risk. Independently, the deterministic bound $R_\bth(\hat\bth)\leq4n$ holds for every $\bth\in\Theta$. Taking the smaller of this bound and the preceding display, and using $1/p\lesssim n$ in the assumed graph regime, yields
\[
R_\bth(\hat\bth)
\lesssim\frac1p+\min\left\{n,\frac{1}{np^2\varepsilon^2}\right\}.
\]
This proves the proposition.
\end{proof}
}

	\section{Omitted Proofs in Section \ref{sec: sparse GLMs}}\label{sec: sparse glm proofs}

\begingroup

Throughout this section, $s$ denotes a generic sparsity level and $P_{\bbeta}$ denotes the law of one observation from \eqref{eq: sparse glm definition}.
Boundedness of $\psi''$ gives the conditional residual bound
\begin{align}
    \E\!\left[
    \exp\{t(y-\psi'(\bm x^\top\bbeta))/c(\sigma)\}
    \,\middle|\,\bm x
    \right]
    \leq \exp\{c_2t^2/(2c(\sigma))\},                       \label{eq: sparse glm residual mgf}
\end{align}
which will also be used in the sparse upper-bound proof.

\subsection{Proof of Proposition \ref{prop: stratified sparse attack}}\label{sec: proof of prop: stratified sparse attack}

\begin{proof}[Proof of Proposition \ref{prop: stratified sparse attack}]
For a deterministic data set $\bm z$, let
\begin{align*}
    q_v(\bm z)=\Pr\{M(\bm z)\in A_v\},
\end{align*}
where the probability is only over the internal randomness of $M$.
The GLM laws have common support.  If
$L_{v,k,\ell}=dP_{v,k,\ell}/dP_{v,k,\ell-1}$, then
$L_{v,k,\ell}\geq0$ and
$\E_{P_{v,k,\ell-1}}L_{v,k,\ell}=1$.  Consequently,
\begin{align*}
    \E_{P_{v,k,\ell-1}}|T_{v,k,\ell}|
    &=\E_{P_{v,k,\ell-1}}|L_{v,k,\ell}-1|\\
    &\leq \E_{P_{v,k,\ell-1}}(L_{v,k,\ell}+1)=2.
\end{align*}
Since $0\leq q_v\leq1$, this verifies the integrability of every
attack expectation below.

We first establish the one-record identity.  Under
$R_{v,k,\ell,i-1}$, the $i$th observation has law
$P_{v,k,\ell-1}$.  Multiplying its density by $L_{v,k,\ell}$ changes
precisely that coordinate to $P_{v,k,\ell}$ and leaves the other
coordinates unchanged.  Therefore, by the definition of the
Radon--Nikodym derivative,
\begin{align}
&\E_{R_{v,k,\ell,i-1}}
\A_{v,k,\ell,i}(z_i,M(\bm Z))\notag\\
&\quad=\int q_v(\bm z)\{L_{v,k,\ell}(z_i)-1\}
dR_{v,k,\ell,i-1}(\bm z)\notag\\
&\quad=\int q_v(\bm z)\,dR_{v,k,\ell,i}(\bm z)
      -\int q_v(\bm z)\,dR_{v,k,\ell,i-1}(\bm z)\notag\\
&\quad=R_{v,k,\ell,i}\{M(\bm Z)\in A_v\}
      -R_{v,k,\ell,i-1}\{M(\bm Z)\in A_v\}.             \label{eq: one-record support swap identity}
\end{align}
For a fixed support swap $\ell$, summing
\eqref{eq: one-record support swap identity} over $i$ is an exact
telescoping sum:
\begin{align}
&\sum_{i=1}^n\E_{R_{v,k,\ell,i-1}}
\A_{v,k,\ell,i}(z_i,M(\bm Z))\notag\\
&\qquad=P_{v,k,\ell}^{\otimes n}\{M(\bm Z)\in A_v\}
       -P_{v,k,\ell-1}^{\otimes n}\{M(\bm Z)\in A_v\}.   \label{eq: one-swap support telescope}
\end{align}
Here $R_{v,k,\ell,n}=P_{v,k,\ell}^{\otimes n}$ and
$R_{v,k,\ell,0}=P_{v,k,\ell-1}^{\otimes n}$.
Summing \eqref{eq: one-swap support telescope} over $\ell$ is a
second telescoping sum.  Since $P_{v,k,0}=P_k$ and
$P_{v,k,L_{vk}}=P_v$, it gives
\begin{align*}
\sum_{\ell=1}^{L_{vk}}\sum_{i=1}^n
\E_{R_{v,k,\ell,i-1}}\A_{v,k,\ell,i}
=Q_v(A_v)-Q_k(A_v),
\end{align*}
including the case $v=k$, for which $L_{vk}=0$ and both sides are
zero.  Averaging this identity over all ordered pairs $(v,k)$ proves
\eqref{eq: support swap telescope}.

We next prove completeness.  Recall that $A_v$ is the set of outputs decoded
to the nearest packing point $\bbeta^v$, with the smallest index used to break
distance ties.  If $M(\bm Z)\notin A_v$, this definition implies that there is
some $u\neq v$ such that
$\|M(\bm Z)-\bbeta^u\|_2\leq
\|M(\bm Z)-\bbeta^v\|_2$.  The $2r$ separation and the triangle
inequality then give
\begin{align*}
    2r
    &\leq\|\bbeta^u-\bbeta^v\|_2\\
    &\leq\|M(\bm Z)-\bbeta^u\|_2
          +\|M(\bm Z)-\bbeta^v\|_2\\
    &\leq2\|M(\bm Z)-\bbeta^v\|_2.
\end{align*}
Thus $A_v^c$ is contained in
$\{\|M-\bbeta^v\|_2^2\geq r^2\}$, and Markov's inequality yields,
for every $v$,
\begin{align*}
    Q_v(A_v)
    \geq1-\frac{\E_{P_v^{\otimes n}}
    \|M(\bm Z)-\bbeta^v\|_2^2}{r^2}
    \geq1-\frac{R}{r^2}.
\end{align*}
Averaging over $v$ gives
\begin{align}
    \frac1N\sum_{v=1}^NQ_v(A_v)\geq1-\frac{R}{r^2}.       \label{eq: support decoder success lower}
\end{align}
On the other hand, the cells form a measurable partition of the output
space.  Hence $\sum_vQ_k(A_v)=1$ for every fixed $k$, and
\begin{align}
    \frac1{N^2}\sum_{v,k=1}^NQ_k(A_v)
    =\frac1{N^2}\sum_{k=1}^N1
    =\frac1N.                                               \label{eq: support decoder cross average}
\end{align}
Substituting \eqref{eq: support decoder success lower} and
\eqref{eq: support decoder cross average} into
\eqref{eq: support swap telescope} proves
\eqref{eq: support swap completeness}.

It remains to prove soundness.  Fix $v,k$, and choose a coupling
$\bm Z\sim P_v^{\otimes n}$ and $\bm Z'\sim P_k^{\otimes n}$ so that
$\E d_{\rm Ham}(\bm Z,\bm Z')\leq D$.  Put
$H_{vk}=d_{\rm Ham}(\bm Z,\bm Z')$ and
$m=\lfloor20D\rfloor$.  If $D>0$, then $m+1>20D$, so Markov's
inequality and the fact that $H_{vk}$ is integer-valued give
\begin{align*}
    \Pr(H_{vk}>m)
    =\Pr(H_{vk}\geq m+1)
    \leq\frac{\E H_{vk}}{m+1}
    <\frac1{20}.
\end{align*}
If $D=0$, then $H_{vk}=0$ almost surely and the same conclusion holds.
Let $G_{vk}=\{H_{vk}\leq m\}$.

Iterating the definition of $(\varepsilon,\delta)$-differential
privacy along a path of at most $m$ adjacent data sets gives,
pointwise on $G_{vk}$,
\begin{align*}
    q_v(\bm Z)\leq e^{m\varepsilon}q_v(\bm Z')+\delta_m,
    \qquad
    \delta_m=\delta\sum_{j=0}^{m-1}e^{j\varepsilon},
\end{align*}
with the empty sum equal to zero.  To see the displayed expression for
$\delta_m$, note that one additional neighboring-data comparison gives
$\delta_m=e^{(m-1)\varepsilon}\delta+\delta_{m-1}$, starting from
$\delta_0=0$.

We integrate the pointwise group-privacy inequality on $G_{vk}$ and
use $0\leq q_v\leq1$ on $G_{vk}^c$:
\begin{align*}
Q_v(A_v)
&=\E\{q_v(\bm Z)\1_{G_{vk}}\}
  +\E\{q_v(\bm Z)\1_{G_{vk}^c}\}\\
&\leq e^{m\varepsilon}
      \E\{q_v(\bm Z')\1_{G_{vk}}\}
      +\delta_m\Pr(G_{vk})+\Pr(G_{vk}^c)\\
&\leq e^{m\varepsilon}Q_k(A_v)+\delta_m+\frac1{20}.
\end{align*}
Thus, for every ordered pair $(v,k)$,
\begin{align}
    Q_v(A_v)\leq e^{m\varepsilon}Q_k(A_v)+\frac1{20}+\delta_m. \label{eq: support swap pair soundness}
\end{align}
Average \eqref{eq: support swap pair soundness} over $v,k$ and use
\eqref{eq: support decoder cross average}.  This gives
\begin{align*}
\frac1N\sum_{v=1}^NQ_v(A_v)
&\leq\frac{e^{m\varepsilon}}N+\frac1{20}+\delta_m.
\end{align*}
Subtracting $N^{-1}$, as in \eqref{eq: support swap telescope}, proves
\eqref{eq: support swap soundness}.

Finally, suppose, toward a contradiction, that $R<r^2/4$.  Since
$N\geq4$, completeness gives
\begin{align}
    \overline\A
    \geq1-\frac14-\frac1N\geq\frac12.                    \label{eq: support swap contradiction lower}
\end{align}
Because $m\leq20D$, choosing the universal constant in
$\varepsilon D\leq c\log N$ sufficiently small ensures
\begin{align*}
    \frac{e^{m\varepsilon}-1}{N}
    \leq\frac{e^{20\varepsilon D}}N<\frac13
    \qquad(N\geq4).
\end{align*}
When $m\geq1$, the inequality $e^\varepsilon-1\geq\varepsilon$
implies
\begin{align*}
    \delta_m
    =\delta\frac{e^{m\varepsilon}-1}{e^\varepsilon-1}
    \leq\frac{\delta}{\varepsilon}e^{m\varepsilon};
\end{align*}
when $m=0$, $\delta_m=0$.  Since $m\leq20D$, the condition
$\delta\leq c\varepsilon e^{-C\varepsilon D}$, with $C>20$ and $c$
sufficiently small, therefore ensures $\delta_m<1/10$.  Soundness now
gives
\begin{align*}
    \overline\A<\frac13+\frac1{20}+\frac1{10}<\frac12,
\end{align*}
contradicting \eqref{eq: support swap contradiction lower}.  Hence
$R\geq cr^2$ for a universal constant $c>0$.
\end{proof}

\subsection{Proof of Theorem \ref{thm: high-dim glm lb}}\label{sec: proof of thm: high-dim glm lb}

We first record the constant-weight packing used in both parts of the
lower bound.

\begin{Lemma}[Constant-weight packing]\label{lm: sparse constant weight packing}
For $2\leq s\leq d/2$ and $H_s=s\log(ed/s)$, there is a set
$\mathcal V\subseteq\{0,1\}^d$ such that every $\bm v\in\mathcal V$
has exactly $s$ nonzero coordinates and, for universal constants
$c,C>0$,
\begin{align}
    4\leq N:=|\mathcal V|,\qquad
    cH_s\leq\log N\leq CH_s,\qquad
    \|\bm v-\bm v'\|_0\geq cs
    \quad(\bm v\neq\bm v').                               \label{eq: sparse packing lemma conclusion}
\end{align}
\end{Lemma}

\begin{proof}
Let $\Omega_s$ be the collection of all $s$-subsets of $[d]$, so
$|\Omega_s|=\binom ds$.  If two supports differ by deleting $j$
coordinates and inserting $j$ coordinates, their indicator vectors
have Hamming distance $2j$.  A ball containing all supports obtainable
from a fixed support using at most $r$ such exchanges therefore has
cardinality
\begin{align*}
    V_r=\sum_{j=0}^r\binom sj\binom{d-s}{j}.
\end{align*}
Choose $r=\lfloor\rho s\rfloor$ for a sufficiently small universal
$\rho>0$.  For $j\leq r$, first choose the union of the $j$ deleted
and $j$ inserted coordinates and then distinguish the two groups.  This
gives
\begin{align*}
    \binom sj\binom{d-s}{j}
    \leq\binom d{2j}\binom{2j}j
    \leq4^j\binom d{2j}.
\end{align*}
Take $\rho\leq1/4$.  Then $2r\leq d/2$, so the binomial coefficients
$\binom d{2j}$ are increasing for $j\leq r$.  Using
$\binom ab\leq(ea/b)^b$ gives, whenever $r\geq1$,
\begin{align*}
    V_r
    \leq(r+1)4^r\binom d{2r}
    \leq(r+1)4^r\left(\frac{ed}{2r}\right)^{2r}.
\end{align*}
If $\rho s\geq2$, then
$\rho s/2\leq\lfloor\rho s\rfloor=r\leq\rho s$.  Substitution in the
preceding display gives
\begin{align*}
    \log V_r
    \leq C\rho s\log\!\left(\frac{ed}{\rho s}\right)
    \leq C'\rho s\log(ed/s).
\end{align*}
	The final inequality uses that $\rho$ is fixed and $d/s\geq2$.  If
	$r=0$, then $V_r=1$ and the required upper bound on $\log V_r$ is
	immediate. When $r=1$ and $\rho s<2$, only finitely many values of
	$s$ are possible and the preceding bound is uniform in $d$;
	adjusting the universal constants gives the same conclusion.
A greedy packing that successively selects a support and removes its
$r$-exchange ball has size at least $\binom ds/V_r$, and any two
selected supports have exchange distance greater than $r$.
Moreover,
\begin{align*}
    \log\binom ds\geq s\log(d/s)\geq c_0s\log(ed/s),
\end{align*}
where the second inequality again uses $d/s\geq2$.  Taking $\rho$
small enough gives $\log N\geq cH_s$ and Hamming separation at least
$cs$.  The upper bound follows from
$N\leq\binom ds\leq(ed/s)^s$.  Enlarging or decreasing the universal
constants handles the finitely many smallest values and ensures
$N\geq4$.
\end{proof}

\begin{Lemma}[GLM divergence and score-transport bounds]\label{lm: sparse glm divergence}
Under the assumptions of Theorem \ref{thm: high-dim glm lb}, for all $\bbeta,\bbeta'$,
\begin{align}
    D_{\rm KL}(P_{\bbeta}\|P_{\bbeta'})
    &\leq \frac{C}{c(\sigma)}
    \|\bbeta-\bbeta'\|_2^2,                                 \label{eq: sparse glm one sample kl}\\
    d_{\rm TV}(P_{\bbeta},P_{\bbeta'})
    &\leq \frac{C}{\sqrt{c(\sigma)}}
    \|\bbeta-\bbeta'\|_2.                                   \label{eq: sparse glm one sample tv}
\end{align}
Moreover, $P_{\bbeta}^n$ and $P_{\bbeta'}^n$ have a coupling whose expected Hamming distance is at most
\begin{align}
    \frac{Cn}{\sqrt{c(\sigma)}}\|\bbeta-\bbeta'\|_2.         \label{eq: sparse glm product coupling}
\end{align}
\end{Lemma}

\begin{proof}
Put $\bm\Delta=\bbeta'-\bbeta$.  Since the marginal distribution of $\bm x$ is the same under both parameters and
$\E_\bbeta(y\mid\bm x)=\psi'(\bm x^\top\bbeta)$,
\begin{align*}
    D_{\rm KL}(P_{\bbeta}\|P_{\bbeta'})
    &=\frac1{c(\sigma)}
    \E_{\bm x}\left[
    \psi(\bm x^\top\bbeta')-\psi(\bm x^\top\bbeta)
    -\psi'(\bm x^\top\bbeta)\bm x^\top\bm\Delta
    \right]\\
    &\leq\frac{c_2}{2c(\sigma)}
    \E(\bm x^\top\bm\Delta)^2
    \leq\frac{C}{c(\sigma)}\|\bm\Delta\|_2^2.
\end{align*}
This proves \eqref{eq: sparse glm one sample kl}.  To obtain the total-variation bound directly from the score, let
$\bbeta_t=\bbeta+t\bm\Delta$ and let $p_t$ be the density of
$P_{\bbeta_t}$.  Along this path,
$\partial_t p_t(z)=p_t(z)\langle\bm\Delta,S_{\bbeta_t}(z)\rangle$.
The Fisher-information bound below makes this derivative integrable in
$L^1$, so the fundamental theorem of calculus gives
\begin{align*}
    p_1(z)-p_0(z)
    =\int_0^1p_t(z)\langle\bm\Delta,S_{\bbeta_t}(z)\rangle\,dt.
\end{align*}
For the GLM,
\begin{align*}
    S_{\bbeta}(y,\bm x)
    &=\frac{y-\psi'(\bm x^\top\bbeta)}{c(\sigma)}\bm x,\\
    \mathcal I(\bbeta)
    &=\frac1{c(\sigma)}\E\{\psi''(\bm x^\top\bbeta)
    \bm x\bm x^\top\}\preceq\frac{C}{c(\sigma)}\bm I_d.
\end{align*}
Consequently, Tonelli's theorem and Cauchy--Schwarz yield
\begin{align*}
    d_{\rm TV}(P_{\bbeta},P_{\bbeta'})
    &\leq\frac12\int_0^1
    \E_{\bbeta_t}|\langle\bm\Delta,S_{\bbeta_t}\rangle|\,dt\\
    &\leq\frac12\int_0^1
    \{\bm\Delta^\top\mathcal I(\bbeta_t)\bm\Delta\}^{1/2}\,dt
    \leq\frac{C}{\sqrt{c(\sigma)}}\|\bm\Delta\|_2,
\end{align*}
which proves \eqref{eq: sparse glm one sample tv}.  Couple the $n$ coordinate pairs independently by maximal couplings.  The expected number of unequal coordinates is $n\,d_{\rm TV}(P_{\bbeta},P_{\bbeta'})$, which proves
\eqref{eq: sparse glm product coupling}.
\end{proof}

\begin{proof}[Proof of Theorem \ref{thm: high-dim glm lb}]
Apply Lemma \ref{lm: sparse constant weight packing} with $s=s^*$.
It gives a set $\mathcal V\subseteq\{0,1\}^d$ such that every
$\bm v\in\mathcal V$ has exactly $s^*$ ones and
\begin{align}
    cH\leq\log N\leq CH,\qquad N=|\mathcal V|\geq4,\qquad
    \|\bm v-\bm v'\|_0\geq cs^*
    \quad(\bm v\neq\bm v').                                 \label{eq: sparse constant weight packing}
\end{align}
For an amplitude $a>0$, set
\begin{align*}
    \bbeta^v=a\bm v,\qquad v\in\mathcal V.
\end{align*}
Then
\begin{align}
    ca^2s^*\leq\|\bbeta^v-\bbeta^{v'}\|_2^2
    \leq 2a^2s^*,                                          \label{eq: sparse packing separation}
\end{align}
and $\|\bbeta^v\|_2^2=a^2s^*$.  The upper bound in
\eqref{eq: sparse packing separation} uses that two binary vectors of
weight $s^*$ differ in at most $2s^*$ coordinates.

\textit{The non-private term.}
Choose
\begin{align*}
    a_{\rm np}^2=c_4\frac{c(\sigma)H}{ns^*}
\end{align*}
with $c_4$ sufficiently small.  Lemma \ref{lm: sparse glm divergence},
\eqref{eq: sparse packing separation}, and tensorization of relative
entropy give
\begin{align*}
    \max_{v\neq v'}
    D_{\rm KL}(P_{\bbeta^v}^n\|P_{\bbeta^{v'}}^n)
    &\leq \frac{Cn}{c(\sigma)}
       \|\bbeta^v-\bbeta^{v'}\|_2^2\\
    &\leq Cn a_{\rm np}^2s^*/c(\sigma)
     \leq Cc_4H
     \leq c\log N.
\end{align*}
The theorem's unsaturated condition contains, in particular,
$c(\sigma)H/n\leq c_2'B^2$.  Hence, after decreasing $c_4$ if
necessary,
\begin{align*}
	    \|\bbeta^v\|_2^2=a_{\rm np}^2s^*\leq B^2,
\end{align*}
so every packing point belongs to $\Theta_{s^*}(B)$.

We spell out the reduction to testing.  Let $V$ be uniform on
$\mathcal V$, generate $\bm Z\sim P_{\bbeta^V}^{\otimes n}$, and,
for any estimator $M$, let $\widehat V$ be a nearest packing point to
$M(\bm Z)$.  If $\widehat V\neq V$, then the triangle inequality and
the nearest-point property imply
\begin{align*}
    \min_{u\neq v}\|\bbeta^u-\bbeta^v\|_2
    \leq2\|M(\bm Z)-\bbeta^V\|_2.
\end{align*}
Thus, with
$r_{\rm np}^2=\frac14\min_{u\neq v}
\|\bbeta^u-\bbeta^v\|_2^2\gtrsim a_{\rm np}^2s^*$,
\begin{align*}
    \sup_v\E_v\|M-\bbeta^v\|_2^2
    \geq r_{\rm np}^2\Pr(\widehat V\neq V).
\end{align*}
Fano's inequality and the preceding pairwise KL bound give
\begin{align*}
    \Pr(\widehat V\neq V)
    \geq1-\frac{I(V;\bm Z)+\log2}{\log N}
    \geq c,
\end{align*}
where the final constant is positive after choosing $c_4$ small; here
$N\geq4$ absorbs the $\log2$ term.  It follows that
\begin{align}
    \inf_M\sup_{\bbeta\in\Theta_{s^*}(B)}
    \E_\bbeta\|M-\bbeta\|_2^2
    \gtrsim \frac{c(\sigma)H}{n}.                            \label{eq: sparse fano nonprivate term}
\end{align}

\textit{The privacy term.}
Write $\ell_n=\log(en)$ and choose
\begin{align*}
    a_{\rm priv}^2
    =c_5\frac{c(\sigma)H^2}{\ell_n n^2\varepsilon^2s^*}
\end{align*}
with $c_5$ sufficiently small.  For every $v\neq v'$, Lemma
\ref{lm: sparse glm divergence} and
\eqref{eq: sparse packing separation} give a coupling of
$P_{\bbeta^v}^{\otimes n}$ and
$P_{\bbeta^{v'}}^{\otimes n}$ whose expected Hamming distance is at
most
\begin{align*}
    D
    &\leq\frac{Cn}{\sqrt{c(\sigma)}}
      \|\bbeta^v-\bbeta^{v'}\|_2
      \leq Cn a_{\rm priv}\sqrt{s^*/c(\sigma)}.
\end{align*}
Substitution of the definition of $a_{\rm priv}$ shows explicitly
that
\begin{align}
    \varepsilon D
    \leq C\sqrt{c_5}\frac{H}{\sqrt{\ell_n}}
    \leq c\log N,                                         \label{eq: sparse privacy coupling condition}
\end{align}
where the second inequality follows from $\log N\geq cH$ and a
sufficiently small $c_5$.  The theorem's condition
\begin{align*}
    \delta\leq c_0\varepsilon
    \exp\{-c_1H/\sqrt{\ell_n}\}
\end{align*}
implies
$\delta\leq c\varepsilon\exp\{-C\varepsilon D\}$ after choosing
$c_1$ larger than the universal constant in
\eqref{eq: sparse privacy coupling condition} and then choosing $c_0$
small enough.

As in the non-private part, the theorem's unsaturated condition implies
\begin{align*}
    a_{\rm priv}^2s^*
    =c_5\frac{c(\sigma)H^2}{\ell_n n^2\varepsilon^2}
    \leq B^2
\end{align*}
after decreasing $c_5$.  Therefore the packing points belong to
$\Theta_{s^*}(B)$.  Every intermediate point on a shortest support-swap
path also has exactly $s^*$ coordinates, each equal to
$a_{\rm priv}$, and hence has the same Euclidean norm.  Thus all
intermediate GLM laws used by the attack remain in the stated parameter
space.

By \eqref{eq: sparse packing separation}, the packing is $2r_{\rm
priv}$-separated for a choice satisfying
\begin{align*}
    r_{\rm priv}^2\asymp a_{\rm priv}^2s^*.
\end{align*}
All hypotheses of Proposition \ref{prop: stratified sparse attack} have
now been verified: $N\geq4$, the coupling condition is
\eqref{eq: sparse privacy coupling condition}, and the required
approximate-privacy condition follows from the theorem's assumption.
The proposition therefore gives
\begin{align}
    \inf_{M\in\mathcal M_{\varepsilon,\delta}}
    \sup_{\bbeta\in\Theta_{s^*}(B)}
    \E_\bbeta\|M-\bbeta\|_2^2
    &\gtrsim r_{\rm priv}^2\notag\\
    &\asymp a_{\rm priv}^2s^*
     =\frac{c_5c(\sigma)H^2}
     {\ell_n n^2\varepsilon^2}.                            \label{eq: sparse support swap privacy term}
\end{align}
	This completes the soundness--completeness comparison for the
	support-swap score attack.

The non-private lower bound applies to all estimators and hence also to
the private subclass.  Taking the maximum of
\eqref{eq: sparse fano nonprivate term} and
\eqref{eq: sparse support swap privacy term}, and using
$\max\{u,v\}\geq(u+v)/2$, proves \eqref{eq: high-dim glm lb}.
\end{proof}

\subsection{Proof of Proposition \ref{prop: sparse glm empirical curvature}}\label{sec: proof of prop: sparse glm empirical curvature}

\begin{proof}[Proof of Proposition \ref{prop: sparse glm empirical curvature}]
Fix a deterministic $L_\psi\geq L_{\psi,0}^{\rm sp}$, where the lower
threshold is chosen in Step 1 below.
Write
\begin{align*}
q=s+s^*,\qquad k=2s+s^*,\qquad m=q+k=3s+2s^*,
\end{align*}
and introduce the larger parameter class
\begin{align*}
\mathcal B_q(3)=\{\bbeta\in\mathbb R^d:
\|\bbeta\|_0\leq q,\ \|\bbeta\|_2\leq3\}.
\end{align*}
Every parameter appearing in the definition of
$\mathcal E_{\rm sc}(\alpha,\gamma;L_\psi)$ belongs to this class.
Indeed, if
$\bm u$ is $s$-sparse and
$\|\bm u-\bbeta^*\|_2\leq2$, then
$\|\bm u\|_2\leq3$ and
$|\supp(\bm u)\cup\supp(\bbeta^*)|\leq q$; the same statements hold
along the line segment from $\bbeta^*$ to $\bm u$.  It is therefore
enough to prove uniform Hessian bounds over $\mathcal B_q(3)$ and over
all coordinate sets of size at most $k$.

For $\bbeta\in\mathcal B_q(3)$, define the population and empirical
Hessians
\begin{align*}
H_{L_\psi}(\bbeta)&=\E\{\psi_{L_\psi}''(\bm x^\top\bbeta)
\bm x\bm x^\top\},\\
H_{n,L_\psi}(\bbeta)&=\frac1n\sum_{i=1}^n
\psi_{L_\psi}''(\bm x_i^\top\bbeta)\bm x_i\bm x_i^\top.
\end{align*}
Notice that
$H_{n,L_\psi}(\bbeta)=\nabla^2\L_{n,L_\psi}^R(\bbeta)$:
truncating the responses changes the gradient but not the Hessian.

\textit{Step 1: uniform population curvature.}
Assumption (D) and the standard moment bound for a sub-Gaussian random
variable give, for every unit vector $\bm v$,
\begin{align}
\E(\bm v^\top\bm x)^2&\geq C^{-1},&
\E|\bm v^\top\bm x|^4&\leq C_1K_{\bm x}^4.                 \label{eq: sparse curvature moments}
\end{align}
Also, uniformly over $\bbeta\in\mathcal B_q(3)$,
\begin{align}
\Pr(|\bm x^\top\bbeta|>L)
\leq2\exp\left\{-\frac{cL^2}{K_{\bm x}^2\|\bbeta\|_2^2}\right\}
\leq2\exp\left\{-\frac{cL^2}{9K_{\bm x}^2}\right\},      \label{eq: sparse predictor tail}
\end{align}
where the assertion is immediate when $\bbeta=\bm0$.
Cauchy--Schwarz, \eqref{eq: sparse curvature moments}, and
\eqref{eq: sparse predictor tail} imply
\begin{align*}
&\E\{(\bm v^\top\bm x)^2
\1(|\bm x^\top\bbeta|>L)\}\\
&\qquad\leq
\{\E|\bm v^\top\bm x|^4\}^{1/2}
\Pr(|\bm x^\top\bbeta|>L)^{1/2}
\leq C_2K_{\bm x}^2e^{-c_2L^2/K_{\bm x}^2}.
\end{align*}
Choose a fixed $L_0$ so that the last expression is at most
$(2C)^{-1}$ and set $L_{\psi,0}^{\rm sp}=\max\{1,L_0\}$.
Assumption (G) and continuity give
\begin{align*}
a_0:=\inf_{|t|\leq L_0}\psi''(t)>0.
\end{align*}
For every $L_\psi\geq L_{\psi,0}^{\rm sp}$,
$\psi_{L_\psi}''(t)=\psi''(t)$ on $[-L_0,L_0]$ and
$0\leq\psi_{L_\psi}''(t)\leq c_2$ everywhere. It follows that,
uniformly over unit $\bm v$ and $\bbeta\in\mathcal B_q(3)$,
\begin{align*}
\bm v^\top H_{L_\psi}(\bbeta)\bm v
&\geq a_0\E\{(\bm v^\top\bm x)^2
\1(|\bm x^\top\bbeta|\leq L_0)\}
\geq\frac{a_0}{2C}=:\alpha_0,\\
\bm v^\top H_{L_\psi}(\bbeta)\bm v
&\leq c_2\E(\bm v^\top\bm x)^2
\leq c_2C=:\gamma_0.
\end{align*}
Thus
\begin{align}
\alpha_0I_d\preceq H_{L_\psi}(\bbeta)\preceq\gamma_0I_d
\quad(\bbeta\in\mathcal B_q(3),\
L_\psi\geq L_{\psi,0}^{\rm sp}).                           \label{eq: sparse population curvature}
\end{align}

\textit{Step 2: a uniform sparse-design norm bound.}
After increasing the constant in the sample-size assumption, we may
suppose $m\leq d$.  Define
\begin{align*}
R_m=C_3K_{\bm x}
\{m\log(ed/m)+\log n\}^{1/2}.
\end{align*}
We claim that, for a sufficiently large fixed $C_3$,
\begin{align}
\mathcal E_m=
\left\{\max_{1\leq i\leq n}\max_{K\subset[d],\ |K|\leq m}
\|\bm x_{i,K}\|_2\leq R_m\right\}                         \label{eq: sparse design norm event}
\end{align}
has probability at least $1-n^{-4}$.  To verify the claim, pad any set
with fewer than $m$ coordinates to one of size $m$.  For a fixed
$m$-element set $K$, let $\mathcal N_K$ be a $1/2$-net of its unit
sphere with $|\mathcal N_K|\leq5^m$.  The net property and (D) give
\begin{align*}
\Pr(\|\bm x_K\|_2>2t)
&\leq\sum_{\bm v\in\mathcal N_K}
\Pr(|\bm v^\top\bm x|>t)\\
&\leq2\cdot5^m\exp(-ct^2/K_{\bm x}^2).
\end{align*}
A union bound over $i\in[n]$ and at most
$\binom dm\leq(ed/m)^m$ supports proves
\eqref{eq: sparse design norm event} and the asserted probability.
We will also use the deterministic moment consequence of (D)
\begin{align}
\sup_{|K|\leq m}\E\|\bm x_K\|_2^3
\leq C_4K_{\bm x}^3m^{3/2}.                                 \label{eq: sparse design third moment}
\end{align}
Indeed, the coordinatewise sub-Gaussian fourth-moment bound gives
$\E\|\bm x_K\|_2^4\leq C K_{\bm x}^4m^2$, and Lyapunov's
inequality yields \eqref{eq: sparse design third moment}.

\textit{Step 3: concentration on sparse nets.}
Fix $\bbeta$, a coordinate set $I$ with $|I|\leq k$, and a unit vector
$\bm v$ supported on $I$.  The variables
\begin{align*}
Z_i=\psi_{L_\psi}''(\bm x_i^\top\bbeta)(\bm v^\top\bm x_i)^2
\end{align*}
are sub-exponential with norm at most $C c_2K_{\bm x}^2$, uniformly in
$\bbeta,I,$ and $\bm v$.  Bernstein's inequality therefore gives
\begin{align}
\Pr\left(\left|\bm v^\top
\{H_{n,L_\psi}(\bbeta)-H_{L_\psi}(\bbeta)\}\bm v\right|>\frac{\alpha_0}{64}\right)
\leq2e^{-c_{\rm B}n}.                                        \label{eq: sparse fixed hessian concentration}
\end{align}

Let
\begin{align*}
\widetilde c_3=c_3+c_2\|\chi'\|_\infty,\qquad
L_m=\widetilde c_3\{R_m^3+C_4K_{\bm x}^3m^{3/2}\},\qquad
r_m=\min\left\{1,\frac{\alpha_0}{64L_m}\right\},
\end{align*}
with $r_m=1$ if $L_m=0$.  For every coordinate set $J$ of size at most
$q$, take an $r_m$-net of the radius-three Euclidean ball on $J$; it
may be chosen to have cardinality at most $(1+6/r_m)^q$.  For every
coordinate set $I$ of size at most $k$, take a $1/4$-net of its unit
sphere with cardinality at most $9^k$.  Padding smaller supports if
necessary, the total logarithm of the number of triples consisting of
a parameter-net point, a direction-net point, and their two supports
is at most
\begin{align}
k\log(ed/k)+k\log9
+q\log(ed/q)+q\log(1+6/r_m)
\leq C_5(s+s^*)\log(ed).                                    \label{eq: sparse hessian net entropy}
\end{align}
For the last inequality, note that $m\lesssim s+s^*$,
$n\leq d$, and the definition of $R_m$ gives
$\log(1+6/r_m)\lesssim\log(ed)$.

Apply \eqref{eq: sparse fixed hessian concentration} to all these
triples.  The sample-size condition and a union bound show that, except
on an event of probability at most $n^{-4}$,
\begin{align*}
\left|\bm v^\top\{H_{n,L_\psi}(\widetilde\bbeta)
-H_{L_\psi}(\widetilde\bbeta)\}\bm v\right|
\leq\frac{\alpha_0}{64}
\end{align*}
simultaneously over the nets.  The standard quadratic-form net
inequality
\begin{align*}
\|A\|_{\rm op}\leq2\max_{\bm v\in\mathcal N_I}
|\bm v^\top A\bm v|
\end{align*}
for a symmetric matrix on $I$ and a $1/4$-net $\mathcal N_I$ therefore
gives
\begin{align}
\max_{|I|\leq k}\max_{\widetilde\bbeta\ \text{ in the parameter nets}}
\|[H_{n,L_\psi}(\widetilde\bbeta)
-H_{L_\psi}(\widetilde\bbeta)]_{I,I}\|_{\rm op}
\leq\frac{\alpha_0}{32}.                                    \label{eq: sparse hessian bound on nets}
\end{align}

\textit{Step 4: passage from the parameter nets to the full class.}
Fix $\bbeta\in\mathcal B_q(3)$, choose a support $J$ of size at most
$q$ containing its support, and select a net point
$\widetilde\bbeta$ on the same support with
$\|\bbeta-\widetilde\bbeta\|_2\leq r_m$. Equation
\eqref{eq: saturated link score} and $L_\psi\geq1$ give
\begin{align*}
\left|\frac{\d}{\d t}\psi_{L_\psi}''(t)\right|
\leq c_3+c_2\|\chi'\|_\infty=\widetilde c_3.
\end{align*}
For a unit vector $\bm v$ supported on $I$, $|I|\leq k$, the
mean-value theorem therefore yields
\begin{align*}
&\left|\bm v^\top\{H_{n,L_\psi}(\bbeta)
-H_{n,L_\psi}(\widetilde\bbeta)\}\bm v\right|\\
&\quad\leq
\frac{\widetilde c_3}{n}\sum_{i=1}^n
|\bm x_i^\top(\bbeta-\widetilde\bbeta)|
(\bm v^\top\bm x_i)^2\\
&\quad\leq \widetilde c_3r_m\frac1n\sum_{i=1}^n
\|\bm x_{i,I\cup J}\|_2^3
\leq \widetilde c_3r_mR_m^3
\end{align*}
on $\mathcal E_m$.  Similarly, \eqref{eq: sparse design third moment}
implies
\begin{align*}
\left|\bm v^\top\{H_{L_\psi}(\bbeta)
-H_{L_\psi}(\widetilde\bbeta)\}\bm v\right|
\leq C_4\widetilde c_3r_mK_{\bm x}^3m^{3/2}.
\end{align*}
The sum of the last two bounds is at most
$r_mL_m\leq\alpha_0/64$.  Combining this fact with
\eqref{eq: sparse hessian bound on nets} shows that
\begin{align}
\sup_{\bbeta\in\mathcal B_q(3)}\max_{|I|\leq k}
\|[H_{n,L_\psi}(\bbeta)-H_{L_\psi}(\bbeta)]_{I,I}\|_{\rm op}
\leq\frac{3\alpha_0}{64}                                    \label{eq: sparse uniform empirical hessian}
\end{align}
outside an event of probability at most $2n^{-4}\leq n^{-3}$.
Together with \eqref{eq: sparse population curvature}, this proves the
claimed event, for example with
\begin{align*}
\alpha=\frac{\alpha_0}{2},\qquad
\gamma=\gamma_0+\frac{\alpha_0}{2}.
\end{align*}
The constants do not depend on $\bbeta^*$, so the probability bound is
uniform over the stated parameter space.
\end{proof}

\subsection{Proof of Lemmas \ref{lm: repaired noisy projection} and \ref{lm: repaired noisy iht contraction}}\label{sec: proofs of sparse glm auxiliary lemmas}

\begin{Lemma}[Noisy projection inequality]\label{lm: repaired noisy projection}
Let $\widehat{\bm v}$ be the output of Algorithm \ref{algo: noisy hard thresholding} on $\bm v$, let $S$ be its selected support, and let $\bm\theta$ be $s_0$-sparse with $s>s_0$.  Put
\begin{align*}
    I=S\cup\supp(\bm\theta),\qquad
    \mathcal W^2=\sum_{r=1}^s\|\bm w_r\|_\infty^2+
    \|\widetilde{\bm w}_S\|_2^2.
\end{align*}
Then
\begin{align}
    \|\widehat{\bm v}-\bm\theta\|_2
    \leq
    \left(1+C\sqrt{\frac{s_0}{s-s_0}}\right)
    \|(\bm v-\bm\theta)_I\|_2+C\mathcal W.                  \label{eq: repaired noisy projection}
\end{align}
\end{Lemma}

\subsubsection{Proof of Lemma \ref{lm: repaired noisy projection}}

\begin{proof}
Write $\overline{\bm v}=\bm v_S$ for the selected vector before the
final perturbation, with zeros outside $S$.  We may restrict the
selection procedure to the coordinates in $I$.  Indeed, at round $r$
the coordinate selected by the original procedure belongs to
$S\subseteq I$ and maximizes the noisy criterion over all coordinates
not selected previously.  It therefore also maximizes that criterion
over the smaller set of remaining coordinates in $I$.  Induction over
$r=1,\ldots,s$ shows that the restricted procedure selects the same
ordered sequence and hence the same set $S$.

Let $T$ contain the $s$ largest coordinates of $\bm v_I$ in absolute
value, with deterministic tie breaking.  Since $|S|=|T|=s$,
$|S\setminus T|=|T\setminus S|$.  Match every
$j\in T\setminus S$ with a distinct coordinate
$k(j)\in S\setminus T$, and let $r(j)$ be the round at which $k(j)$
was selected.  The coordinate $j$ is never selected, so it is still
available at round $r(j)$.  The selection rule therefore gives
\begin{align*}
    |v_j|+w_{r(j),j}
    \leq |v_{k(j)}|+w_{r(j),k(j)}.
\end{align*}
It follows that
\begin{align*}
    |v_j|
    \leq |v_{k(j)}|
       +|w_{r(j),k(j)}|+|w_{r(j),j}|
    \leq |v_{k(j)}|+2\|\bm w_{r(j)}\|_\infty.
\end{align*}
For every $c>0$, the inequality
$(x+y)^2\leq(1+1/c)x^2+(1+c)y^2$ now yields
\begin{align*}
    v_j^2
    \leq (1+1/c)v_{k(j)}^2
    +4(1+c)\|\bm w_{r(j)}\|_\infty^2.
\end{align*}
The sets $I\setminus S$ and $T^c=I\setminus T$ decompose as
\begin{align*}
I\setminus S
&=(T\setminus S)\mathbin{\dot\cup}
  \bigl(I\setminus(S\cup T)\bigr),\\
T^c
&=(S\setminus T)\mathbin{\dot\cup}
  \bigl(I\setminus(S\cup T)\bigr).
\end{align*}
Summing the matched inequalities, adding the common coordinates in
$I\setminus(S\cup T)$, and using that each selection round is used at
most once give
\begin{align}
\|\overline{\bm v}-\bm v_I\|_2^2
&=\sum_{j\in I\setminus S}v_j^2\notag\\
&\leq(1+1/c)\sum_{j\in T^c}v_j^2
 +4(1+c)\sum_{r=1}^s\|\bm w_r\|_\infty^2.                \label{eq: noisy projection intermediate}
\end{align}

We next compare the deterministic tail in
\eqref{eq: noisy projection intermediate} with the approximation error
to $\bm\theta$.  Put $m=|I|$ and
$S_0=\supp(\bm\theta)$.  The set $T^c$ contains the $m-s$ smallest
squared coordinates of $\bm v_I$, whereas $I\setminus S_0$ contains
at least $m-s_0$ coordinates.  If $m=s$, then $T^c$ is empty and the
desired tail bound is immediate.  If $m>s$, choose any
$A\subseteq I\setminus S_0$ with $|A|=m-s_0$.  The average of the
$m-s$ smallest values cannot exceed the average over $A$, so
\begin{align*}
    \frac{1}{m-s}\sum_{j\in T^c}v_j^2
    \leq
    \frac{1}{m-s_0}\sum_{j\in A}v_j^2
    \leq
    \frac{\|\bm v_I-\bm\theta\|_2^2}{m-s_0}.
\end{align*}
Because $S\subseteq I=S\cup S_0$, we have $s\leq m\leq s+s_0$.
The function $x\mapsto(x-s)/(x-s_0)$ is increasing for $x\geq s$,
and hence
\begin{align*}
    \frac{m-s}{m-s_0}
    \leq\frac{(s+s_0)-s}{(s+s_0)-s_0}
    =\frac{s_0}{s}.
\end{align*}
We have therefore proved, including the case $m=s$, that
\begin{align}
    \sum_{j\in T^c}v_j^2
    \leq\frac{s_0}{s}\|\bm v_I-\bm\theta\|_2^2.           \label{eq: noisy projection tail comparison}
\end{align}

Finally,
$\widehat{\bm v}=\overline{\bm v}+\widetilde{\bm w}_S$ and
$\bm\theta$ is supported on $I$.  The triangle inequality,
\eqref{eq: noisy projection intermediate}, and
\eqref{eq: noisy projection tail comparison} give
\begin{align*}
    \|\widehat{\bm v}-\bm\theta\|_2
    &\leq\|\overline{\bm v}-\bm v_I\|_2
    +\|\bm v_I-\bm\theta\|_2
    +\|\widetilde{\bm w}_S\|_2\\
    &\leq
    \left\{1+\sqrt{(1+1/c)s_0/s}\right\}
    \|(\bm v-\bm\theta)_I\|_2
    +2\sqrt{1+c}\left\{\sum_{r=1}^s
    \|\bm w_r\|_\infty^2\right\}^{1/2}
    +\|\widetilde{\bm w}_S\|_2.
\end{align*}
Taking $c=1$, using
$\sqrt{s_0/s}\leq\sqrt{s_0/(s-s_0)}$, and recalling the definition of
$\mathcal W$ prove \eqref{eq: repaired noisy projection}.
\end{proof}

\begin{Lemma}[Deterministic contraction about the truth]\label{lm: repaired noisy iht contraction}
Suppose \eqref{eq: sparse glm curvature event} holds and no design
coordinate is clipped, so that
\(\overline\L_{n,L_\psi}^R(\bbeta)=\L_{n,L_\psi}^R(\bbeta)\)
for every \(\bbeta\).
Let
$\bm g^*=\nabla\L_{n,L_\psi}^R(\bbeta^*)$, choose
$\eta^0=2/(\alpha+\gamma)$, and let $\bbeta^{t+1}$ be one step of Algorithm \ref{algo: private sparse glm}.  There is a sufficiently large universal $C_0$ such that, when
$s\geq C_0(\gamma/\alpha)^2s^*$ and
$\|\bbeta^t-\bbeta^*\|_2\leq2$,
\begin{align}
    \|\bbeta^{t+1}-\bbeta^*\|_2^2
    \leq q\|\bbeta^t-\bbeta^*\|_2^2
    +Cs\|\bm g^*\|_\infty^2+C\mathcal W_t^2,                \label{eq: repaired noisy iht recursion}
\end{align}
where $q\in(0,1)$ depends only on $\alpha/\gamma$ and
$\mathcal W_t$ is defined as in Lemma \ref{lm: repaired noisy projection}.
\end{Lemma}

\subsubsection{Proof of Lemma \ref{lm: repaired noisy iht contraction}}

\begin{proof}
Let $S^t$ and $S^{t+1}$ denote the coordinate sets selected by NoisyHT
when producing $\bbeta^t$ and $\bbeta^{t+1}$, respectively; take
$S^0=\varnothing$.  Thus $\bbeta^t$ is supported on $S^t$, even if a
selected coordinate happens to equal zero after the final perturbation.
Set $\bm e^t=\bbeta^t-\bbeta^*$ and
\begin{align*}
    I^t=S^t\cup S^{t+1}\cup\supp(\bbeta^*).
\end{align*}
Since $|S^t|,|S^{t+1}|\leq s$,
$|I^t|\leq2s+s^*$.  Moreover, $\bm e^t$ is supported on $I^t$.
The integral mean-value formula therefore gives
\begin{align*}
    [\nabla\L_{n,L_\psi}^R(\bbeta^t)
    -\nabla\L_{n,L_\psi}^R(\bbeta^*)]_{I^t}
    =\bm H_{I^t,I^t}^t\bm e_{I^t}^t,
\end{align*}
where
\begin{align*}
    \bm H^t=\int_0^1\nabla^2\L_{n,L_\psi}^R
    (\bbeta^*+u\bm e^t)\,\d u.
\end{align*}
The assumed bound $\|\bm e^t\|_2\leq2$ ensures that every point on the
line segment in this integral is covered by
\eqref{eq: sparse glm curvature event}.  Integration preserves the
Loewner inequalities, so every eigenvalue of
$\bm H_{I^t,I^t}^t$ lies in $[\alpha,\gamma]$.

Put
\begin{align*}
    \rho=\frac{\gamma-\alpha}{\gamma+\alpha}.
\end{align*}
Because $\eta^0=2/(\alpha+\gamma)$, spectral calculus gives
\begin{align*}
\|\bm I-\eta^0\bm H_{I^t,I^t}^t\|_{\rm op}
&=\max_{\lambda\in[\alpha,\gamma]}
  \left|1-\frac{2\lambda}{\alpha+\gamma}\right|\\
&=\rho.
\end{align*}
Consequently,
\begin{align}
    \left\|
    \bm e_{I^t}^t-\eta^0
    [\nabla\L_{n,L_\psi}^R(\bbeta^t)
    -\nabla\L_{n,L_\psi}^R(\bbeta^*)]_{I^t}
    \right\|_2
    \leq \frac{\gamma-\alpha}{\gamma+\alpha}
    \|\bm e^t\|_2.                                         \label{eq: restricted gradient contraction}
\end{align}

Apply Lemma \ref{lm: repaired noisy projection} to
\begin{align*}
    \bm v^t=\bbeta^t-\eta^0\nabla\L_{n,L_\psi}^R(\bbeta^t),
    \qquad \bm\theta=\bbeta^*.
\end{align*}
The set denoted by $I$ in that lemma is now
$J^t=S^{t+1}\cup\supp(\bbeta^*)$, which is contained in $I^t$.
Decompose the restricted pre-thresholding error as
\begin{align*}
(\bm v^t-\bbeta^*)_{J^t}
&=\bigl[\bm e^t-\eta^0
  \{\nabla\L_{n,L_\psi}^R(\bbeta^t)
  -\nabla\L_{n,L_\psi}^R(\bbeta^*)\}
  \bigr]_{J^t}
  -\eta^0\bm g^*_{J^t}.
\end{align*}
The first term is bounded by the same norm over $I^t$, and the second
term is bounded by
$\eta^0\sqrt{|J^t|}\|\bm g^*\|_\infty$.  Since
$|J^t|\leq s+s^*$, \eqref{eq: restricted gradient contraction} gives
\begin{align}
\|(\bm v^t-\bbeta^*)_{J^t}\|_2
\leq \rho\|\bm e^t\|_2
     +\eta^0\sqrt{s+s^*}\|\bm g^*\|_\infty.               \label{eq: restricted prethreshold error}
\end{align}

Let
\begin{align*}
    \kappa_s=1+C_p\sqrt{\frac{s^*}{s-s^*}},
\end{align*}
where $C_p$ is the universal constant in Lemma
\ref{lm: repaired noisy projection}.  Combining that lemma with
\eqref{eq: restricted prethreshold error} yields
\begin{align}
    \|\bbeta^{t+1}-\bbeta^*\|_2
    \leq \kappa_s\rho\|\bm e^t\|_2
    +\kappa_s\eta^0\sqrt{s+s^*}\|\bm g^*\|_\infty
    +C_p\mathcal W_t.                                      \label{eq: noisy contraction before squaring}
\end{align}

We now verify quantitatively that the first coefficient is below one.
Put $\tau=\alpha/\gamma\in(0,1]$.  Then
\begin{align*}
    \rho=\frac{1-\tau}{1+\tau}\leq1-\tau.
\end{align*}
Choosing $C_0$ sufficiently large in
$s\geq C_0\tau^{-2}s^*$ ensures
\begin{align*}
    C_p\sqrt{\frac{s^*}{s-s^*}}\leq\frac\tau4.
\end{align*}
Since $\rho\leq1$, this gives
\begin{align}
    a:=\kappa_s\rho
    \leq\rho+\frac\tau4
    \leq1-\frac{3\tau}{4}<1.                              \label{eq: contraction coefficient below one}
\end{align}
Write the last two terms on the right of
\eqref{eq: noisy contraction before squaring} as $b_t$.  Because
$s\geq s^*$, $\kappa_s$ and $\eta^0$ are bounded by constants depending
only on $\alpha/\gamma$ and the fixed curvature scale, and hence
\begin{align*}
    b_t^2\leq C\{s\|\bm g^*\|_\infty^2+\mathcal W_t^2\}.
\end{align*}
For $0\leq a<1$ and $x,b\geq0$, completing the square gives the useful
form of Young's inequality
\begin{align*}
    (ax+b)^2\leq ax^2+\frac{b^2}{1-a}.
\end{align*}
Apply this inequality to \eqref{eq: noisy contraction before
squaring}.  \furtherrev{Set $q_0=1-3\alpha/(4\gamma)$. By
\eqref{eq: contraction coefficient below one}, $a\leq q_0<1$; hence
the resulting recursion remains valid with $q=q_0$ (also when $a=0$).}
Thus \eqref{eq: contraction coefficient below one} proves
\eqref{eq: repaired noisy iht recursion}.  All constants depend only
on the fixed curvature constants.
\end{proof}

\subsection{Proof of Theorem \ref{thm: glm upper bound}}\label{sec: proof of thm: glm upper bound}

\begin{proof}[Proof of Theorem \ref{thm: glm upper bound}]
\textit{Privacy.}
Condition on all releases before iteration $t$.  The iterate
$\bbeta^t$ is then fixed.  For a datum $z=(y,\bm x)$, put
$\overline{\bm x}=\mathsf C_{L_{\bm x}}^{(\infty)}(\bm x)$ and define
its clipped, response-truncated gradient contribution by
\begin{align*}
    \bm h_{\bbeta^t}(z)
    =\{\psi_{L_\psi}'(\overline{\bm x}^\top\bbeta^t)-\Pi_R(y)\}
    \overline{\bm x}.
\end{align*}
The clipping map is defined for every \(\bm x\in\mathbb R^d\) and
satisfies \(\|\overline{\bm x}\|_\infty\leq L_{\bm x}\).
Consequently, \eqref{eq: saturated score bound} and
$|\Pi_R(y)|\leq R$ give the deterministic bound
\begin{align*}
    \|\bm h_{\bbeta^t}(z)\|_\infty
    \leq(R+S_\psi)L_{\bm x}.
\end{align*}
If $\bm Z$ and $\bm Z'$ differ in one datum, the two empirical
gradients therefore satisfy
\begin{align*}
\|\nabla\overline\L_{n,L_\psi}^R(\bbeta^t;\bm Z)
-\nabla\overline\L_{n,L_\psi}^R(\bbeta^t;\bm Z')\|_\infty
&\leq\frac{2(R+S_\psi)L_{\bm x}}n
 \leq\frac{B_0}{n}.
\end{align*}
Consequently, the pre-thresholding vector has sensitivity
\begin{align*}
    \left\|
    \{\bbeta^t-\eta^0\nabla\overline\L_{n,L_\psi}^R
    (\bbeta^t;\bm Z)\}
    -\{\bbeta^t-\eta^0\nabla\overline\L_{n,L_\psi}^R
    (\bbeta^t;\bm Z')\}
    \right\|_\infty
    \leq \frac{\eta^0B_0}{n}.
\end{align*}
Lemma \ref{lm: noisy hard thresholding privacy}, used with privacy
parameters $(\varepsilon/T,\delta/T)$, makes the conditional release
$\bbeta^{t+1}$ $(\varepsilon/T,\delta/T)$-differentially private.
This reasoning is valid for every realization of the preceding
releases, so adaptive basic composition over the $T$ iterations gives
$(\varepsilon,\delta)$-differential privacy.  The final projection is
post-processing and does not change the guarantee.

\textit{Accuracy.}
Fix $\xi=n^{-4}$.  We first collect simultaneous bounds for the score,
design clipping, predictor saturation, response truncation, and privacy
noise.

We begin with the design.  Assumption (D), applied to the coordinate
vector \(\bm e_j\), implies
\begin{align*}
\Pr\{|x_{ij}|>u\}
\leq2\exp\{-cu^2/K_{\bm x}^2\}
\qquad(i\in[n],\ j\in[d]).
\end{align*}
Set
\begin{align*}
\mathcal E_{\bm x}
=\left\{\max_{i\in[n],j\in[d]}|x_{ij}|\leq L_{\bm x}\right\}.
\end{align*}
Because \(L_{\bm x}=C_{\bm x}K_{\bm x}\sqrt{\log(ed)}\), a union
bound gives
\begin{align*}
\Pr(\mathcal E_{\bm x}^c)
&\leq2nd\exp\{-cC_{\bm x}^2\log(ed)\}.
\end{align*}
Choose \(C_{\bm x}\) so that \(cC_{\bm x}^2\geq8\).  Since
\(d\geq n\geq2\), the last display is at most
\(2nd(ed)^{-8}\leq n^{-4}=\xi\), after increasing
\(C_{\bm x}\) once more if necessary.  We have therefore proved
\begin{align}
\Pr(\mathcal E_{\bm x})\geq1-\xi.                         \label{eq: repaired sparse no covariate clipping}
\end{align}
On \(\mathcal E_{\bm x}\), every coordinatewise clipping operation is
the identity.  In particular,
\begin{align}
\overline\L_{n,L_\psi}^R(\bbeta)=\L_{n,L_\psi}^R(\bbeta)
\quad\text{for every }\bbeta\in\mathbb R^d.              \label{eq: sparse clipped loss coupling}
\end{align}

For the score, write
\begin{align*}
    \nabla\L_n(\bbeta^*)
    =\frac1n\sum_{i=1}^n
    \{\psi'(\bm x_i^\top\bbeta^*)-y_i\}\bm x_i.
\end{align*}
For completeness, we derive its concentration without conditioning on
a bounded-design event.  Put
\begin{align*}
r_i=y_i-\psi'(\bm x_i^\top\bbeta^*),
\qquad Z_{ij}=-r_ix_{ij}.
\end{align*}
The conditional mean of \(r_i\) is zero.  Replacing \(t\) by
\(c(\sigma)\lambda\) in \eqref{eq: sparse glm residual mgf} shows that,
conditionally on \(\bm x_i\),
\begin{align*}
\E\{\exp(\lambda r_i)\mid\bm x_i\}
\leq\exp\{c_2c(\sigma)\lambda^2/2\}.
\end{align*}
The Chernoff bound applied to \(r_i\) and \(-r_i\) therefore gives
\begin{align}
\Pr\{|r_i|>u\mid\bm x_i\}
\leq2\exp\{-u^2/(2c_2c(\sigma))\}.                       \label{eq: sparse conditional residual tail}
\end{align}
For every \(p\geq1\), integrate this tail using
\(\E|W|^p=\int_0^\infty p u^{p-1}\Pr(|W|>u)\,\d u\).
The change of variables
\(v=u^2/\{2c_2c(\sigma)\}\) yields
\begin{align}
\{\E(|r_i|^p\mid\bm x_i)\}^{1/p}
&\leq C\sqrt{c(\sigma)}\{\Gamma(p/2+1)\}^{1/p}\\
&\leq C\sqrt{c(\sigma)p}.                                \label{eq: sparse conditional residual moments}
\end{align}
Here the final inequality follows from the standard bound
\(\Gamma(p/2+1)^{1/p}\leq C\sqrt p\).

Assumption (D) similarly implies
\begin{align*}
\{\E|x_{ij}|^p\}^{1/p}\leq CK_{\bm x}\sqrt p.
\end{align*}
Using the tower property and \eqref{eq: sparse conditional residual moments},
we consequently obtain
\begin{align}
\{\E|Z_{ij}|^p\}^{1/p}
&=\left[\E\left\{|x_{ij}|^p
\E(|r_i|^p\mid\bm x_i)\right\}\right]^{1/p}\\
&\leq C\sqrt{c(\sigma)p}\{\E|x_{ij}|^p\}^{1/p}\\
&\leq CK_{\bm x}\sqrt{c(\sigma)}\,p.                    \label{eq: sparse score product moments}
\end{align}
Thus \(Z_{ij}\) is mean zero and sub-exponential with
\(\psi_1\)-norm at most \(CK_{\bm x}\sqrt{c(\sigma)}\).
Independence across observations and Bernstein's inequality now give
\begin{align*}
\Pr\{|[\nabla\L_n(\bbeta^*)]_j|>u\}
\leq2\exp\left[-cn\min\left\{
\frac{u^2}{K_{\bm x}^2c(\sigma)},
\frac{u}{K_{\bm x}\sqrt{c(\sigma)}}\right\}\right].
\end{align*}
The theorem's sample-size condition implies
\(n\geq C\log(2d/\xi)\): indeed, \(d\geq n\),
\(\xi=n^{-4}\), and all the remaining factors in that condition are
bounded below by fixed positive constants.  Hence, for a sufficiently
large fixed \(C\), the choice
\begin{align*}
u=CK_{\bm x}\sqrt{c(\sigma)\log(2d/\xi)/n}
\end{align*}
lies in the quadratic regime of Bernstein's inequality and makes the
last probability at most \(\xi/d\).  A union bound over \(j\in[d]\)
therefore gives
\begin{align}
    \|\nabla\L_n(\bbeta^*)\|_\infty^2
    \leq Cc(\sigma)\frac{\log(d/\xi)}{n}                    \label{eq: repaired sparse score concentration}
\end{align}
with probability at least $1-\xi$.  The fixed factor
\(K_{\bm x}^2\) has been absorbed into \(C\).

For predictor saturation and response truncation, define
\[
\mathcal E_{\rm pred}
=\{|\bm x_i^\top\bbeta^*|\leq L_\psi,\ i\in[n]\}.
\]
Since $\|\bbeta^*\|_2\leq1$, assumption (D) and a union bound give
\begin{align}
\Pr(\mathcal E_{\rm pred}^c)
&\leq2n\exp\{-cL_\psi^2/K_{\bm x}^2\}
\leq\xi,                                                   \label{eq: repaired sparse no score saturation}
\end{align}
after choosing $C_\psi$ sufficiently large. On
$\mathcal E_{\rm pred}$, equation \eqref{eq: saturated link score}
and the mean-value theorem give, for every $i$,
\[
\psi_{L_\psi}'(\bm x_i^\top\bbeta^*)
=\psi'(\bm x_i^\top\bbeta^*),
\qquad
|\psi'(\bm x_i^\top\bbeta^*)|
\leq|\psi'(0)|+c_2L_\psi.
\]
The two-sided Chernoff bound following from
\eqref{eq: sparse glm residual mgf} gives
\begin{align*}
\Pr\{|y_i-\psi'(\bm x_i^\top\bbeta^*)|>u\}
\leq2\exp\{-cu^2/c(\sigma)\}.
\end{align*}
On $\mathcal E_{\rm pred}$, if $|y_i|>R$, then the residual on the
left is larger than
$R-|\psi'(0)|-c_2L_\psi=C_R\sqrt{c(\sigma)\log n}$.
Thus a union bound over $i\in[n]$ and a sufficiently large $C_R$ give
\begin{align}
    \Pi_R(y_i)=y_i\quad\text{for every }i\in[n]              \label{eq: repaired sparse no truncation}
\end{align}
outside an event of probability at most $\xi$, on
$\mathcal E_{\rm pred}$. Consequently, on $\mathcal E_{\rm pred}$ and
the event in \eqref{eq: repaired sparse no truncation},
\begin{align}
\nabla\L_{n,L_\psi}^R(\bbeta^*)
=\nabla\L_n(\bbeta^*).                                    \label{eq: sparse saturated score coupling at truth}
\end{align}

It remains to control the privacy noise simultaneously.  At every
iteration, each coordinate of each Laplace vector in NoisyHT has scale
\begin{align*}
    b
    =\frac{2\eta^0B_0T}{n\varepsilon}
    \sqrt{3s\log(T/\delta)}.
\end{align*}
Indeed, the sensitivity is $\eta^0B_0/n$, and Algorithm
\ref{algo: private sparse glm} calls NoisyHT with privacy parameters
$(\varepsilon/T,\delta/T)$.  Across all iterations there are
\begin{align*}
    U=d(s+1)T
\end{align*}
Laplace coordinates: $sd$ selection-noise coordinates and $d$ final
perturbation coordinates per iteration.  Since
$\Pr\{|W|>bu\}=e^{-u}$ for $W\sim\operatorname{Laplace}(b)$, a union
bound gives, with probability at least $1-\xi$,
\begin{align*}
    \max_{\text{all }U\text{ coordinates}}|W|
    \leq b\log(U/\xi).
\end{align*}
For each $t$, the definition of $\mathcal W_t$ contains $s$ squared
$\ell_\infty$ norms and one squared $\ell_2$ norm supported on $s$
coordinates.  Hence, on this event,
\begin{align*}
    \mathcal W_t^2
    \leq s b^2\log^2(U/\xi)+s b^2\log^2(U/\xi)
    =2s b^2\log^2(U/\xi).
\end{align*}
Substituting the value of $b$ proves
\begin{align}
    \max_{t<T}\mathcal W_t^2
    &\leq
    C s b^2
    \log^2\{2d(s+1)T/\xi\}\notag\\
    &\leq
    C\frac{B_0^2T^2s^2\log(T/\delta)
    \log^2\{2d(s+1)T/\xi\}}
    {n^2\varepsilon^2}.                                    \label{eq: repaired sparse noise concentration}
\end{align}

	Intersect the five events
	\eqref{eq: repaired sparse no covariate clipping},
	\eqref{eq: repaired sparse score concentration},
	$\mathcal E_{\rm pred}$,
	\eqref{eq: repaired sparse no truncation}, and
	\eqref{eq: repaired sparse noise concentration} with
	$\mathcal E_{\rm sc}(\alpha,\gamma;L_\psi)$. Their intersection, denoted by
	$\mathcal E$, has probability at least $1-n^{-3}-5\xi$ by
	Proposition \ref{prop: sparse glm empirical curvature}.  To apply that
	proposition, note that $s\asymp s^*$ and that the theorem's sample-size
	condition, together with $\delta\leq n^{-1}$, $\varepsilon\leq1$, and
	the fixed lower bound on $c(\sigma)$, implies
	$n\geq C(s+s^*)\log(ed)$.  On
	$\mathcal E$, define
\begin{align*}
a_n&=Cc(\sigma)\frac{s\log(d/\xi)}n,\\
b_n&=C\frac{B_0^2T^2s^2\log(T/\delta)
\log^2\{2d(s+1)T/\xi\}}{n^2\varepsilon^2}.
\end{align*}
By \eqref{eq: repaired sparse no covariate clipping} and
\eqref{eq: sparse clipped loss coupling}, the update in Algorithm
\ref{algo: private sparse glm} is exactly an update based on
\(\L_{n,L_\psi}^R\) throughout the entire trajectory. Thus Lemma
\ref{lm: repaired noisy iht contraction} applies.  By
\eqref{eq: sparse saturated score coupling at truth}, its score
satisfies
\(\nabla\L_{n,L_\psi}^R(\bbeta^*)=\nabla\L_n(\bbeta^*)\).
Equations \eqref{eq: repaired sparse score concentration} and
\eqref{eq: repaired sparse noise concentration} show that the additive
part of \eqref{eq: repaired noisy iht recursion} is at most
$a_n+b_n$.

We next check the size of these quantities.  Since
$T\asymp\log n$, $\xi=n^{-4}$, $d\geq n$, $s\leq d$, and
$\delta\leq n^{-1}$,
\begin{align}
\log(d/\xi)&\lesssim\log d,\notag\\
\log(T/\delta)&\lesssim\log(1/\delta),\notag\\
\log\{2d(s+1)T/\xi\}&\lesssim\log d.                     \label{eq: sparse logarithm simplifications}
\end{align}
Moreover, the choice of $R$ and the fixed upper and lower bounds on
$c(\sigma)$ give
\begin{align*}
B_0^2
&=16(R+S_\psi)^2L_{\bm x}^2
\lesssim c(\sigma)\log n\log(ed).
\end{align*}
Hence
\begin{align}
a_n&\lesssim c(\sigma)\frac{s\log d}{n},\notag\\
b_n&\lesssim c(\sigma)
\frac{s^2\log^2d\,\log(ed)\log(1/\delta)\log^3n}
{n^2\varepsilon^2}.                                      \label{eq: sparse additive terms simplified}
\end{align}
Because $s\asymp s^*$ with a constant depending only on the curvature
ratio, the square root of the second bound in
\eqref{eq: sparse additive terms simplified} is, up to a fixed
constant,
\begin{align*}
\sqrt{c(\sigma)}\,
\frac{s^*\log d\sqrt{\log(ed)}
\sqrt{\log(1/\delta)}\log^{3/2}n}{n\varepsilon}.
\end{align*}
The theorem's sample-size assumption makes this quantity, and hence
\(b_n\), arbitrarily small when its implicit constant is sufficiently
large.  The same condition also implies
\(n\geq Cc(\sigma)s^*\log d\), because
\(0<\varepsilon\leq1\), \(\delta\leq n^{-1}\), and
\(c(\sigma)\) is bounded above and away from zero.  It therefore makes
\(a_n\) arbitrarily small as well.  In particular, because
\(1-q>0\) is a fixed curvature-dependent constant, we may ensure
\begin{align}
    a_n+b_n\leq4(1-q).                                    \label{eq: sparse additive terms invariant region}
\end{align}

Starting from
$\|\bbeta^0-\bbeta^*\|_2^2=\|\bbeta^*\|_2^2\leq1$, we now prove by
induction that $\|\bbeta^t-\bbeta^*\|_2^2\leq4$ for every
$t\leq T$.  If the claim holds at time $t$, then the hypothesis
$\|\bbeta^t-\bbeta^*\|_2\leq2$ of Lemma
\ref{lm: repaired noisy iht contraction} is satisfied, and
\eqref{eq: repaired noisy iht recursion} together with
\eqref{eq: sparse additive terms invariant region} gives
\begin{align*}
\|\bbeta^{t+1}-\bbeta^*\|_2^2
\leq4q+4(1-q)=4.
\end{align*}
This closes the induction and justifies applying the contraction lemma
at every iteration.

Iterating \eqref{eq: repaired noisy iht recursion} and summing the
geometric series gives
\begin{align}
    \|\bbeta^T-\bbeta^*\|_2^2
    &\leq q^T\|\bbeta^0-\bbeta^*\|_2^2\notag\\
    &\quad+C\sum_{t=0}^{T-1}q^{T-1-t}
    \{s\|\nabla\L_n(\bbeta^*)\|_\infty^2+\mathcal W_t^2\}
    \notag\\
    &\leq q^T\|\bbeta^0-\bbeta^*\|_2^2
    +\frac{C}{1-q}\left\{s\|\nabla\L_n(\bbeta^*)\|_\infty^2
    +\max_{t<T}\mathcal W_t^2\right\}.                    \label{eq: repaired sparse iterated recursion}
\end{align}
\furtherrev{With $q=q_0=1-3\alpha/(4\gamma)$ as above, the theorem's choice
\[
T=\left\lceil\frac{\log n}{-\log q_0}\right\rceil
\]
satisfies $q^T=q_0^T\leq n^{-1}$.}
Substitution of \eqref{eq: repaired sparse score concentration} and
\eqref{eq: repaired sparse noise concentration} into
\eqref{eq: repaired sparse iterated recursion}, followed by
\eqref{eq: sparse logarithm simplifications}, gives
\begin{align*}
\|\bbeta^T-\bbeta^*\|_2^2
\lesssim c(\sigma)\left\{
\frac{s^*\log d}{n}
+\frac{(s^*\log d)^2\log(ed)\log(1/\delta)\log^3n}
{n^2\varepsilon^2}\right\}.
\end{align*}
Here $q^T\leq n^{-1}$ is absorbed into the first term, and
$s\asymp s^*$.  Since $\bbeta^*$ belongs to the unit Euclidean ball,
the Euclidean projection theorem gives
\begin{align*}
\|\widehat\bbeta-\bbeta^*\|_2
\leq\|\bbeta^T-\bbeta^*\|_2.
\end{align*}
	This proves \eqref{eq: glm upper bound}.  Since $\xi=n^{-4}$, the
	event $\mathcal E$ has probability at least $1-Cn^{-3}$.

Finally, both $\widehat\bbeta$ and $\bbeta^*$ belong to the unit ball, so
$\|\widehat\bbeta-\bbeta^*\|_2^2\leq4$ on every outcome.  Let $r_n$
denote the right side of \eqref{eq: glm upper bound} without its
comparison constant.  Splitting the expectation over $\mathcal E$ and
$\mathcal E^c$ gives
\begin{align*}
\E_{\bbeta^*}\|\widehat\bbeta-\bbeta^*\|_2^2
&=\E\{\|\widehat\bbeta-\bbeta^*\|_2^2\1_{\mathcal E}\}
+\E\{\|\widehat\bbeta-\bbeta^*\|_2^2\1_{\mathcal E^c}\}\\
	&\leq Cr_n+4\Pr(\mathcal E^c)
	\leq Cr_n+Cn^{-3}.
	\end{align*}
	Here $d\geq n$, $s^*\geq1$, and the first term in $r_n$ is at least
	a constant multiple of $\log n/n$, so the last remainder is absorbed
	by $Cr_n$.  This proves \eqref{eq: sparse glm expected upper bound}.
\end{proof}

\endgroup

	\section{Omitted Proofs in Section \ref{sec: nonparametric}}\label{sec: nonparametric proofs}

\subsection{Proof of Proposition \ref{prop: nonparametric attack soundness}}\label{sec: proof of prop: nonparametric attack soundness}
\begin{proof}[Proof of Proposition \ref{prop: nonparametric attack soundness}]
	Let $A'_{i} := \mathcal A(M(\bm X_i', \bm Y_i'), (X_i, Y_i))$, where $(\bm X'_i, \bm Y'_i)$ is an adjacent data set of $(\bm X, \bm Y)$ obtained by replacing $(X_i, Y_i)$ with an independent copy.
	
		For each $A_{i}$ and every $T > 0$, equation \eqref{eq: score attack upper bound intermediate step} gives
		\begin{align*}
			\E A_{i} \leq \E A'_{i} + 2\varepsilon\E|A'_{i}| + 2 \delta T
			+ \int_T^\infty \{\Pro(|A_{i}| > t)+\Pro(|A_i'|>t)\} \d t.
		\end{align*}
	Now observe that, since $M(\bm X'_i, \bm Y'_{i})$ and $(X_i, Y_i)$ are independent by construction, we have
	\begin{align*}
			\E A'_{i} = \left\langle \E \left(M(\bm X'_i, \bm Y'_{i}) - \bth\right),
			\sigma^{-2}\E\left[\left(Y_i - \sum_{j=1}^k \theta_j\varphi_j(X_i)\right)\bm\varphi(X_i)\right] \right\rangle = 0.
		\end{align*}
		Indeed, the residual in the last expectation is $\xi_i$, which is centered and independent of $X_i$. For $\E|A'_i|$, orthonormality gives $\E\{\bm\varphi(X_i) \bm\varphi(X_i)^\top\} = \bm I_k$, and conditional Cauchy--Schwarz yields
	\begin{align*}
		\E |A'_i| \leq \E \left|\left\langle M(\bm X'_i, \bm Y'_i) - \bth,  \sigma^{-2}\xi_i\bm\varphi(X_i)\right\rangle\right| \leq \sigma^{-1} \sqrt{\E_{\bm X, \bm Y|\bth} \|M(\bm X, \bm Y) - \bth\|_2^2}.
	\end{align*} 
		{For the tail terms, the range condition on $M$ and $\tau_j\geq1$ imply that both $M$ and the true parameter have Euclidean norm at most $C$. Hence their distance is at most $2C$. Also $\|\bm \varphi(X_i)\|_2 \leq \sqrt{2k}$. Consequently, both $|A_i|$ and $|A_i'|$ are pointwise bounded by $2C\sqrt{2k}\sigma^{-1}|Z|$, where $Z\sim N(0,1)$. Taking $T=2C\sqrt{2k}\sigma^{-1}\sqrt{2\log(1/\delta)}$ and using the Gaussian tail integral gives}
		\begin{align*}
			\rev{\int_T^\infty \{\Pro(|A_{i}| > t)+\Pro(|A_i'|>t)\} \d t
			\lesssim C\sqrt{k}\sigma^{-1}\delta.}
		\end{align*}
			Combining these bounds gives
		\begin{align*}
		\E A_i \leq \sigma^{-1}\left(2\varepsilon\sqrt{\E_{\bm X, \bm Y|\bth} \|M(\bm X, \bm Y) - \bth\|_2^2} + C_1C\sqrt{k}\delta\sqrt{\log(1/\delta)}\right).
		\end{align*}
	Summing over $i \in [n]$ completes the proof.
\end{proof}

\subsection{Proof of Proposition \ref{prop: nonparametric attack completeness}}\label{sec: proof of prop: nonparametric attack completeness}
\begin{proof}[Proof of Proposition \ref{prop: nonparametric attack completeness}]
			{Let $g(\bth)=\E_{\bm X,\bm Y\mid\bth}M(\bm X,\bm Y)$. Since $M$ takes values in $\Theta_k(\alpha,C)$, $g$ is bounded. The Gaussian regression likelihood has parameter-independent support; on a neighborhood of any $\bth$, its derivative multiplied by the bounded output of $M$ is dominated by an integrable Gaussian envelope. Thus the completeness identity applies. The density $q_B$ vanishes quadratically at $\pm B$, so $g_j(t,\bth_{-j})q_B(t)\to0$ at both endpoints. The score identity and Proposition \ref{prop: score attack stein's lemma} therefore give
	\begin{align*}
		\sum_{i\in[n]}\E_\bth\E_{\bm X,\bm Y\mid\bth}A_i
		&=\sum_{j=1}^k\E_\bth\frac{\partial g_j(\bth)}{\partial\theta_j}\\
		&=k+\sum_{j=1}^k\E_\bth\left[(\theta_j-g_j(\bth))\frac{q_B'(\theta_j)}{q_B(\theta_j)}\right].
	\end{align*}
	A direct calculation gives
	\[
	\E\left[-\theta_j\frac{q_B'(\theta_j)}{q_B(\theta_j)}\right]=1,
	\qquad
	\E\left[\left(\frac{q_B'(\theta_j)}{q_B(\theta_j)}\right)^2\right]=\frac{10}{B^2}.
	\]
	Hence Cauchy--Schwarz and the assumed risk bound yield
	\begin{align*}
		\sum_{i\in[n]}\E_\bth\E_{\bm X,\bm Y\mid\bth}A_i
		&\geq k-\sqrt{\E_\bth\E_{\bm X,\bm Y\mid\bth}\|M-\bth\|_2^2}\sqrt{\frac{10k}{B^2}}\\
		&\geq \frac{k}{2}.
	\end{align*}}
\end{proof}

\subsection{Proof of Proposition \ref{prop: nonparametric coefficient lower bound}}\label{sec: proof of prop: nonparametric coefficient lower bound}
\begin{proof}[Proof of Proposition \ref{prop: nonparametric coefficient lower bound}]
		Let $M$ be any $(\varepsilon,\delta)$-differentially private estimator \furtherrev{with values in the coefficient decision space $\mathbb R^k$}, and put
		\[
		\widetilde M=\Pi_{\Theta_k(\alpha,C)}\circ M.
		\]
		Projection is post-processing and, for every $\bth\in\Theta_k(\alpha,C)$,
		\[
		\|\widetilde M-\bth\|_2\leq\|M-\bth\|_2.
		\]
		It therefore suffices to prove the lower bound for the $\Theta_k(\alpha,C)$-valued estimator $\widetilde M$, which we relabel as $M$. \furtherrev{This projection is used only in the lower-bound reduction; the original function-estimation decision space remains $L^2[0,1]$.} Define
		\[
		\mathcal R(M)=\sup_{\bth\in\Theta_k(\alpha,C)}
		\E_{\bm X,\bm Y\mid\bth}\|M-\bth\|_2^2.
		\]

		If $\mathcal R(M)>kB^2/40$, then
		\[
		\mathcal R(M)\gtrsim kB^2\asymp k^{-2\alpha},
		\]
		which gives the first branch of the claimed minimum. It remains to consider
		$\mathcal R(M)\leq kB^2/40$, in which case Proposition
		\ref{prop: nonparametric attack completeness} gives
		\begin{align}\label{eq: nonparametric completeness half k}
		\overline A:=\sum_{i=1}^n\E_{\bth\sim\bm\pi}
		\E_{\bm X,\bm Y\mid\bth}A_i\geq k/2.
		\end{align}
		Averaging Proposition \ref{prop: nonparametric attack soundness} over the same prior and applying Jensen's inequality yields
		\begin{align}\label{eq: nonparametric averaged soundness}
		\overline A
		&\leq \sigma^{-1}\left\{2n\varepsilon\sqrt{\overline R}
		+C_1Cn\sqrt{k\log(1/\delta)}\,\delta\right\},\\
		\overline R
		&:=\E_{\bth\sim\bm\pi}\E_{\bm X,\bm Y\mid\bth}
		\|M-\bth\|_2^2\leq\mathcal R(M).\notag
		\end{align}
		For $0<\delta\leq cn^{-2}$, the function
		$x\mapsto x\sqrt{\log(1/x)}$ is increasing for all sufficiently small $x$, and therefore
		\begin{align*}
		n\delta\sqrt{k\log(1/\delta)}
		&\leq cn^{-1}\sqrt{k\log(n^2/c)}\\
		&\leq cn^{-1}\sqrt{\log(n^2/c)}\,k.
		\end{align*}
			Because $\sigma\geq\sigma_->0$, choosing $c$ sufficiently small makes the last display at most, uniformly for $n\geq2$,
		$\sigma k/(4C_1C)$. Combining \eqref{eq: nonparametric completeness half k} and \eqref{eq: nonparametric averaged soundness} gives
		\[
		\frac{\sigma k}{4}\leq2n\varepsilon\sqrt{\overline R},
		\qquad\text{and hence}\qquad
		\overline R\geq\frac{\sigma^2k^2}{64n^2\varepsilon^2}.
		\]
		Under the fixed-noise convention $\sigma\geq\sigma_->0$, this implies
		$\mathcal R(M)\gtrsim k^2/(n^2\varepsilon^2)$. Combining the two risk cases and then taking the infimum over private estimators proves the proposition.
\end{proof}

\subsection{Proof of Theorem \ref{thm: nonparametric upper bound}}\label{sec: proof of thm: nonparametric upper bound}
\begin{proof}[Proof of Theorem \ref{thm: nonparametric upper bound}]
	The system $\{\varphi_j\}_{j \in \mathbb N}$ is an orthonormal basis of $L^2[0,1]$, and therefore
	\begin{align}\label{eq: nonparametric estimator integrated risk decomposition}
			\int_0^1 (\tilde f_{K, T}(x) - f(x))^2 \d x = \|\tilde\bth_{K, T} - \bth_K\|_2^2 + \sum_{j > K} \theta_j^2,
	\end{align}
	where $\bth_K = (\theta_1, \theta_2, \cdots, \theta_K)$ is the vector of the first $K$ Fourier coefficients of $f$. Let $\hat\bth_K$ denote the vector of the first $K$ empirical Fourier coefficients, \rev{$\hat\bth_K = (\hat\theta_1, \hat\theta_2, \cdots, \hat\theta_K)$}, and $\hat\bth_{K, T}$ denote the noiseless version of $\tilde\bth_{K, T}$,
	\begin{align*}
		\hat\bth_{K, T} = \frac{1}{n} \sum_{i=1}^n Y_i\1(|Y_i| \leq T) \cdot \bm \varphi(X_i).
	\end{align*}
	We have
	\begin{align}
		\E\|\tilde\bth_{K, T} - \bth_K\|_2^2 &\lesssim \E \|\hat\bth_{K, T} - \bth_K\|_2^2 + \E \|\bm w \|_2^2 \notag\\
		&\lesssim \E \|\hat\bth_{K, T} - \hat\bth_K\|_2^2 + \E \|\hat\bth_K - \bth_K\|_2^2 + \E \|\bm w \|_2^2. \label{eq: nonparametric estimator coefficient risk decomposition}
	\end{align}
	{For the first term, put $V_i=Y_i\1(|Y_i|>T)\bm\varphi(X_i)$. Since $\sup_x\|\bm\varphi(x)\|_2^2\leq2K$,
	\begin{align*}
		\E\left\|\frac1n\sum_{i=1}^nV_i\right\|_2^2
		&=\frac1n\E\|V_1-\E V_1\|_2^2+\|\E V_1\|_2^2\\
		&\leq\frac{2K}{n}\E\!\left[Y_1^2\1(|Y_1|>T)\right]+2K\left(\E\!\left[|Y_1|\1(|Y_1|>T)\right]\right)^2.
	\end{align*}
		The inhomogeneous Sobolev embedding gives $\sup_{f\in\tilde W(\alpha,C)}\|f\|_\infty\leq r_{\alpha,C}$. With $T=r_{\alpha,C}+2\sigma\sqrt{2\log n}$, the event $|Y_1|>T$ implies $|Z|>a_n:=2\sqrt{2\log n}$ for $Z=\xi_1/\sigma$. The standard Gaussian identities
		\[
		\E\{|Z|\1(|Z|>a)\}=2\phi(a),\qquad
		\E\{Z^2\1(|Z|>a)\}=2\{a\phi(a)+\Phi(-a)\}
		\]
		show that the first two truncated moments of $|Y_1|$ are bounded by a fixed polynomial in $\log n$ times $n^{-4}$. Since the stated choice of $K$ satisfies $K\leq n$, the preceding display therefore implies, uniformly over $f\in\tilde W(\alpha,C)$, that
	\begin{align*}
		\E \|\hat\bth_{K, T} - \hat\bth_K\|_2^2\lesssim n^{-2}.
	\end{align*}}
		For the untruncated empirical coefficients, independence across observations and unbiasedness give
		\begin{align*}
		\E\|\hat\bth_K-\bth_K\|_2^2
		&=\frac1n\sum_{j=1}^K\Var\{Y_1\varphi_j(X_1)\}\\
		&\leq\frac1n\E\{Y_1^2\|\bm\varphi(X_1)\|_2^2\}
		\leq\frac{2K}{n}\E Y_1^2
		\lesssim\frac Kn.
		\end{align*}
			The last inequality is uniform over $f\in\widetilde W(\alpha,C)$ because $\|f\|_\infty\leq r_{\alpha,C}$ and $\sigma$ is a fixed model constant. \rev{The coordinates of $\bm w$ are independent $N(0,\tau^2)$ variables, so its second moment is the exact quantity
	\begin{align}\label{eq: nonparametric gaussian noise second moment}
	\E\|\bm w\|_2^2=K\tau^2
	=\frac{16T^2K^2\log(1.25/\delta)}{n^2\varepsilon^2}.
	\end{align}
		} Finally, to bound the right side of \eqref{eq: nonparametric estimator integrated risk decomposition}, by the definition of Sobolev ellipsoid \eqref{eq: Sobolev ellipsoid definition} we have
	\begin{align*}
		\rev{\sum_{j > K} \theta_j^2 \leq \left(\inf_{j>K}\tau_j\right)^{-2}\sum_{j > K} \tau_j^2\theta_j^2 \lesssim K^{-2\alpha}.}
	\end{align*}
		Since $T^2\lesssim\log n$ and $\log(1.25/\delta)\lesssim\log(1/\delta)$ for $0<\delta<1/2$, the preceding bounds and \eqref{eq: nonparametric estimator coefficient risk decomposition} give
	\begin{align*}
		\E \left[\int_0^1 (\tilde f_{K, T}(x) - f(x))^2 \d x \right] \lesssim \frac{K}{n} + \rev{\frac{K^2\log n\log(1/\delta)}{n^2\varepsilon^2}} + K^{-2\alpha}.
	\end{align*}
		\rev{Put $K_0=\min\{n^{1/(2\alpha+1)},(n\varepsilon)^{1/(\alpha+1)}\}$. If $c_1K_0\geq2$, then the stated integer choice satisfies $K\asymp K_0$, while always $K\leq K_0$. Hence
		\begin{align*}
		\frac Kn&\leq n^{-2\alpha/(2\alpha+1)},\\
		K^{-2\alpha}&\lesssim n^{-2\alpha/(2\alpha+1)}
		+(n\varepsilon)^{-2\alpha/(\alpha+1)},\\
		\frac{K^2}{n^2\varepsilon^2}
		&\leq(n\varepsilon)^{-2\alpha/(\alpha+1)}.
		\end{align*}
	If $c_1K_0<2$, then at least one of the two target rates is bounded below by a fixed positive constant. The preceding bound with $K=1$ is $O\{1+\log n\log(1/\delta)/(n^2\varepsilon^2)\}$, and its privacy term is bounded by a constant multiple of the privacy rate displayed in the theorem. Enlarging the comparison constant therefore covers this remaining regime. These observations prove \eqref{eq: nonparametric upper bound}.}
\end{proof}

\end{document}